\documentclass[twoside,11pt]{article} 
\usepackage[margin=1.0in]{geometry}

\usepackage[square,numbers]{natbib}
\RequirePackage[colorlinks,citecolor=blue,urlcolor=blue]{hyperref}
\usepackage{mathrsfs,amsthm,amsmath}
\usepackage{amssymb,multirow}
\usepackage{amsfonts}
\usepackage{stmaryrd}
\usepackage{dsfont}

\usepackage{ucs}
\usepackage[utf8x]{inputenc}
\usepackage{amsmath}
\usepackage{amsfonts}
\usepackage{authblk}
\usepackage{amssymb}
\usepackage{amsthm}
\usepackage{fontenc}
\usepackage{appendix}
\usepackage{graphicx}

\usepackage{bbm}
\usepackage{bm}
\RequirePackage[OT1]{fontenc}
\usepackage{tikz}
\usepackage[misc,geometry]{ifsym}

\let\svtikzpicture\tikzpicture
\def\tikzpicture{\noindent\svtikzpicture}
\usetikzlibrary{shapes,arrows}

\definecolor {bianco}{RGB}{255,255,255}

\tikzstyle{block} = [rectangle, fill=bianco,
    text width=6em, text centered, rounded corners, minimum height=4em]
\tikzstyle{line} = [draw, -latex']

\def\rd{\mathbb{R}^d}

\def\naturals{\mathbb{N}}

\def\reals{\mathbb{R}}

\def\ps{\Theta}
\def\reals{\mathbb{R}}
\def\probabilityspace{(\Omega, \Fcr, \P)}
\def\ss{\mathbb{X}}

\def\ssa{\mathscr{X}}
\def\psa{\mathscr{T}}

\allowdisplaybreaks

\newcommand{\Space}{\mathbb{M}}

\newcommand{\Hunoz}{\mathrm{H}^1_{\ast}(0,1)}
\newcommand{\Huno}{\mathrm{H}^1(0,1)}

\newcommand{\T}{\mathbb{T}}
\newcommand{\X}{\mathbb{X}}

\newcommand{\Sd}{\mathbb{S}}

\newcommand{\Hbb}{\mathbb{H}}
\newcommand{\B}{\mathbb{B}}
\newcommand{\D}{\mathrm{D}}
\newcommand{\Kbb}{\mathbb{K}}
\newcommand{\Vbb}{\mathbb{V}}
\newcommand{\Dcal}{\mathcal{D}}

\newcommand{\R}{\mathbb{R}}
\newcommand{\E}{\mathds{E}}

\newcommand{\N}{\mathbb{N}}
\renewcommand{\P}{\mathds{P}}

\newcommand{\Bcr}{\mathscr{B}}
\newcommand{\Fcr}{\mathscr{F}}

\newcommand{\Wdue}{\mathcal{W}_2}
\newcommand{\Wp}{\mathcal{W}_p}

\def\rd{\mathbb{R}^d}
\def\naturals{\mathbb{N}}
\def\reals{\mathbb{R}}

\newcommand{\kksf}{\mathsf{K}}

\newcommand{\pms}{\mathcal{P}(\Theta)}

\newcommand{\empiric}{\mathfrak{e}_n}
\newcommand{\pfrak}{\mathfrak{p}}
\newcommand{\ee}{\textsf{E}}

\newcommand{\ud}{\mathrm{d}}

\newcommand{\ind}{\mathds{1}}
\newcommand{\Cpw}{\mathfrak{C}_2}
\newcommand{\Bdual}{\mathbb{B}^{\ast}}

\newtheorem{thm}{Theorem}[section]

\newtheorem{defi}[thm]{Definition}

\newtheorem{prp}[thm]{Proposition}
\newtheorem{lm}[thm]{Lemma}
\newtheorem{cor}[thm]{Corollary}
\newtheorem{remark}[thm]{Remark}

\newtheorem{assum}[thm]{Assumptions}

\providecommand{\keywords}[1]
{
  \small	
  \textbf{\textit{Keywords:}} #1
}

\begin{document}

\title{Strong posterior contraction rates via Wasserstein dynamics}


\author[1]{Emanuele Dolera\thanks{emanuele.dolera@unipv.it}}
\author[2]{Stefano Favaro\thanks{stefano.favaro@unito.it}}
\author[3]{Edoardo Mainini\thanks{mainini@dime.unige.it}}
\affil[1]{\small{Department of Mathematics, University of Pavia, Italy}}
\affil[2]{\small{Department of Economics and Statistics, University of Torino and Collegio Carlo Alberto, Italy}}
\affil[3]{\small{Department of Mechanical Engineering, University of Genova, Italy}}

\maketitle

\begin{abstract}
In Bayesian statistics, posterior contraction rates (PCRs) quantify the speed at which the posterior distribution concentrates on arbitrarily small neighborhoods of a true model, in a suitable way, as the sample size goes to infinity. In this paper, we develop a new approach to PCRs, with respect to strong norm distances on parameter spaces of functions. Critical to our approach is the combination of a local Lipschitz-continuity for the posterior distribution with a dynamic formulation of the Wasserstein distance, which allows to set forth an interesting connection between PCRs and some classical problems arising in mathematical analysis, probability and statistics, e.g., Laplace methods for approximating integrals, Sanov's large deviation principles in the Wasserstein distance, rates of convergence of mean Glivenko-Cantelli theorems, and estimates of weighted Poincar\'e-Wirtinger constants. We first present a theorem on PCRs for a model in the regular infinite-dimensional exponential family, which exploits sufficient statistics of the model, and then extend such a theorem to a general dominated model. These results rely on the development of novel techniques to evaluate Laplace integrals and weighted Poincar\'e-Wirtinger constants in infinite-dimension, which are of independent interest. The proposed approach is applied to the regular parametric model, the multinomial model, the finite-dimensional and the infinite-dimensional logistic-Gaussian model and the infinite-dimensional linear regression. In general, our approach leads to optimal PCRs in finite-dimensional models, whereas for infinite-dimensional models it is shown explicitly how the prior distribution affect PCRs.
\end{abstract}



\section{Introduction}
Bayesian consistency guarantees that the posterior distribution concentrates on arbitrarily small neighborhoods of the true model, in a suitable way, as the sample size goes to infinity (\citet{Doo(49),Sch(65),Fre(63),Fre(65),Dia(86),Bar(99),GGV(00),W(04)}). See \citet[Chapter 6 and Chapter 7]{GV(00)} for a general overview on Bayesian consistency. Posterior contractions rates (PCRs) strengthen the notion of Bayesian consistency, as they quantify the speed at which such small neighborhoods of the true model may decrease to zero meanwhile still capturing most of the posterior mass. The problem of establishing optimal PCRs in finite-dimensional (parametric) Bayesian models have been first considered in \citet{IbHa(81)} and \citet{LeCam(86)}. However, it is in the works of \citet{GGV(00)} and \citet{ShWa(01)} that the problem of establishing PCRs have been investigated in a systematic way, setting forth a general approach to provide PCRs in both finite-dimensional and infinite-dimensional (nonparametric) Bayesian models. Since then, several methods have been proposed to obtain more explicit and also sharper PCRs. Among them, we recall the metric entropy approach, in combination with the definition of specific tests (\citet{Sch(65),GGV(00)}), the methods based on bracketing numbers and entropy integrals (\citet{ShWa(01)}), the martingale approach (\citet{W(04),WLP(07)}), the Hausdorff entropy approach \citet{Xing(10)}, and some approaches based on the Wasserstein distance (\citet{deBlasi(20),Cam(22)}). See \citet[Chapter 8 and Chapter 9]{GV(00)}, and references therein, for a comprehensive and up-to-date account on PCRs.

\subsection{Our contributions}

In this paper, we develop a new approach to PCRs, in the spirit of the seminal work of \citet{GGV(00)}. We consider a dominated statistical model as a family $\mathscr M = \{f_{\theta}\}_{\theta \in \Theta}$ of densities, with the parameter space $\Theta$ being a (possibly infinite-dimensional) separable Hilbert space. We focus on posterior Hilbert neighborhoods of a given true parameter, say $\theta_0$, measuring PCRs in terms of strong norm distances on parameter spaces of functions, such as Sobolev-like norms. This 
assumption on $\Theta$ yields a stronger metric structure on $\mathscr M$, as a subset of the space of densities, usually not equivalent to those considered so far by the literature on nonparametric density estimation (see, e.g., \citet{Gho(99),GiNi(11),Scri(06), ShWa(01), vdVvZ(08), W(04),WLP(07)}), based on the choice of (pseudo-)distances such as the $\mathrm L^p$-norm, the Hellinger, the Kullback-Leibler, and the chi-square. To the best of our knowledge, we are not aware of works in the Bayesian literature that deal with strong PCRs for density estimation by using constructive tests, as prescribed by the standard theory, even if this line of research could be pursued as well. As far as we know, the standard nonparametric approach covers the case of (semi-)metrics which are dominated
by the Hellinger distance (see, e.g., \citet[Proposition D.8]{GV(00)}). 

We present a theorem on PCRs for the regular infinite-dimensional exponential family of statistical model, and a theorem on PCRs for a general dominated statistical models. The former may be viewed as a special case of the latter, allowing to exploit sufficient statistics arising from the infinite-dimensional exponential family.  Critical to our approach is an assumption of local Lipschitz-continuity for the posterior distribution, with respect to the observations or a sufficient statistics of them. Such a property is typically known as ``Bayesian well-posedness'' (\citet[Section 4.2]{Stew(10)}), and it has been investigated in depth in  \citet{DM(20a),DM(20b)}. By combining the local Lipschitz-continuity with the dynamic formulation of the Wasserstein distance (\citet{BB(00), AGS(08)}), referred to as Wasserstein dynamics, we set forth a connection between the problem of establishing PCRs and some classical problems arising in mathematical analysis, probability and statistics, e.g., Laplace methods for approximating integrals (\citet{Bre(94),Won(01)}), Sanov's large deviation principle in Wasserstein distance (\citet{BGV(07),Jing(20)}), rates of convergence of mean Glivenko-Cantelli theorems (\citet{Ajt(84),AST(19),DY(95),DR(19),FG(15),Ta(94a),Ta(94b),BoLe(19),WB(19),Jing(20)}), and estimates of weighted Poincar\'e-Wirtinger constants (\citet[Chapter 4]{BGL(14)}, \citet[Chapter 15]{HKM(93)}). In particular, our study leads to introduce new results on Laplace methods for approximating integrals and the estimation of weighted Poincar\'e-Wirtinger constants in infinite dimension, which are of independent interest.

Some applications of our main theorems  are presented for the regular parametric model, the multinomial model, the finite-dimensional and the infinite-dimensional logistic-Gaussian model and the infinite-dimensional linear regression. It turns out that our main results lead to optimal PCRs in finite dimension, whereas in infinite dimension it is shown explicitly how the prior distribution affects PCRs. Among the applications of our results, the infinite-dimensional logistic-Gaussian model is arguably the best setting to motivate the use of strong norm distances. In such a setting our approach is of interest when the ultimate goal of the inferential procedure is the estimation of some functional $\Phi(f_{\theta})$ of the density \citep[Chapter 6]{Sil(82)} for which the mapping $f \mapsto \Phi(f)$ is not continuous with respect to the aforesaid metrics on densities, whereas $\theta \mapsto \Phi(f_{\theta})$ turns out to be even locally Lipschitz-continuous with respect to the Hilbertian metric on $\Theta$. Thus, strong norms allow to consider larger classes of functionals of density functions, and then possibly a broader range of analyses. Another motivation in the use of strong norms comes from the theory of density estimation under penalized loss functions, with penalizations depending on derivatives of the density, according to the original Good-Gaskins proposal \citep{GG(71), Silv(82)}. As these penalized loss functions are used to derive smoother estimators, it sounds interesting to derive relative PCRs under the same loss functions.  

\subsection{Related works}

The most popular classical (frequentist) approaches to density estimation are developed within the following frameworks: i) a parameter space that is the space of density functions, typically endowed with the $\mathrm{L}^p$ norm or the Hellinger distance (\citet{Tsy(94)}), usually associated to the notion of ``strong consistency''; ii) a parameter space that is the space of density functions endowed with the Wasserstein distance, under which the parameter space is metrized according to a (concrete) metric structure on the space of the observations (\citet{Bert(21)}), usually associated to the notion of ``weak consistency''. Both these 
frameworks are different from the one we consider in this paper, and therefore a comparison of our PCRs with optimal minimax rates from \citet{Tsy(94)} and \citet{Bert(21)} it is not directly possible. Within the classical literature, \citet{SFGHK(17)} considered our statistical framework and provided rates of consistency under the infinite-dimensional exponential family of statistical models, though without any formal statement on their minimax optimality. In principle, our approach to PCRs may be developed within the aforementioned popular statistical frameworks for density estimation. However, since our approach relies on properties of the Wasserstein distance that are well-known for parameter spaces with a linear structure, i.e. Wasserstein dynamics, the framework considered in this paper is the most natural and convenient to start with. As for the other statistical frameworks for density estimation, we conjecture that our approach to PCRs requires a suitable formulation of Wasserstein dynamics for parameter spaces with a nonlinear structure. While such a formulation is available from \citet{Gig(09)} and \citet{Gig(12)}, it is still not clear to us how to exploit it to deal with PCRs.

\subsection{Organization of the paper}

The paper is structured as follows. In Section \ref{sect:main_problem} we recall the definition of PCR, presenting an equivalent definition in terms of the Wasserstein distance, and we outline the main steps of our approach to PCRs. Section \ref{sec_mainres} contains the main results of our work, that is a theorem on PCRs for the regular infinite-dimensional exponential family of statistical models, and a generalization of it for general dominated statistical models. In Section \ref{sect:illustrations} we present some applications of our results for the regular parametric model, the multinomial model, the finite-dimensional and the infinite-dimensional logistic-Gaussian model and the infinite-dimensional linear regression. Section \ref{sec:discuss} contains a discussion of some directions for future work, especially with respect to the application of our approach to other nonparametric models, such as the popular class of hierarchical (mixture) models. Proofs of our results are deferred to appendices.


\section{A new approach to PCRs}\label{sect:main_problem}

We consider $n\geq1$ observations to be modeled as part of a sequence $X^{(\infty)} := \{X_i\}_{i \geq 1}$ of exchangeable random variables, with the $X_i$'s taking values in a measurable space $(\ss, \ssa)$. Let $(\ps, \ud_{\ps})$ be metric space, referred to as the parameter space, endowed with its Borel $\sigma$-algebra $\psa$. Moreover, let $\pi$  be a probability measure on $(\ps, \psa)$, referred to as the prior measure, and let $\mu(\cdot\,|\, \cdot) : \ssa\times\ps \rightarrow [0,1]$ be a probability kernel, referred to as the statistical model. The Bayesian approach relies on modeling the parameter of interest as a $\ps$-valued random variable, say $T$, with probability distribution $\pi$. At the core of Bayesian inferences lies the posterior distribution, that is the conditional distribution of $T$ given a random sample $(X_1, \dots, X_n)$, whenever both $T$ and the sequence $X^{(\infty)}$ are supported on a common probability space $\probabilityspace$. The minimal regularity conditions that are maintained, and possibly strengthened, throughout the paper are the following: the set $\ss$ is a separable topological space, with $\ssa$ coinciding with the ensuing Borel $\sigma$-algebra, and $(\ps, \psa)$ is a standard Borel space.  In this setting, the posterior distribution can be represented through a probability kernel $\pi_n(\cdot\,|\, \cdot): \psa\times\ss^n \rightarrow [0,1]$ that satisfies the disintegration
\begin{equation} \label{disintegration}
\P[X_1\in A_1, \dots, X_n\in A_n, T\in B] = \int_{A_1\times\dots\times A_n} \pi_n(B\,|\,x^{(n)}) \alpha_n(\ud x^{(n)})
\end{equation}
for all sets $A_1, \dots, A_n \in \ssa$ and $B \in \psa$ and $n\geq1$, where $x^{(n)} := (x_1, \dots, x_n)$ and
\begin{equation} \label{lawobservations}
\alpha_n(A_1\times\dots\times A_n) := \int_{\ps}\left[\prod_{i=1}^n \mu(A_i\,|\, \theta)\right] \pi(\ud\theta), 
\end{equation}
so that $\P[T \in B\,|\,X_1, \dots, X_n] = \pi_n(B\,|\,X_1, \dots, X_n)$ is valid $\P$-a.s. for any $B \in \psa$. 

\begin{remark}
When the statistical model $\mu(\cdot\,|\, \cdot)$ is dominated by some $\sigma$-finite measure $\lambda$ on $(\ss, \ssa)$, with a relative family of $\lambda$-densities $\{f(\cdot\, |\, \theta)\}_{\theta \in \ps}$, then (a version of) the posterior distribution is given by the Bayes formula, that is we write
$$
\pi_n(B\,|\,x^{(n)}) = \frac{\int_B [\prod_{i=1}^n f(x_i\, |\, \theta)] \pi(\ud\theta)}{\int_{\ps} [\prod_{i=1}^n f(x_i\, |\, \theta)] \pi(\ud\theta)}
$$ 
for any set $B \in \psa$ and $\alpha_n$-a.e. $x^{(n)}$, while $\alpha_n$ turns out to be absolutely continuous with respect to the product measure $\lambda^{\otimes_n}$ with density function of the form
\begin{equation} \label{density_Bayes}
\rho_n(x_1, \dots, x_n) := \int_{\ps}\left[\prod_{i=1}^n f(x_i\,|\, \theta)\right] \pi(\ud\theta)\ . 
\end{equation}
\end{remark}

We say that the posterior distribution is (weakly) consistent at $\theta_0\in\ps$ if, as $n\rightarrow+\infty$, $\pi_n(U_0^c\,|\, \xi_1, \dots, \xi_n) \rightarrow 0$ holds in probability for any neighborhood $U_0$ of $\theta_0$, where $\xi^{(\infty)} := \{\xi_i\}_{i \geq 1}$ stands for a sequence of $\ss$-valued independent random variables identically distributed as $\mu_0(\cdot) := \mu(\cdot| \theta_0)$ (\citet[Definition 6.1]{GV(00)}). The non uniqueness of the posterior distribution $\pi_n$ requires additional regularity assumptions in order that $\pi_n(\cdot\,|\, \xi_1, \dots, \xi_n)$ is well-defined. PCRs strengthen the notion of Bayesian consistency, in the sense that they quantify the speed at which such neighborhoods may decrease to zero meanwhile still capturing most of the posterior mass. In particular, the definition of PCR can be stated as follows (\citet[Definition 8.1]{GV(00)}).

\begin{defi} \label{def:consistency}
A sequence $\{\epsilon_n\}_{n \geq 1}$ of positive numbers is a PCR at $\theta_0$ if, as $n\rightarrow+\infty$,
\begin{equation} \label{post_consistency}
\pi_n\left(\left\{\theta \in \ps : \ud_{\ps}(\theta,\theta_0) \geq M_n\epsilon_n\right\}\,|\, \xi_1, \dots, \xi_n \right) \rightarrow 0
\end{equation}
holds in probability for every sequence $\{M_n\}_{n \geq 1}$ of positive numbers such that $M_n \rightarrow \infty$. 
\end{defi}

Now, we present our approach to PCRs based on the Wasserstein distance. This is a new approach, which relies on four main steps that are outlined hereafter. The first step of our approach originates from a reformulation of Definition \ref{def:consistency} in terms of the so-called $p$-Wasserstein distance, for $p \geq 1$. In particular, to recall this concept in full generality, we denote by $(\Space, \ud_{\Space})$ an abstract separable metric space, and we denote by $\mathcal P(\Space)$ the relative space of all probability measures on $(\Space, \Bcr(\Space))$. Then, the $p$-Wasserstein distance is defined as
\begin{equation} \label{Wasserstein}
\Wp^{(\mathcal P(\Space))}(\gamma_1; \gamma_2) := \inf_{\eta \in \mathcal{F}(\gamma_1,\gamma_2)} \left(\int_{\Space^2}  [\ud_{\Space}(x, y)]^p\ \eta(\ud x\ud y) \right)^{1/p}
\end{equation}
for any $\gamma_1, \gamma_2 \in \mathcal P_p(\Space)$, where 
$$
\mathcal P_p(\Space):= \left\{\gamma \in \mathcal P(\Space)\ :\ \int_{\Space} [\ud_{\Space}(x, x_0)]^p \gamma(\ud x) < +\infty\, \ \text{for\ some\ } x_0 \in \Space\right\}
$$ 
and $\mathcal{F}(\gamma_1, \gamma_2)$ is the class of all probability measures on $(\Space^2, \Bcr(\Space^2))$ with $i$-the marginal $\gamma_i$, for $i=1,2$. See \citet[Chapter 7]{AGS(08)} and \citet[Proposition 7.1.5]{AGS(08)}. If we let $(\Space, \ud_{\Space}) = (\ps, \ud_{\ps})$, then we can reformulate Definition \ref{def:consistency} according to the next lemma; the proof is deferred to Appendix \ref{proof:lm1}

\begin{lm} \label{lm:Wpcr}
Assume that $\pi \in \mathcal P_p(\ps)$ and that $\mu_0^{\otimes_n} \ll \alpha_n$ is valid for any $n \in \naturals$. Then, $\pi_n(\cdot\,|\, \xi_1, \dots, \xi_n)$ is a well-defined random probability measure belonging to $\mathcal P_p(\ps)$ 
with $\P$-probability one, and
\begin{equation} \label{Wpcr}
\epsilon_n = \E\left[ \Wp^{(\mathcal P(\ps))}(\pi_n(\cdot\,|\, \xi_1, \dots, \xi_n); \delta_{\theta_0}) \right]
\end{equation}
gives a PCR at $\theta_0$, where $\delta_{\theta_0}$ denotes the degenerate distribution at $\theta_0$.
\end{lm}

The second step of our approach relies on the assumption of the existence of a suitable sufficient statistics. In particular, we assume the existence of another metric space, say $(\Sd, \ud_{\Sd})$, and the existence of a measurable map, say $\mathfrak{S}_n : \ss^n \rightarrow \Sd$, in such a way that the kernel $\pi_n(\cdot| \cdot)$ in \eqref{disintegration} can be represented by means of another kernel, say $\pi_n^{\ast}(\cdot\,|\, \cdot): \psa\times\Sd \rightarrow [0,1]$, according to the identity 
\begin{equation} \label{representation_kernel}
\pi_n(\cdot\,|\,x_1, \dots, x_n) = \pi_n^{\ast}\left(\cdot\,|\,\mathfrak{S}_n(x_1, \dots, x_n)\right)
\end{equation}
for all $(x_1, \dots, x_n) \in \ss^n$.  See \citet{FLR(00)}, and references therein, for the existence of sufficient statistics in relationship with the exchangeability assumption. Of course, when the statistical model $\mu(\cdot| \cdot)$ is
dominated, the existence of the sufficient statistics $\mathfrak{S}_n$ is implied by standard assumptions on the statistical model, such as the well-known Fisher-Neyman factorization criterion. 

The third step of our approach relies on the large $n$ asymptotic behavior of the random variable $\hat{S}_n := \mathfrak{S}_n(\xi_1, \dots, \xi_n)$. In particular, we assume the existence of a weak law of large numbers for $\hat{S}_n$, which means that there exists some (non random) $S_0 \in \Sd$ for which $\hat{S}_n \rightarrow S_0$ holds true in $\P$-probability, as $n \rightarrow +\infty$. Hereafter, for any sequence $\{\delta_n\}_{n \geq 1}$ of positive numbers, we denote by 
\begin{equation} \label{Sanov_new}
J_n(\delta_n) := \P[\ud_{\Sd}(\hat{S}_n, S_0) \geq \delta_n]
\end{equation}
the probability that $\hat{S}_n$ lies outside a $\delta_n$-neighborhood of $S_0$. Usually, $J_n(\delta_n)$ can be evaluated by means of concentration inequalities and large deviation principles.

Based on \eqref{representation_kernel}, the fourth step of our approach relies on a form of local Lipschitz-continuity for the kernel $\pi_n^{\ast}(\cdot\,|\, \cdot)$, which holds under suitable assumptions on the model $\mu(\cdot\,|\,\cdot)$ and the prior $\pi$. It corresponds to the existence of two sequences of positive numbers, say $\{\delta_n\}_{n \geq 1}$ and $\{L_0^{(n)}\}_{n \geq 1}$ such that, for each $n \in \N$,
\begin{equation} \label{Lipschitz_kernel}
\Wp^{(\mathcal P(\ps))}\left(\pi_n^{\ast}(\cdot\,|\,S_0); \pi_n^{\ast}(\cdot\,|\,S') \right) \leq L_0^{(n)} \ud_{\Sd}(S_0, S')
\end{equation}
holds for any $S'$ belonging to $\mathcal{U}_{\delta_n}(S_0) := \{S \in \Sd\text{ : } \ud_{\Sd}(S_0, S) < \delta_n\}$. We refer to \citet{DM(20a),DM(20b)} for a detailed treatment of the property of local Lipschitz-continuity, for fixed $n\in \N$, providing some quantitative estimates for $L_0^{(n)}$. Then, according to Lemma \ref{lm:Wpcr}, under the validity of \eqref{representation_kernel} and \eqref{Lipschitz_kernel}, we write
\begin{align} \label{splitPCR}
\epsilon_n &\leq \Wp^{(\mathcal P(\ps))}(\pi_n^{\ast}(\cdot\,|\, S_0); \delta_{\theta_0})\\
&\notag\quad + L_0^{(n)}\E[ \ud_{\Sd}(\hat{S}_n, S_0)]\\
&\notag\quad\quad+ \E[ \Wp^{(\mathcal P(\ps))}(\pi_n^{\ast}(\cdot\,|\, S_0); \pi_n^{\ast}(\cdot\,|\, \hat{S}_n)) \mathds{1}\{\hat{S}_n \not\in \mathcal{U}_{\delta_n}(S_0)\}].
\end{align}

Under additional assumptions, in Section \ref{sec_mainres} we develop a careful analysis of the three terms on the right-hand side of \eqref{splitPCR}, in order to show that they can be bounded by terms of more explicit quantities that behave like $n^{-\alpha}$, for some $\alpha > 0$. In particular, the first term is a non-random quantity which is equal to
\begin{equation} \label{Wasserstein_nonrandom}
\Wp^{(\mathcal P(\ps))}(\pi_n^{\ast}(\cdot\,|\, S_0); \delta_{\theta_0}) = \left(\int_{\ps} \ud_{\ps}^p(\theta, \theta_0) \pi_n^{\ast}(\ud\theta\,|\, S_0) \right)^{1/p},
\end{equation}
and it measures the speed of shrinkage of $\pi_n^{\ast}(\cdot| S_0)$ at $\theta_0$. Its evaluation is a pure analytical problem, which relies on an extension to infinite-dimensional spaces of the classical Laplace methods of approximating integrals. In \eqref{splitPCR}, the term
\begin{equation} \label{GlivenkoCantelli}
\varepsilon_{n,p}(\Sd,S_0) := \E[ \ud_{\Sd}(\hat{S}_n, S_0)]
\end{equation}
provides the speed of convergence of the mean law of large numbers, which is well-known, at least for the situations considered throughout this paper. The term $\mathds{1}\{\hat{S}_n \not\in \mathcal{U}_{\delta_n}(S_0)\}$ in \eqref{splitPCR} hints at an application of a large deviation principle. As for the $L_0^{(n)}$'s in \eqref{splitPCR}, the bounds provided in \citet{DM(20a),DM(20b)} show that they can be expressed in terms of weighted Poincar\'e-Wirtinger constants. As we will show below, a proper choice of the sequency $\{\delta_n\}_{n \geq 1}$ should entail that $\{L_0^{(n)}\}_{n \geq 1}$ is bounded or, at least, diverges at a controlled rate.

Critical to our analysis of the term $L_0^{(n)}$ is the so-called dynamic formulation of the $p$-Wasserstein distance, which is referred to as Wasserstein dynamics (\citet{BB(00)}). In particular, assume that $\Space$ is the norm-closure of some nonempty, open and connected subset of a separable Hilbert space $\Hbb$, and endowed with scalar product $\langle\cdot,\cdot\rangle$ and norm $\|\cdot\|$. Then, for any $\gamma_0, \gamma_1 \in \mathcal P_p(\Space)$
\begin{displaymath}
\left[\Wp^{(\mathcal P(\Space))}(\gamma_0; \gamma_1)\right]^p = \inf_{\{\gamma_t\}_{t \in [0,1]} \in AC^p[\gamma_0; \gamma_1]} \int_0^1 \int_{\Space} \|\mathbf v_t(x)\|^p \gamma_t(\ud x) \ud t,
\end{displaymath}
where $AC^p[\gamma_0; \gamma_1]$ is the space of all absolutely continuous curves in $\mathcal P_p(\Space)$ with $\mathrm L^p(0,1)$ metric derivative (w.r.t. $\mathcal W_p$) connecting $\gamma_0$ to $\gamma_1$, and 
$[0,1]\times\Space \ni(t,x)\mapsto\mathbf v_t(x)\in\Hbb$ is a Borel function such that for almost every $t\in (0,1)$ it holds
\begin{equation}\label{cont_eq}
\frac{\ud}{\ud s} \int_{\Space} \psi(x) \gamma_s(\ud x) \ _{\big| s=t} = \int_{\Space} \langle \mathbf v(x), \mathrm{D} \psi(x)\rangle\gamma_t(\ud x)\qquad\forall \ \psi\in C^1_b(\Space).
\end{equation}
Here, $\mathrm{D}\psi $ denotes the Riesz representative of the Frech\'et differential of the function $\psi$, and $\psi\in C^1_b(\Space)$ means that $\psi$ is the restriction to $\Space$ of a function in the class $C^1_b(\Hbb)$, that is $\psi$ is a bounded continuous function with bounded continuous Fr\'echet derivative on $\Hbb$. See \citet[Chapter 2]{dPZ(14)} for  spaces of continuous functions defined on Hilbert spaces, and \citet[Chapter 8]{AGS(08)} for a detailed account on the partial differential equation \eqref{cont_eq}.

For any fixed $t$ and given $\gamma_t$, it is natural to look for a solution $\mathbf v_t(\cdot)$ of Equation \eqref{cont_eq} in the form of a gradient, and therefore we may interpret \eqref{cont_eq} as an abstract elliptic equation, for which it is well-known that a critical role is played by Poincar\'e inequalities in the context of proving the existence and regularity of a solution.

\begin{defi}
We say that a probability measure $\mu$ on $(\Space, \Bcr(\Space))$ satisfies a weighted Poincar\'e inequality of order $p$ if there exists a constant $\mathcal C_p$ for which
\begin{equation}\label{eqpw}
\inf_{a\in\mathbb R}\left(\int_\Space |\psi(x) - a|^p\,\mu(\ud x)\right)^{\frac1p} \le \mathcal C_p \left(\int_\Space \|\mathrm{D} \psi(x)\|^p\,\mu(\ud x) \right)^{\frac1p}
\end{equation}
holds for every $\psi\in C^1_b(\Space)$. We denote by $\mathfrak C_p[\mu]$ the best constant $C_p$ in \eqref{eqpw}. In particular, for $p=2$ the best constant $\mathfrak C_2[\mu]$ may be characterized by means of
\begin{align*}
\left(\frac1{\mathfrak C_2[\mu]}\right)^2 = \inf_{\{\psi\in C^1_b(\Space)\,:\,\int_\Space \psi(x)\,\mu(\ud x)=0\,\wedge\,\int_\Space |\psi(x)|^2\,\mu(\ud x)=1\}}\int_\Space \|\mathrm{D}\psi(x)\|^2\,\mu(\ud x) . 
\end{align*}
Finally, if $\Space = \Hbb$ and $\mu$ is absolutely continuous with respect to a non-degenerate Gaussian measure, then the Fr\'echet derivative in \eqref{eqpw} can be replaced by the Malliavin derivative $\Dcal$,  
yielding the following weaker definition
\begin{align*}
\left(\frac1{\mathfrak C_2^{(M)}[\mu]}\right)^2 = \inf_{\{\psi\in C^1_b(\Space)\,:\,\int_\Space \psi(x)\,\mu(\ud x)=0\,\wedge\,\int_\Space |\psi(x)|^2\,\mu(\ud x)=1\}}\int_\Space \|\mathrm{\Dcal}\psi(x)\|^2\,\mu(\ud x) . 
\end{align*}
\end{defi}

We refer to the monographs of \citet{Bog(10),dP(14), dPZ(14)} for a detailed account of Malliavin calculus and related Sobolev spaces. 


\section{Main results on PCRs}\label{sec_mainres}
Following the approach to PCRs outlined in Section \ref{sect:main_problem}, we present two main results: i) a theorem on PCRs for the regular infinite-dimensional exponential family of statistical models; ii) a theorem on PCRs for a general dominated statistical model. 

\subsection{PCRs for the regular infinite-dimensional exponential family} \label{sect:statementsNew}

It is useful to recall the definition and some basic properties of the infinite-dimensional exponential family. In general,  classical results on exponential family may be extended to the infinite-dimensional setting through suitable arguments of convex analysis (\citet{BBC(01),BC(17)}).

\begin{defi} \label{def:Reg_Exp_Fam}
Let $\lambda$ be a $\sigma$-finite measure on $(\ss, \ssa)$, $\B$ be a separable Banach space with dual $\B^{\ast}$, and $ _{\B^{\ast}}\langle \cdot, \cdot \rangle_{\B}$ be the pairing between $\B$ and $\B^{\ast}$. Also, let $\Gamma$ be a nonempty open subset of $\B^{\ast}$, and let $\beta : \ss \rightarrow \B$ be a measurable map. If the interior $\Lambda$ of the convex hull of the support of $\lambda \circ \beta^{-1}$ is nonempty and
\begin{equation}
\int_{\ss} \exp\{ _{\B^{\ast}}\langle \gamma, \beta_x \rangle_{\B}\} \lambda(\ud x) < +\infty
\end{equation}
holds for any $\gamma \in \Gamma$, then the regular infinite-dimensional exponential family is a statistical model defined through the family of $\lambda$-densities $\{\varphi(\cdot\,|\,\gamma)\}_{\gamma  \in \Gamma}$, where 
\begin{equation} \label{density_exp}
\varphi(x\,|\,\gamma) = \exp\{ _{\B^{\ast}}\langle \gamma, \beta_x \rangle_{\B} - M_{\varphi}(\gamma)\} 
\end{equation}
with 
\begin{equation} \label{M_exp}
M_{\varphi}(\gamma) := \log\left(\int_{\ss} \exp\{ _{\B^{\ast}}\langle \gamma, \beta_x \rangle_{\B}\} \lambda(\ud x) \right).
\end{equation}
\end{defi}

\citet[Theorems 1.13, Theorem 2.2 and Theorem 2.7]{brown} state that $M_{\varphi}$ is a strictly convex function on $\Gamma$, lower semi-continuous on $\B^{\ast}$, of class $C^{\infty}(\Gamma)$ and analytic. In addition, \citet[Corollary 5.3]{barndorff} implies that $M_{\varphi}$ is steep (essentially smooth). Therefore, from \citet[Theorem 3.6]{brown} it holds that
\begin{equation} \label{definition_S}
\mathcal S : \gamma \mapsto \nabla M_{\varphi}(\gamma) = \int_{\X} \beta_x \varphi(x\,|\,\gamma) \lambda(\ud x)
\end{equation}
defines a smooth injective map from $\Gamma$ into $\B$, with dense range. Finally, \cite[Corollary 2.5]{brown} entails the identifiability of the model characterized by the densities \eqref{density_exp}.

To introduce the setting of our theorem on PCR, it is useful to express the statistical model $\mu(\cdot\,|\,\cdot)$ in terms of an infinite-dimensional exponential family. In this regard, we introduce a further measurable mapping $g : \ps \rightarrow \Gamma$ and write
\begin{equation} \label{mu_n_exp}
\mu(\ud x\,|\,\theta) = \varphi(x\,|\, g(\theta)) \lambda(\ud x)\ .  
\end{equation}
In the setting of \eqref{mu_n_exp}, we observe that the identity \eqref{representation_kernel} is satisfied with $\Sd = \B$, $\ud_{\Sd}(s_1, s_2) = \|s_1 - s_2\|_{\B}$ and $\mathfrak{S}_n(x_1, \dots, x_n) = n^{-1} \sum_{1\leq i\leq n} \beta_{x_i}$. Therefore, we write
\begin{equation} \label{S_n_exp}
\hat{S}_n = \frac 1n \sum_{i=1}^n \beta_{\xi_i} \ .
\end{equation}
Note that Equation \eqref{mu_n_exp} arises naturally from the assumption that the statistical model $\mu(\cdot\,|\,\cdot)$ is dominated, which provides a family $\{f(\cdot\,|\,\theta)\}_{\theta \in \Theta}$ of density functions. Accordingly, by assuming that $\ss$ is endowed of a richer metric structure, if $f(\cdot\,|\,\theta) > 0$ and $x \mapsto f(x\,|\,\theta)$ is continuous for any $\theta \in \Theta$, we write
\begin{displaymath}
\log f(x\,|\,\theta) = \int_{\ss} \log f(y\,|\,\theta) \delta_x(\ud y).
\end{displaymath}
The functions $g$ and $\beta$ then arise from the mapping $y \mapsto \log f(y\,|\,\theta)$ and the measure $\delta_x$ through (some sort of) integration-by-parts, if this is admitted, or through classical Fourier transformation arguments, such as the Plancherel formula.

According to \citet[Corollary 7.10]{LT(91)}, under the assumption
\begin{equation} \label{first_mom_exp}
\E\left[\|\beta_{\xi_i}\|_{\B}\right] = \int_{\ss} \|\beta_x\|_{\B}\ \mu_0(\ud x) = \int_{\ss} \|\beta_x\|_{\B}\ \varphi(x\,|\, g(\theta_0)) \lambda(\ud x) < +\infty\ ,
\end{equation}
we set
\begin{equation} \label{S_0_exp}
S_0 = \int_{\ss} \beta_x \mu_0(\ud x) = \int_{\ss} \beta_x \varphi(x\,|\, g(\theta_0)) \lambda(\ud x) 
\end{equation}
in the sense of Bochner integral, and conclude the strong law of large numbers, i.e. $\hat{S}_n \rightarrow S_0$ holds $\P$-a.s., as $n\rightarrow+\infty$. Now, we set $M(\theta) = M_{\varphi}(g(\theta))$ and then define 
\begin{equation} \label{kernel_exp}
\pi_n^{\ast}(\ud\theta\,|\, b) := \frac{\exp\{n[ _{\B^{\ast}}\!\langle g(\theta), b\rangle_{\B} - M(\theta)]\} \pi(\ud\theta)}{\int_{\ps} \exp\{n[\langle g(\tau), b\rangle_{\B} - M(\tau)]\} \pi(\ud\tau)},
\end{equation}
provided that
\begin{equation} \label{int_exp_post}
\int_{\ps} \exp\{n\ _{\B^{\ast}}\!\langle g(\tau), b\rangle_{\B} \} \pi(\ud\tau) < +\infty
\end{equation}
for any $n \in \N$ and $b \in \B$. We remark that \eqref{int_exp_post} is a necessary assumption for the existence of the posterior distribution.  
Now, we state the theorem on PCRs in the setting of infinite-dimensional exponential families; the proof is deferred to Appendix \ref{proof:new_thm}.

\begin{thm} \label{new_therem_PCR}
Let $p \geq 1$ be a fixed number.
Let $\{\varphi(\cdot\,|\,\gamma)\}_{\gamma  \in \Gamma}$ be a regular infinite-dimensional exponential family according to \emph{Definition} \ref{def:Reg_Exp_Fam}. Let $\ps$ be an open, connected subset of some separable Hilbert
space $\Hbb$. Let $g : \ps \rightarrow \Gamma$ be a measurable mapping for which representation \eqref{mu_n_exp} is in force. For a fixed $\theta_0 \in \ps$, suppose that:
\begin{enumerate}    
\item[i)] \eqref{first_mom_exp} is valid;    
\item[ii)] \eqref{int_exp_post} holds for any $n \in \N$ and $b \in \B$;
\item[iii)] $\int_{\ps} \|\theta\|^{ap} \pi(\ud\theta) < +\infty$ for some $a > 1$;
\item[iv)] there exists a sequence $\{\delta_n\}_{n \geq 1}$ of positive numbers for which \eqref{Lipschitz_kernel} is valid for any $n \in \N$, 
with $\Sd = \B$, $\ud_{\Sd}(s_1, s_2) =  \|s_1 - s_2\|_{\B}$ and suitable positive constants $L_0^{(n)}$. 
\end{enumerate}    
Then, for the PCR $\epsilon_n$ at $\theta_0$ it holds
\begin{align}\label{Main_PCR1}
\epsilon_n &\lesssim \left(\int_{\ps} \|\theta - \theta_0\|_{\ps}^p \pi_n^{\ast}(\ud\theta\,|\, S_0) \right)^{\frac 1p}\\
&\notag\quad+ \|\theta_0\|_{\ps}\ \P\left[\hat{S}_n \not\in \mathcal{U}_{\delta_n}(S_0) \right] \nonumber \\
&\notag\quad\quad + \left( \E\left[ \int_{\ps} \|\theta\|_{\ps}^{ap} \pi_n^{\ast}(\ud\theta\,|\, \hat{S}_n) \right] \right)^{\frac{1}{ap}} \left(\P\left[\hat{S}_n \not\in \mathcal{U}_{\delta_n}(S_0)\right]\right)^{1 - \frac{1}{ap}} \\
&\notag\quad\quad\quad+ L_0^{(n)} \E[\|\hat{S}_n- S_0\|_{\B}] 
\end{align}
where $\pi_n^{\ast}(\cdot\,|\,\cdot)$ is given by  \eqref{kernel_exp}, $S_0$ and $\hat{S}_n$ are as in \eqref{S_0_exp} and \eqref{S_n_exp}, respectively, and $\mathcal{U}_{\delta_n}(S_0) := \{S \in \Sd\ |\ \|S_0-S\|_{\B} < \delta_n\}$.
\end{thm}

Theorem \ref{new_therem_PCR} provides an implicit form for PCRs. That is, the large $n$ asymptotic behaviour of the terms on the right-hand side of \eqref{Main_PCR1} must be further investigated to obtain a more explicit expression for the corresponding PCR. In this regard, it is useful to rewrite $\pi^{\ast}_{n}$ in terms of the Kullback-Leibler divergence. That is, if $\mathcal S \circ g$ is injective and $b$ belongs to the range of $\mathcal S \circ g$, then
\begin{equation} \label{representation_kernel_kullback}
\pi_n^{\ast}(\ud\theta\,|\, b) = \frac{\exp\{-n \kksf(\theta \,|\, \theta_b)\}\pi(\ud\theta)}{\int_{\ps} \exp\{-n \kksf(\tau\, |\, \theta_b)\}\pi(\ud\tau)}
\end{equation}
where $\theta_b = (\mathcal S \circ g)^{-1}(b)$ and 
\begin{equation} \label{definition_kullback}
\kksf(\theta\,|\,\theta') := \int_{\ss} \left[\ln\left(\frac{f(x\,|\,\theta')}{f(x\,|\,\theta)}\right)\right] f(x\,|\,\theta')\ud x\ . 
\end{equation}
denotes the Kullback-Leibler divergence. See Appendix \ref{proof:lm2} for the proof of Equation \eqref{representation_kernel_kullback}. 
It is natural to expect that the main contribution to PCRs arises from the first and the fourth term on the right-hand side of \eqref{Main_PCR1}, which provide general algebraic rates of convergence to zero. Hereafter, we investigate the large $n$ asymptotic behaviour of the terms on the right-hand of \eqref{Main_PCR1}. More explicit results in terms of PCRs will be presented in Section \ref{sect:illustrations} with respect to the application of Theorem \ref{new_therem_PCR} in the context of the regular parametric model, the multinomial model, the finite-dimensional and the infinite-dimensional logistic-Gaussian model and the infinite-dimensional linear regression.

\subsubsection{First term on the right-hand of \eqref{Main_PCR1}} \label{sect:first_term}
We start by considering the large $n$ asymptotic behaviour of the first term on the right-hand side of \eqref{Main_PCR1}. In particular, from \eqref{representation_kernel_kullback}, we can rewrite this terms as
\begin{align} \label{Gibbs}
&\int_{\ps} \|\theta - \theta_0\|_{\ps}^p \pi_n^{\ast}(\ud\theta\,|\, S_0)\\
&\quad= \frac{\int_{\ps} \|\theta - \theta_0\|_{\ps}^p \exp\{n[ _{\B^{\ast}}\!\langle g(\theta), S_0\rangle_{\B} - M(\theta)]\} \pi(\ud\theta)}{\int_{\ps} 
\exp\{n[\langle g(\tau), S_0\rangle - M(\tau)]\} \pi(\ud\tau)} \nonumber \\
&\notag\quad=\frac{\int_{\ps} \|\theta - \theta_0\|_{\ps}^p \exp\{-n \kksf(\theta\, |\, \theta_0)\}\pi(\ud\theta)}{\int_{\ps} \exp\{-n \kksf(\theta\, |\, \theta_0)\}\pi(\ud\theta)}\ .
\end{align}
The last expression of  \eqref{Gibbs} shows the ratio of two Laplace integrals, and therefore the Laplace method of approximating integrals can be applied. In the finite-dimensional setting, i.e. $\ps \subseteq \rd$, the Laplace approximation method is well-known (\citet{Bre(94),Won(01)}), and it leads to the following proposition.

\begin{prp} \label{prp:Laplace_finite}
In the case that $\ps \subseteq \rd$, assume that $\pi$ has a continuous density $q$ with respect to the Lebesgue measure, with $q(\theta_0) > 0$, and that $\theta \mapsto \kksf(\theta\,|\, \theta_0)$ is a $C^2$-function
with a strictly positive definite Hessian at $\theta_0$, which coincides with the Fisher information matrix $\mathrm{I}[\theta_0]$ at $\theta_0$. Let $\int_{\ps}  |\theta|^p\pi(\ud \theta) < +\infty$ be fulfilled for some
$p \geq 1$.
Finally, suppose that for any $\delta > 0$ there exists $c(\delta) >0$ such that 
\begin{equation} \label{Kullback_deltapositive}
\inf_{|\theta - \theta_0| \geq \delta} \kksf(\theta \,|\, \theta_0) \geq c(\delta).
\end{equation}
 Then, for any $p > 0$, there hold
\begin{displaymath}
\int_{\ps}  |\theta - \theta_0|^p e^{-n \kksf(\theta\,|\, \theta_0)} \pi(\ud \theta) \sim \frac 12\left(\frac{2}{n}\right)^{\frac{d+p}{2}} \!\!\!\!\! \Gamma\left(\frac{d+p}{2}\right) 
\end{displaymath}
and
\begin{displaymath}
\frac{\int_{S_d(1)} \{\langle z, \mathrm{I}[\theta_0]^{-1} z\rangle \}^{p/2} \ud S(z)}{\sqrt{\mathrm{det}[\mathrm{I}(\theta_0)]}} \\
\int_{\ps} e^{-n \kksf(\theta\, |\, \theta_0)} \pi(\ud \theta) \sim \left(\frac{2\pi}{n}\right)^{d/2} \frac{1}{\sqrt{\mathrm{det}[\mathrm{I}(\theta_0)]}}
\end{displaymath}
as $n \rightarrow +\infty$, where $S_d(1) := \{z \in \rd\ |\ \|z\| = 1\}$, $\ud S$ denotes the surface measure and $\langle \cdot, \cdot \rangle$ stands for the standard scalar product in $\rd$. Thus, under these assumptions, 
\begin{equation} \label{Laplace1}
\int_{\ps} \|\theta - \theta_0\|_{\ps}^p \pi_n^{\ast}(\ud\theta\,|\, S_0) = O(n^{-p/2})
\end{equation}
as $n\rightarrow +\infty$. 
\end{prp}  

It is interesting to observe that the inequality \eqref{Kullback_deltapositive} is a sort of strengthening of the  so-called Shannon-Kolmogorov information inequality. See, e.g., \citet[Chapter 17]{Fer(96)}. In particular, because of \eqref{Kullback_deltapositive}, integrals on the whole $\Theta$ can be reduced to integrals over balls centered at $\theta_0$, as integration over the complement of any such ball yields exponentially small 
quantities with respect to $n$.

According to Proposition \ref{prp:Laplace_finite}, in the finite-dimensional setting the prior distribution does not affect the large $n$ asymptotic behaviour of the first term on the right-hand side of \eqref{Main_PCR1}. Differently from the standard finite-dimensional setting, the literature on the Laplace approximation method in the infinite-dimensional setting appears to be not well developed. That is, to the best of our knowledge, infinite-dimensional Laplace approximations are limited to the case in which the measure $\pi$ is a Gaussian measure (\citet{AS(19),AS(99)}). Unfortunately, this literature does not cover the case in which the Hessian of the map $\theta \mapsto \kksf(\theta \,|\, \theta_0)$ at $\theta_0$ is not coercive (uniformly elliptic), which is precisely the case of interest in our specific problem. The next proposition covers this critical gap; the proof is deferred to Appendix \ref{proof:Laplace_infinite}. The proposition is of independent interest in the context of the classical Laplace method.

\begin{prp}\label{prp:Laplace_infinite}
Let $\ps$ be a separable Hilbert space with scalar product $\langle \cdot, \cdot\rangle$, and let $\pi$ be the non-degenerate Gaussian measure $\mathcal{N}(m,Q)$, with $m \in \ps$ and 
$Q$ a trace-class operator. For fixed $\theta_0 \in \ps$,
assume that $\theta \mapsto \kksf(\theta\,|\,\theta_0)$ belongs to $C^{2+q}(\ps)$ for some $q \in (0,1]$, and that its Hessian at $\theta_0$, which coincides with the Fisher information operator 
$\mathrm{I}(\theta_0)$ at $\theta_0$, is a 
compact self-adjoint linear operator from $\ps$ into itself, with trivial kernel. Suppose there exists an orthonormal Fourier basis $\{\mathbf e_k\}_{k \geq 1}$ of $\ps$ which diagonalizes simultaneously both $Q$ and $\mathrm{I}(\theta_0)$, so that
\begin{equation} \label{eigenvalues_logistic}
Q[\mathbf e_k] = \lambda_k \mathbf e_k \qquad\qquad \mathrm{I}(\theta_0)[\mathbf e_k] = \gamma_k \mathbf e_k 
\end{equation}
are valid with two suitable sequences $\{\lambda_k\}_{k \geq 1}$ and $\{\gamma_k\}_{k \geq 1}$ of strictly positive numbers that go to zero as $k \rightarrow +\infty$, with $\{\lambda_k\}_{k \geq 1} \in \ell_1$. Finally, assume there exist two other Hilbert spaces $\Kbb$ and $\Vbb$ such that 
\begin{enumerate}
\item[i)] $\Vbb \subset \ps \subset \Kbb$ with continuous, dense embeddings;
\item[ii)] an interpolation inequality like
\begin{equation} \label{interpolation}
\|\theta\|_{\ps} \lesssim \|\theta\|_{\Kbb}^{1/r} \|\theta\|_{\Vbb}^{1/s} 
\end{equation}
holds for any $\theta \in \Vbb$ with conjugate exponents $r,s > 1$ such that $r < 1 + q/2$;
\item[iii)] for all $\theta \in \Vbb$, the inequalities
\begin{align} 
\kksf(\theta\, |\, \theta_0) &\geq \phi(\|\theta- \theta_0\|_{\Kbb})  \label{lower_kullback} \\
\langle \theta-\theta_0, \mathrm{I}(\theta_0)[\theta-\theta_0]\rangle &\gtrsim \|\theta- \theta_0\|_{\Kbb}^2 \label{lower_fisher}
\end{align}
are valid with some monotone non-decreasing function $\phi : [0,+\infty) \to [0,+\infty)$ such that $\phi(x) = O(x^2)$ as $x \to 0^+$; 
\item[iv)] $\pi(\Vbb) = 1$ and $\int_{\Vbb} e^{t \|\theta\|_{\Vbb}} \pi(\ud\theta) < +\infty$ for some $t > 0$.
\end{enumerate}
Then, as $n\rightarrow +\infty$, the following expansion
\begin{equation} \label{Laplace2}
\int_{\ps} \|\theta - \theta_0\|_{\ps}^2 \pi_n^{\ast}(\ud\theta\,|\, S_0) = O\left(\sum_{k=1}^{\infty} \frac{\lambda_k}{n \lambda_k \gamma_k + 1} \right) + O\left(\sum_{k=1}^{\infty} \frac{\omega_k^2}{(n \lambda_k \gamma_k + 1)^2} \right)
\end{equation}
holds with the sequence $\{\omega_k\}_{k \geq 1} \in \ell_2$ given by $(\theta_0 - m) = \sum_{k=1}^{\infty} \omega_k\mathbf e_k$.
\end{prp} 

\begin{remark}
In the infinite-dimensional setting, the assumption \eqref{Kullback_deltapositive} is, in general, too strong. Conditions \eqref{lower_kullback}--\eqref{lower_fisher}, combined with the interpolation \eqref{interpolation}, constitute 
a reasonable set of assumptions that allow a quite general treatment in the applications. It is worth noticing that \eqref{Kullback_deltapositive}, as well as \eqref{lower_kullback}, is expressed in the form of a lower 
bound for $\kksf(\theta\, |\, \theta_0)$. These bounds are conceptually opposite with respect to the so-called ``prior mass condition'' required in the standard theory (\citet[Theorem 8.9, inequality (8.4)]{GV(00)}), which
is usually proved by means of upper bounds for $\kksf(\theta\, |\, \theta_0)$. See, e.g. the upper bounds for $\kksf(\theta\, |\, \theta_0)$ in Lemma 2.5 of \citet{GV(00)}.  
\end{remark} 

\begin{remark}
With respect to Proposition \ref{prp:Laplace_finite}, the statement of Proposition \ref{prp:Laplace_infinite} is confined to the case $p=2$. There are no technical limitations for treating the more general case $p \neq 2$, though $p=2$ yields to a more readable (conclusive) result. 
\end{remark}

\begin{remark} \label{rmk:Carlen}
Assumption \eqref{eigenvalues_logistic} is not necessary to obtain PCRs. However, without this assumption, the resulting PCR would have a complicated form, which may be recovered from the proof. For example, let $\ps_N$ be the finite-dimensional subspace of $\ps$ obtained by the linear span of $\{\mathbf e_1, \dots, \mathbf e_N\}$, let $Q_N$ denote the $N\times N$ matrix that represents the restriction of $Q$ to $\ps_N$, after projecting the range of such restriction again on $\ps_N$, and let $\mathrm{I}_N(\theta_0)$ denote the $N\times N$ matrix associated to the restriction just explained of the operator $\mathrm{I}(\theta_0)$ to $\ps_N$. If $Q_N$ and $\mathrm{I}_N(\theta_0)$
are non-singular, then the first term on the right-hand side of \eqref{Laplace2} can be replaced by 
\begin{equation} \label{Trace_truncated}
\lim_{N\to +\infty} \mathrm{Tr}\left[\left(n\mathrm{I}_N(\theta_0) +  Q_N^{-1}\right)^{-1}\right],
\end{equation}
which is not as clear as the series $\sum_{k=1}^{\infty} \lambda_k/(n \lambda_k \gamma_k + 1)$. An analogous operation can be performed with respect to the second term on the right-hand side of \eqref{Laplace2}. 

Moreover, the above argument can be reinforced by resorting to some \emph{trace inequalities}, as explained in \cite{Car(10)}. In particular, we assume there exists another compact, self-adjoint operator 
$\mathrm{I}^{\ast}$
such that $\mathrm{I}(\theta_0) \geq \mathrm{I}^{\ast}$ in the sense of quadratic forms, i.e. 
$$
\langle \theta, \mathrm{I}(\theta_0)[\theta]\rangle \geq \langle \theta, \mathrm{I}^{\ast}[\theta]\rangle
$$
for any $\theta \in \ps$. Whence, upon denoting by $\mathrm{I}^{\ast}_N$ the restriction of $\mathrm{I}^{\ast}$ to $\ps_N$ as above, we have $\mathrm{I}_N(\theta_0) \geq \mathrm{I}^{\ast}_N$
and, consequently, $n\mathrm{I}_N(\theta_0) +  Q_N^{-1} \geq n \mathrm{I}^{\ast}_N +  Q_N^{-1}$. By the L\"owner--Heinz theorem, the mapping $t \mapsto -t^{-1}$ is \emph{operator monotone}, yielding that
$$
\mathrm{Tr}\left[\left(n\mathrm{I}_N(\theta_0) +  Q_N^{-1}\right)^{-1}\right] \leq \mathrm{Tr}\left[\left(n \mathrm{I}^{\ast}_N +  Q_N^{-1}\right)^{-1}\right] \ .
$$
See again \cite{Car(10)} for the details. Therefore, if the orthonormal Fourier basis $\{\mathbf e_k\}_{k \geq 1}$ of $\ps$ diagonalizes simultaneously both $Q$ and $\mathrm{I}^{\ast}$ (instead of
$\mathrm{I}(\theta_0)$), so that
\begin{equation} \label{eigenvalues_logistic2}
Q[\mathbf e_k] = \lambda_k \mathbf e_k \qquad\qquad \mathrm{I}^{\ast}[\mathbf e_k] = \gamma_k^{\ast} \mathbf e_k 
\end{equation}
are valid with suitable strictly positive $\gamma_k^{\ast}$'s that go to zero as $k \rightarrow +\infty$, then by Proposition \ref{prp:Laplace_infinite}
\begin{equation} \label{Laplace2_ast}
\int_{\ps} \|\theta - \theta_0\|_{\ps}^2 \pi_n^{\ast}(\ud\theta\,|\, S_0) \lesssim  \sum_{k=1}^{\infty} \frac{\lambda_k}{n \lambda_k \gamma_k^{\ast} + 1} 
+ \sum_{k=1}^{\infty} \frac{\omega_k^2}{(n \lambda_k \gamma_k^{\ast} + 1)^2} \ . 
\end{equation}
\end{remark} 

Proposition \ref{prp:Laplace_infinite} shows that the large $n$ asymptotic behavior of the first term on the right-hand side of \eqref{Main_PCR1} is worse than $1/n$, which is the large $n$ asymptotic behaviour obtained in Proposition \ref{prp:Laplace_finite} with $p=2$. For example, by taking the first term on the right-hand side of \eqref{Laplace2} into account,
if $\lambda_k \sim k^{-(1+a)}$ and $\gamma_k \sim k^{-b}$ as $k\rightarrow+\infty$, for some $a,b > 0$, a straightforward calculation shows that
\begin{equation} \label{elem_asymptotic1}
\sum_{k=1}^{\infty} \frac{k^{-(1+a)}}{n k^{-(1+a+b)} + 1} \sim n^{-{\frac{a}{1+a+b}}}
\end{equation}
holds as $n\rightarrow+\infty$. As for the second term on the right-hand side of \eqref{Laplace2}, it can be made identical to zero by choosing $m = \theta_0$, that is by means of centering the Gaussian prior  at $\theta_0$. However, if $\lambda_k \sim k^{-(1+a)}$, $\gamma_k \sim k^{-b}$ and $\omega^2_k \sim k^{-(1+c)}$ as $k\rightarrow+\infty$, for some choice of $a,b,c > 0$ with $c < 2(1+a+b)$, then
$$
\sum_{k=1}^{\infty} \frac{k^{-(1+c)}}{(n k^{-(1+a+b)} + 1)^2} \sim n^{-{\frac{c}{1+a+b}}}
$$
holds as $n\rightarrow+\infty$. Therefore, if $c < a$ this second term is slower than the one in \eqref{elem_asymptotic1}, whilst if $c > a$ it is negligible with respect to that 
term. Again on \eqref{elem_asymptotic1}, it is interesting to notice what happens if the eigenvalues $\lambda_k$'s approach zero very rapidly, like
$\lambda_k \sim e^{-k}$, for example. Another straightforward calculation shows that
$$
\sum_{k=1}^{\infty} \frac{e^{-k}}{n e^{-k} k^{-b} + 1} \sim \frac{(\log n)^{b+1}}{n}
$$
holds as $n\rightarrow+\infty$. A refinement of this argument  entails that the large $n$ asymptotic behavior of the right-hand side of \eqref{Laplace2} can be made arbitrarily close to the rate $1/n$, for example by choosing $\lambda_k \sim e^{-k^r}$ and $\gamma_k \sim k^{-b}$ and $\omega^2_k \sim k^{-(1+c)}$ for some $r, b,c > 0$, with arbitrarily large $r$. By recalling that the first term on the  right-hand side of \eqref{Main_PCR1} coincides with the square root of the left-hand side of \eqref{Main_PCR1}, this argument shows that the PCR is arbitrarily close to $1/\sqrt{n}$. It is reasonable to guess that the minimax (classical) risk should go to zero as fast as $1/\sqrt{n}$, though we are not aware of any result proving such a behaviour.

A merit of Proposition \ref{prp:Laplace_infinite} is to show explicitly that, within the infinite-dimensional setting, PCRs are influenced by three quantities that do not appear in finite-dimensional setting of Proposition \ref{prp:Laplace_finite}: i) the rate of approach to zero of the sequence $\{\lambda_k\}_{k \geq 1}$, which measures the ``regularity of the prior''; ii) the rate of approach to zero of the sequence $\{\gamma_k\}_{k \geq 1}$, which measures the ``regularity of the model''; iii) the rate of approach to zero of the sequence $\{\omega_k\}_{k \geq 1}$, which measures how close is $\theta_0$ to $m$. Finally, we notice that the space $\Vbb$ is linked with the Cameron-Martin space associated to $\pi$, which must be included in $\Vbb$.

\subsubsection{Second and third term on the right-hand of \eqref{Main_PCR1}} Now, we consider the large $n$ asymptotic behaviour of the second term and of the third term on the right-hand side of \eqref{Main_PCR1}. Both these terms depend explicitly on
\begin{align} \label{LargeDev}
\P\left[\hat{S}_n \not\in \mathcal{U}_{\delta_n}(S_0) \right] = \P\left[\|\hat{S}_n - S_0\|_{\B} \geq \delta_n \right] = \P\left[\|\hat{S}_n - \E[\hat{S}_n]\|_{\B} \geq \delta_n \right]. 
\end{align}
Note that the tail probability in \eqref{LargeDev} is directly related to classical concentration inequalities for sum or random variables. Besides well-know Bernstein-type concentration inequalities for real-valued random variables (\citet{BLM(13),DZ(98)}), some useful generalizations or extension can be found in, e.g., \citet{GiNi(16)}, \citet{LT(91)}, \citet{Pin(86)} and \citet{Yur(70)}. In particular, for a suitable choice of the sequence $\{\delta_n\}_{n \geq 1}$, such that a constant sequence or a vanishing sequence at an algebraic rate, the term \eqref{LargeDev} goes to zero at suitable exponential rates, and therefore it provides a negligible contribution in the right-hand side of \eqref{Main_PCR1}. 

The third term  on the right-hand side of \eqref{Main_PCR1} includes the posterior moment $\E[\int_{\ps} \|\theta\|_{\ps}^{ap} \pi_n^{\ast}(\ud\theta\,|\, \hat{S}_n)] = \E[ \int_{\ps} \|\theta\|_{\ps}^{ap} \pi_n(\ud\theta\,|\, \xi_1, \dots, \xi_n)]$. In particular, an application of H\"older's inequality shows that such a moment is bounded from above by
\begin{align*}
&\left(\int_{\ps} \|\theta\|^{\rho ap}\pi(\ud\theta)\right)^{1/\rho}\left(\int_{\ss^n} \left[\frac{\prod_{i=1}^n f(x_i\,|\,\theta_0)}{\rho_n(x_1, \dots, x_n)} \right]^{\rho'} \!\!\!\!\!\rho_n(x_1, \dots, x_n) \prod_{i=1}^{n}\lambda(\ud x_i)\right)^{1/\rho'} 
\end{align*}
for conjugate exponents $\rho, \rho' > 1$, provided that $\int_{\ps} \|\theta\|^{\rho ap}\pi(\ud\theta) < +\infty$. It is useful to recall that the density function $\rho_n$ has been defined in \eqref{density_Bayes}.
Accordingly, the second factor above coincides with the $\rho'$-th moment of a martingale, since
\begin{align*}
&\int_{\ss^n} \left[\frac{\prod_{i=1}^n f(x_i\,|\,\theta_0)}{\rho_n(x_1, \dots, x_n)} \right]^{\rho'} \!\!\!\!\!\rho_n(x_1, \dots, x_n) \prod_{i=1}^{n}\lambda(\ud x_i)= \E\left[\left(\frac{\prod_{i=1}^n f(X_i\,|\,\theta_0)}{\rho_n(X_1, \dots, X_n)} \right)^{\rho'} \right]
\end{align*}
and
\begin{displaymath}
\E\left[ \frac{\prod_{i=1}^{n+1} f(X_i\,|\,\theta_0)}{\rho_{n+1}(X_1, \dots, X_{n+1})}\,\Big|\, X_1, \dots, X_n  \right] = \frac{\prod_{i=1}^n f(X_i\,|\,\theta_0)}{\rho_{n}(X_1, \dots, X_{n})}\ . 
\end{displaymath}
At this stage, a possible resolutive strategy may rely on well-known bounds for moments of martingales (\citet{DFJ(68)}). As for the term $\E[\|\hat{S}_n- S_0\|_{\B}]$, by means of a direct application of Lyapuonov's inequality, we can write that
\begin{displaymath}
\E[\|\hat{S}_n- S_0\|_{\B}] \leq \left( \E[\|\hat{S}_n- S_0\|_{\B}^2] \right)^{1/2}
\end{displaymath} 
and the right-hand side typically goes to zero as $1/\sqrt{n}$. Besides the obvious case in which $\B$ coincides with a separable Hilbert space, we refer to \citet{Nem(00)}, \citet{Mas(07)} and \citet{MR(13)} for the case in which we have $\B = \ell_p(\rd)$.

\subsubsection{Fourth term on the right-hand of \eqref{Main_PCR1}} \label{sect:fourth_term}
Finally, we consider the large $n$ asymptotic behaviour of the fourth term on the right-hand side of \eqref{Main_PCR1}. In particular, this term involves the constant $L_0^{(n)}$, whose treatment requires to recall some fundamental notions of infinite-dimensional calculus. Given $g : \ps \to \Bdual$, the Fr\'echet differential $\mathfrak{D}_{\theta}[g]$ of $g$ is now meant as 
a bounded linear operator from $\ps$ to $\Bdual$ such that $g(\theta + \delta) = g(\theta) + \mathfrak{D}_{\theta}[g](\delta) + o(\|\delta\|_{\ps})$, as $\delta \to 0$ in $\ps$, and
$$
\|\mathfrak{D}_{\theta}[g]\|_{\ast} := \sup_{\|\delta\|_{\ps}\ \leq 1} \|\mathfrak{D}_{\theta}[g](\delta)\|_{\Bdual}\ .
$$

Here, we consider the case $p=2$. It should be recalled that the theory of weighted Poincar\'e constant has been 
mainly focused on the two cases $p=1$ and $p=2$ (see, e.g., \citet{BBCG(08)}). We choose only the latter case in order to avoid other technical problems connected with the Wasserstein dynamic when $p=1$.
See, e.g., the first comment opening Section 8.3 of \citet{AGS(08)}. Therefore, in order to obtain an explicit upper bound for the constant $L_0^{(n)}$ it is useful to consider the following proposition; the proof is deferred to Appendix \ref{proof:Bayes_Club}

\begin{prp} \label{prop:Bayes_Club}
In addition to the assumptions of \emph{Theorem} \ref{new_therem_PCR}, suppose that $g \in C^1(\ps; \B^{\ast})$, that $\int_{\ps} \|\mathfrak{D}_{\theta}[g]\|_{\ast}^2\ \pi(\ud\theta) < +\infty$, and that 
map $\mathcal S \circ g$ is continuous. Then, for the constant $L_0^{(n)}$ in \eqref{Main_PCR1} we can put
\begin{equation} \label{L0_bound1}
L_0^{(n)} = n \sup_{S \in \mathcal{U}_{\delta_n}(S_0)} \left\{ \Cpw[\pi_n^{\ast}(\cdot\,|\, S)] \right\}^2 \left( \int_{\ps} \|\mathfrak{D}_{\theta}[g]\|_{\ast}^2\ \pi_n^{\ast}(\ud\theta\,|\, S) \right)^{1/2}\ .
\end{equation}
In addition, if 
\begin{equation} \label{Laplace_rafforzato}
\sup_{n\in \N}\ \ \sup_{\theta' \in \mathcal{V}_{\delta_n}(\theta_0)}\
\frac{\int_{\ps} \|\mathfrak{D}_{\theta}[g]\|_{\ast}^2 \exp\{-n \kksf(\theta\, |\, \theta')\}\pi(\ud\theta)}{\int_{\ps} \exp\{-n \kksf(\theta\, |\, \theta')\}\pi(\ud\theta)} =: \mathcal B(g) < +\infty
\end{equation}
holds with $\mathcal{V}_{\delta_n}(\theta_0) := (\mathcal S \circ g)^{-1}(\mathcal{U}_{\delta_n}(S_0))$, then
\begin{equation} \label{L0_bound2}
L_0^{(n)} \leq \sqrt{\mathcal B(g)}\  n \sup_{S \in \mathcal{U}_{\delta_n}(S_0)} \left\{ \Cpw[\pi_n^{\ast}(\cdot\,|\, S)] \right\}^2\ . 
\end{equation}
\end{prp}

\begin{remark} 
When $\pi$ is a Gaussian measure on the infinite-dimensional Hilbert space $\ps$, an analogous statement can be formulated with the Fr\'echet derivative replaced by the Malliavin derivative. Hence,  
for the constant $L_0^{(n)}$ in \eqref{Main_PCR1} we can set
\begin{equation} \label{L0_bound1_M}
L_0^{(n)} = n \sup_{S \in \mathcal{U}_{\delta_n}(S_0)} \left\{ \Cpw^{(M)}[\pi_n^{\ast}(\cdot\,|\, S)] \right\}^2 \left( \int_{\ps} \|\Dcal_{\theta}[g]\|_{\ast}^2\ \pi_n^{\ast}(\ud\theta\,|\, S) \right)^{1/2},
\end{equation}
and if 
\begin{equation} \label{Laplace_rafforzato_M}
\sup_{n\in \N}\ \sup_{\theta' \in \mathcal{V}_{\delta_n}(\theta_0)}\ \
\frac{\int_{\ps} \|\Dcal_{\theta}[g]\|_{\ast}^2 \exp\{-n \kksf(\theta\, |\, \theta')\}\pi(\ud\theta)}{\int_{\ps} \exp\{-n \kksf(\theta\, |\, \theta')\}\pi(\ud\theta)} =: \mathcal B_M(g) < +\infty
\end{equation}
holds with $\mathcal{V}_{\delta_n}(\theta_0) := (\mathcal S \circ g)^{-1}(\mathcal{U}_{\delta_n}(S_0))$, then the following inequality holds true
\begin{equation} \label{L0_bound2_M}
L_0^{(n)} \leq \sqrt{\mathcal B_M(g)}\  n \sup_{S \in \mathcal{U}_{\delta_n}(S_0)} \left\{ \Cpw^{(M)}[\pi_n^{\ast}(\cdot\,|\, S)] \right\}^2\ . 
\end{equation}
\end{remark}

Denote by $\Rightarrow$ the weak convergence of probability measures on $(\ps, \mathscr{B}(\ps))$. Verifying the validity of \eqref{Laplace_rafforzato} represents a strengthening of the fact that, as $n\rightarrow+\infty$ 
\begin{displaymath}
\frac{\exp\{-n \kksf(\theta \,|\, \theta')\}\pi(\ud\theta)}{\int_{\ps} \exp\{-n \kksf(\tau\, |\, \theta')\}\pi(\ud\tau)} \Rightarrow \delta_{\theta'}. 
\end{displaymath}
This may be proved by means of the same arguments as in the proofs of Proposition \ref{prp:Laplace_finite} and Proposition \ref{prp:Laplace_infinite}. According to Proposition \ref{prop:Bayes_Club}, to conclude it remains to make more explicit the large $n$ asymptotic behaviour of the weighted Poincar\'e-Wirtinger constant $\Cpw[\pi_n^{\ast}(\cdot\,|\, S')]$. In the finite-dimensional setting, i.e. $\ps \subseteq \rd$, the representation \eqref{Gibbs} shows that the posterior distribution, or better $\pi_n^{\ast}(\cdot\,|\, \cdot)$, characterizes Gibbsean (Boltzmann) probability distributions. Properties of the Kullback-Leibler divergence entail that the mapping $\theta \mapsto \kksf(\theta\, |\, \theta')$ is non-negative and vanishes iff $\theta = \theta'$ \citet[Chapter 17]{Fer(96)}. Moreover, under standard regularity assumptions for $f(\cdot\ |\ \cdot)$ (\citet[Chapter 18]{Fer(96)}), the aforesaid mapping proved also to be strictly convex, at least in finite dimension. In this context, there are several conditions that entail the upper bound
\begin{displaymath}
[\mathfrak C_2(\pi_n^{\ast}(\cdot\,|\, \mathcal S \circ g(\theta')))]^2 \leq \frac{C(\theta')}{n}
\end{displaymath}
for every $n \in \naturals$ and positive constant $C(\theta')$. In particular, the simplest condition to quote is the so-called Bakry-Emery condition, characterized by the fact that 
\begin{equation}\label{prelcond}
\mathrm{Hess}[\kksf(\cdot\, |\, \theta')](\theta) \geq \rho \mathrm{Id}
\end{equation}
for some $\rho > 0$, with $ \mathrm{Id}$ being the identity matrix, uniformly with respect to $\theta \in \ps$, in conjunction with the hypothesis that $\pi(\ud\theta) = e^{-U(\theta)}\ud\theta$ for some $U \in C^2(\ps)$. Some generalizations of the condition \eqref{prelcond} are given in the next proposition, which specifies some results that have first appeared in \citet{BBCG(08)}.
\begin{prp}[\citet{DM(20b)}]\label{francesi}
Let $U$ and $G$ be elements of $\mathrm{C}^2(\ps)$, bounded from below and such that $\mathrm{Hess}(G(\theta))\ge \alpha \mathrm{Id}$ and $\mathrm{Hess}(U(\theta))\ge h \mathrm{Id}$ (in the sense of quadratic forms) whenever $|\theta|\le R$, for some $\alpha>0$, $R>0$ and $h\in\reals$. 
\begin{itemize}
\item[(1)] If, in addition, there exist $c>0$ and $\ell\in\reals$ such that 
$\theta\cdot\nabla G(\theta)\ge c|\theta|$ and $\theta\cdot \nabla U(\theta)\ge \ell|\theta|$ whenever $|\theta|\ge R$, then
\begin{align*}
&\left[ \mathfrak C_2\left( \frac{e^{-nG(\theta) - U(\theta)} \ud\theta}{\int_{\ps} e^{-nG(\tau) - U(\tau)} \ud\tau}\right)\right]^2\\ 
&\quad\le \frac{\alpha n+h+(cn+\ell-d_R+ nG_R+U_R)\,C_R}{(\alpha n+h)\,(cn+\ell-1-d_R)} \sim \frac 1n
\end{align*}
for every $n>(-h/\alpha)\vee((d_R+1-\ell)/c)$, where $d_R:=(d-1)/R$, $G_R:=\sup_{B_R}|\nabla G|$, $U_R:=\sup_{B_R}|\nabla U|$ and $C_R$ is an explicit universal constant only depending on $R$. 
\item[(2)] If, in addition, there exist $c_1>0$, $c_2>0$ such that
$$
|\nabla G(\theta)|^2\ge 2c_1+c_2\, [\Delta G(\theta)+\nabla G(\theta)\cdot\nabla U(\theta)]_+
$$
whenever $|\theta|\ge R$, then
\[
\left[ \mathfrak C_2\left( \frac{e^{-nG(\theta) - U(\theta)} \ud\theta}{\int_{\ps} e^{-nG(\tau) - U(\tau)} \ud\tau}\right)\right]^2  \le \frac{\alpha n+h+e^{\omega_R}(c_1n+G_R^*+W_R)} {(\alpha n+h)c_1n} \sim \frac 1n
\]  
for every $n>(1+1/c_2)\vee(-h/\alpha)$, where $G_R^*:=\sup_{B_R}|\Delta G|$, $W_R:=\sup_{B_R}|\nabla U||\nabla G|$ and $\omega_R:=\sup_{B_R}G-\inf_{\ps}G$. 
\end{itemize}
\end{prp}

According to Proposition \ref{francesi}, in the finite-dimensional setting the prior distribution does not affect the large $n$ asymptotic behaviour of the weighted Poincar\'e-Wirtinger constant $\Cpw[\pi_n^{\ast}(\cdot\,|\, S')]$. A similar phenomenon has been observed in the study of the first term on the right-hand side of \eqref{Main_PCR1}. Differently from the finite-dimensional setting, the literature on weighted Poincar\'e-Wirtinger constants in the the infinite-dimensional setting appears to be not well developed . To the best of our knowledge, in the infinite-dimensional setting, upper bounds on weighted Poincar\'e-Wirtinger constants are limited to the case of Gibbsean (Boltzmann) measures, that is  measures of the form $\exp\{-nG(\theta)\}\pi(\ud\theta)$ with $G$ being a smooth convex function and $\pi$ being an infinite-dimensional Gaussian measure (\citet[Chapters 10-11]{DaPrato(06)}). While this is the case of interest in our problem, the upper bounds available in the literatures are not sharp for large values of $n$, and therefore they can not be applied. The next proposition covers this critical gap by providing results involving Malliavin calculus; the proof is 
deferred to Appendix \ref{proof:Poincare_infinite}. The proposition is of independent interest in the context of weighted Poincar\'e-Wirtinger constants.

\begin{prp}\label{prp:Poincare_infinite}
Let $\ps$ be a separable Hilbert space, and let $\pi$ be the non-degenerate Gaussian measure $\mathcal{N}(m,Q)$, with $m \in \ps$ and $Q$ a trace-class operator. Let $\mathrm{G}_0 :\ps \to \ps$ be a compact linear operator, with
trivial kernel. Let $G$ be an element of $\mathrm{C}^2(\ps)$, bounded from below and such that $\mathrm{Hess}(G(\theta))\ge \mathrm{G}_0$ (in the sense of operators) whenever $\|\theta\|_{\ps} \le R$,
for some $R>0$. Suppose there exists a Fourier orthonormal basis $\{\mathbf e_k\}_{k \geq 1}$ of $\ps$ which diagonalizes simultaneously both $Q$ and $\mathrm{G}_0$, that is
\begin{equation} \label{factor_Poincare_Infinite}
Q[\mathbf e_k] = \lambda_k \mathbf e_k \qquad\qquad \mathrm{G}_0[\mathbf e_k] = \eta_k \mathbf e_k 
\end{equation}
for two suitable sequences $\{\lambda_k\}_{k \geq 1}$ and $\{\eta_k\}_{k \geq 1}$ of strictly positive numbers that go to zero as $k \rightarrow +\infty$, with $\{\lambda_k\}_{k \geq 1} \in \ell_1$. 
\begin{itemize}
\item[(1)] Suppose, in addition, there exists $c>0$ such that $\theta\cdot\Dcal_{\theta} G\ge c\|\theta\|_{\ps}$ whenever $\|\theta\|_{\ps} \ge R$. Then, for every $n>\mathrm{Tr}[Q](1 + 1/R)/c$, it holds
\begin{align}
&\left[ \mathfrak C_2^{(M)}\left( \frac{e^{-nG(\theta)} \pi(\ud\theta)}{\int_{\ps} e^{-nG(\tau)} \pi(\ud\tau)}\right)\right]^2\nonumber \\ 
&\quad \lesssim \frac{1 + C_R(1+\tau_n+nG_R) \max_{k \in \N} \left\{\frac{\lambda_k}{n\lambda_k \eta_k + 1}\right\} }{cn - \mathrm{Tr}[Q](1 + 1/R)} \nonumber\\
&= O\left(\max_{k \in \N} \left\{\frac{\lambda_k}{n\lambda_k \eta_k + 1}\right\}\right) \label{tesi_francesi_infinite}
\end{align}
where $G_R:=\sup_{B_R} \|\Dcal_{\theta} G\|$ and $C_R$ is an explicit universal constant only depending on $R$. 
\item[(2)] Suppose, in addition, there exist $c_1>0$, $c_2>0$ such that
\begin{equation} \label{growth_francesi}
\|\Dcal_{\theta} G\|^2 \ge 2c_1+c_2\, [\mathfrak{L}_{\pi} G(\theta)]_+
\end{equation}
whenever $\|\theta\|\ge R$, where $\Dcal_{\theta}$ and $\mathfrak{L}_{\pi}$ denote the Malliavin derivative and the Malliavin-Laplace operator associated to $\pi$, respectively.   
Then, for every $n> 1 + 1/c_2$, it holds
\begin{align*}
&\left[ \mathfrak C_2^{(M)}\left( \frac{e^{-nG(\theta)} \pi(\ud\theta)}{\int_{\ps} e^{-nG(\tau)} \pi(\ud\tau)}\right)\right]^2 \\
&\quad \lesssim \frac{1 + e^{\omega_R}(C_1 n + G^*_R) \max_{k \in \N} \left\{\frac{\lambda_k}{n\lambda_k \eta_k + 1}\right\} }{c_1 n} \\
&= O\left(\max_{k \in \N} \left\{\frac{\lambda_k}{n\lambda_k \eta_k + 1} \right\}\right)
\end{align*}
where $\omega_R := \sup_{B_R} G - \inf_{\ps} G$ and $G^*_R := \sup_{B_R} [|\mathfrak{L}_{\pi}[G]| + \| \Dcal_{\theta}[G] \|^2]$.
\end{itemize}
\end{prp} 

Proposition \ref{prp:Poincare_infinite} shows that the large $n$ asymptotic behavior of the Poincar\'e-Wirtinger constant $\Cpw^{(M)}[\pi_n^{\ast}(\cdot\,|\, S')]$ is worse than $1/n$, which is the large $n$ asymptotic behaviour obtained in Proposition \ref{francesi}. By straightforward calculations, as $n\rightarrow+\infty$
\begin{equation}\label{elem_asymptotic2}
\max_{k \in \N} \left\{ \frac{k^{-(1+a)}}{n k^{-(1+a+b)} + 1} \right\} \sim n^{-\frac{a+1}{1+a+b}} 
\end{equation}
holds for any $a,b > 0$. A particular merit of Proposition \ref{prp:Poincare_infinite} consists in showing explicitly that, within the infinite-dimensional setting, PCRs are influenced by two quantities that do not appear in finite-dimensional setting of Proposition \ref{francesi}: i) the rate of approach to zero of the sequence $\{\lambda_k\}_{k \geq 1}$, which measures the ``regularity of the prior''; ii) the rate of approach to zero of the sequence $\{\eta_k\}_{k \geq 1}$, which measures another ``regularity of the model''. A similar phenomenon has been observed in the study of the first term on the right-hand side of \eqref{Main_PCR1}. To conclude, we observe that, under the assumptions of Proposition \ref{prp:Laplace_infinite}, we can apply Equation \eqref{representation_kernel_kullback} to rewrite the right-hand side of \eqref{L0_bound2} as follows
\begin{align} \label{kullback_trick}
&\sup_{S' \in \mathcal{U}_{\delta_n}(S_0)} \left\{ \Cpw^{(M)}[\pi_n^{\ast}(\cdot\,|\, S')] \right\}^2\\
&\notag\quad = \sup_{\theta' \in (\mathcal S \circ g)^{-1}(\mathcal{U}_{\delta_n}(S_0))} 
\left[ \mathfrak C_2^{(M)}\left( \frac{e^{-n\kksf(\theta\,|\,\theta')} \pi(\ud\theta)}{\int_{\ps} e^{-n\kksf(\tau\,|\,\theta')} \pi(\ud\tau)}\right)\right]^2,
\end{align}
and then observe that the role of $\theta'$ is now confined to the multiplicative constants that appear on the right-hand sides of the various inequalities that we considered. Thus, in order to handle the supremum, it is enough to check the boundedness of such multiplicative constants by standard arguments of continuity.

We conclude this section by summarizing our results on the large $n$ asymptotic behaviours of the terms on the right-hand side of \eqref{Main_PCR1}. The second and the third term go to zero exponentially fast, and this holds true independently on the dimension of the statistical model. This confirms that the main contribution to the PCR arise from the first and the fourth term, which give generally algebraic rates of convergence to zero. In the finite-dimensional setting, the first and the fourth term go to zero as $n^{-1}$, which is the optimal rate. In the infinite-dimensional setting, the first and the second term go to zero according to Proposition \ref{prp:Laplace_infinite} and Proposition \ref{prp:Poincare_infinite}. At least when $\eta_k \sim \gamma_k \sim k^{-b}$, $\lambda_k \sim k^{-(1+a)}$ and $\omega_k^2 \sim k^{-(1+c)}$ with $a,b,c >0$, Equation \eqref{elem_asymptotic1} and Equation \eqref{elem_asymptotic2} show that the first term on the right-hand side of \eqref{Main_PCR1} is asymptotically equivalent to
\begin{displaymath}
n^{-\frac{a}{2(a+b+1)}} + n^{-\frac{c}{2(a+b+1)}},
\end{displaymath}
whereas the fourth term on the right-hand side of \eqref{Main_PCR1} is asymptotically equivalent to 
\begin{displaymath}
n^{-\frac{a+1-b}{2(a+b+1)}},
\end{displaymath}
at least assuming that $\E[\|\hat{S}_n- S_0\|_{\B}] \sim n^{-1/2}$. This completes our analysis of PCRs in the setting of infinite-dimensional exponential families. Some applications of these results will be presented in Section \ref{sect:illustrations} with respect to specific statistical models.

\subsection{PCRs for a general dominated statistical model} \label{sect:statements1}

We present a more general version of Theorem \ref{new_therem_PCR}, which relies on the assumption that both the sample space $\ss$ and parameter space $\ps$ have richer analytical structures. As in Section  \ref{sect:statementsNew}, we confine to the case $p=2$.  In particular, the setting that we consider may be summarized through the following assumptions.

\begin{assum} \label{ass:thm1} 
The set $\ss$, the parameter space $\ps$ and the statistical model  $\mu(\cdot\,|\,\cdot)$ are such that
\begin{enumerate}
\item[i)] $\ss$ coincides with an open, connected subset of $\rm$ with Lipschitz boundary, and $\ssa = \mathscr{B}(\ss)$. With minor changes of notation, $\ss$ could also coincide with a smooth Riemannian manifold without boundary 
of dimension $m \in \naturals$.
\item[ii)] $\ps$ coincides with an open, connected subset of a separable Hilbert space of dimension $d \in \naturals \cup \{+\infty\}$.
\item[iii)] $\mu(\cdot\,|\,\cdot)$ is dominated by the $m$-dimensional Lebesgue measure, i.e. $\mu(A\,|\,\theta) = \int_A f(x\,|\,\theta) \ud x$ for every $A \in \ssa$, where 
$x \mapsto f(x\,|\,\theta) > 0$ is a probability density function for any $\theta \in \ps$.
\item[iv)] $(x,\theta) \mapsto f(x\,|\,\theta) \in C^2(\ss\times\ps)$;
\item[v)] the model $\{f(\cdot\,|\,\theta)\}_{\theta \in \ps}$ is $C^2$-regular at $\theta_0$ (as in \cite[Theorem 18]{Fer(96)})
\item[vi)] for any $\theta$, there exist positive constants $b(\theta), c(\theta)$ for which 
\begin{equation} \label{buond_statmodel}
|\log f(x\,|\,\theta)| \leq b(\theta)(1 + |x|^2) \qquad\mbox{and}\qquad |\nabla_x \log f(x\,|\,\theta)|\le c(\theta)(1+|x|)
\end{equation}
hold for every $x \in \ss$;
\item[vii)] $\pi \in \mathcal P_2(\ps)$, with full support; 
\item[viii)] $\mu_0 \in \mathcal P_2(\ss)$.
\end{enumerate}
\end{assum}

The setting of infinite-dimensional exponential families, considered in Section \ref{sect:fourth_term}, is a popular example that satisfies Assumptions \ref{ass:thm1}. Now, we state the theorem on PCRs in the setting of Assumptions \ref{ass:thm1}; the proof is deferred to Appendix \ref{proof:thm1}

\begin{thm} \label{main_thm1}
Within the setting specified by Assumptions \ref{ass:thm1}, 
\eqref{representation_kernel} is fulfilled with
\begin{equation} \label{Bayes_measure}
\pi_n^{\ast}(\ud\theta\,|\,\gamma) := \frac{\exp\{ n\int_{\ss} \log f(y \,|\, \theta) \gamma(\ud y)\} }{\int_{\ps} \exp\{ n\int_{\ss} \log f(y \,|\, \tau) \gamma(\ud y)\} \pi(\ud \tau)} \pi(\ud\theta)
\end{equation}
where $\gamma \in \Sd = \mathcal P_2(\ss)$. 
Moreover, $\eqref{Lipschitz_kernel}$ holds relatively to a suitable choice of a $\Wdue^{(\pms)}$-neighborhood $V_0^{(n)}$ of $\mu_0(\cdot) := \mu(\cdot\,|\,\theta_0)$, provided that
\begin{align} \label{Ln}
L_0^{(n)} &:= n \sup_{\gamma \in V_0^{(n)}} [\mathfrak C_2(\pi_n^{\ast}(\cdot\,|\,\gamma))]^2\\
&\notag\quad\times \left(\int_{\ps}\int_{\ss} \Big\|\mathfrak D_{\theta} \frac{\nabla_x f(x\,|\,\theta)}{f(x\,|\,\theta)}\Big\|^2
\gamma(\ud x) \pi_n^{\ast}(\ud \theta\,|\,\gamma) \right)^{1/2} < +\infty
\end{align}
for any $n \in \naturals$. Thus, the assumptions of Lemma \ref{lm:Wpcr} are fulfilled and a PCR at $\theta_0$ is given by 
\begin{align} \label{PCR_main1}
\epsilon_n &=  3\left( \frac{\int_{\ps} \|\theta-\theta_0\|^2 e^{-n \kksf(\theta\,|\,\theta_0)} \pi(\ud\theta)}{\int_{\ps} e^{-n \kksf(\theta\,|\,\theta_0)} \pi(\ud\theta)} \right)^{1/2}  \\
&+ L_0^{(n)} \varepsilon_{n,2}(\ss,\mu_0) + 2\|\theta_0\| \P[\empiric^{(\xi)} \not\in V_0^{(n)}] \nonumber \\
&+ \E\left[\left( 2 \frac{\int_{\ps} \|\theta\|^2 \left[\prod_{i=1}^n f(\xi_i\,|\,\theta) \right] \pi(\ud\theta)}{\int_{\ps} \left[\prod_{i=1}^n f(\xi_i\,|\,\theta) \right]\pi(\ud\theta)} \right)^{1/2} 
\mathds{1}\{\empiric^{(\xi)}\not\in V_0^{(n)}\} \right]\nonumber 
\end{align}
where $\empiric^{(\xi)} := n^{-1} \sum_{1\leq i\leq n} \delta_{\xi_i}$ and $\varepsilon_{n,p}(\ss,\mu_0) := \E[ \Wp^{(\mathcal P(\ss))}(\mu_0; \empiric^{(\xi)})]$ is the speed of mean Glivenko-Cantelli.
\end{thm}

From Theorem \ref{main_thm1}, we observe that if $V_0^{(n)} = \mathcal P_2(\ss)$ makes $L_0^{(n)}$ finite for every $n \in \naturals$, then the expression on the right-hand side of \eqref{PCR_main1} reduces to the first two terms. Similarly to Theorem \ref{new_therem_PCR}, Theorem \ref{main_thm1} provides an implicit form for the PCR, thus requiring to further investigate the large $n$ asymptotic behaviour of the terms on the right-hand side of \eqref{PCR_main1}. The posterior distribution appears in \eqref{Ln} and \eqref{PCR_main1}, meaning that further work is required to obtain more explicit terms. In general, it is possible to get rid of $\pi_n^{\ast}$ in  \eqref{Ln} and \eqref{PCR_main1}, thus reducing \eqref{Ln} and \eqref{PCR_main1} to expressions that involve only the statistical model and the prior distribution. The first term on the right-hand side of \eqref{PCR_main1} has the same form as in \eqref{Main_PCR1}, meaning that the Laplace method plays a critical role in the study of these term. Such a term can be handled as described in Proposition \ref{prp:Laplace_finite} and Proposition \ref{prp:Laplace_infinite}. With regards to $\varepsilon_{n,2}(\ss,\mu_0)$, we recall from \citet[Theorem 1]{FG(15)} that, if $\int_{\ss} |x|^q \mu_0(\ud x) < +\infty$ for some $q > 2$, then
\begin{align*}
&\varepsilon_{n,2}(\ss,\mu_0)\\
&\quad \leq C(q,m) \left(\int_{\ss} |x|^q \mu_0(\ud x)\right)^{1/q}\\
&\quad\quad\times \left\{\begin{array}{ll}
n^{-1/4} + n^{-(q-2)/(2q)} &\mbox{if}\ m=1,2,3\ \mbox{and}\ q\neq 4\\
n^{-1/4}\sqrt{\log(1+n)} + n^{-(q-2)/(2q)} &\mbox{if}\ m=4\ \mbox{and}\ q\neq 4\\
n^{-1/m} + n^{-(q-2)/(2q)} &\mbox{if}\ m>4\ \mbox{and}\ q\neq m/(m-2)
\end{array}\right.
\end{align*}
with some positive constant $C(q,m)$. Under some more restrictive assumptions, $\varepsilon_{n,2}(\ss,\mu_0)$ is of order $O(n^{-1/2})$, which is optimal in the dimension $1$ (\citet[Section 5]{BoLe(19)}). In the dimension $2$, the optimal rate is $\sqrt{(\log n)/n}$ (\citet{AST(19)}), whereas for $m\geq3$ the optimal rate is $n^{-1/m}$ (\citet{Ta(94b)}). Lastly, when $\ss$ has infinite dimension, logarithmic rates have been obtained in \citet{Jing(20)}. With regards to $\P[\empiric^{(\xi)} \not\in V_0^{(n)}]$, we refer to \citet[Theorem 2.7]{BGV(07)}. In particular, if $\int_{\ss} |x|^q \mu_0(\ud x) < +\infty$ for some $q \geq 1$, then
$$
\P[\Wdue(\empiric^{(\xi)}; \mu_0) > t] \leq B(q,m) t^{-q} \times \left\{\begin{array}{ll}
n^{-q/4} &\mbox{if}\ q >4\\
n^{1 - q/2} &\mbox{if}\ q \in [2,4) 
\end{array}\right.
$$
for any $t >0$ and $n \in \naturals$, with some positive constant $B(q,m)$. Exponential bounds can be also obtained upon requiring that $\int_{\ss} e^{\alpha |x|} \mu_0(\ud x) < +\infty$ for some $\alpha > 0$. See \citet[Theorem 2.8]{BGV(07)}. In the next corollary we show that, under additional assumptions, similar bounds hold true for the other terms appearing on the right-hand side of \eqref{PCR_main1}; the proof is deferred to Appendix \ref{proof:coro1}.

\begin{cor} \label{coro1}
In addition to the hypotheses of Theorem \ref{main_thm1}, suppose that there exist constants some $C>0$ and $\beta \geq 2$ for which 
\begin{equation} \label{bound_Kullback}
\kksf(\theta|\theta_0) \geq C \min\{\|\theta - \theta_0\|^{\beta}, 1\}
\end{equation}
holds for all $\theta \in \ps$. Moreover, assume that  
\begin{displaymath}
\int_{\ss} |x|^q \mu_0(\ud x) < +\infty
\end{displaymath}
holds for some constants $q > 4$, and 
\begin{equation} \label{bound_moment_posterior}
 \mathfrak{M}_{n,r} := \E\left[\int_{\ps} \|\theta\|^r \pi_n(\ud\theta\,|\,\xi_1, \dots, \xi_n) \right] < +\infty
\end{equation}
for all $n \in \naturals$ and some $r > 2$. 
Then, if the neighborhood $V_0^{(n)}$ has the form 
$$
\{\gamma \in \mathcal P_2(\ss)\ |\ \Wdue(\gamma; \mu_0) \leq Kn^{-a}\}
$$ 
for some $K >0$ and $a \in [0, 1/4)$,
for the PCR given in \eqref{PCR_main1} we obtain the new bound
\begin{align}\label{PCR_main2}
\epsilon_n &\leq C_1 n^{-1/\beta} + L_0^{(n)} \varepsilon_{n,2}(\ss,\mu_0) + 2\|\theta_0\| K^{-q} B(q,m) n^{q(a - 1/4)} \\
& + C_2 n^{[q(r-1)(a - 1/4)]/r} \mathfrak{M}_{n,r}  \nonumber
\end{align}
with suitable positive constants $C_1$ and $C_2$.
\end{cor}

From Corollary \ref{coro1}, the posterior distribution appears in  \eqref{bound_moment_posterior} and \eqref{PCR_main2}. With regards to \eqref{bound_moment_posterior}, this term is typically available in an explicit form, even if the posterior is not explicit. In general, a possible strategy may rely on well-known bounds for moments of martingales. With regards to $L_0^{(n)}$, this term can be handled as described in Proposition \ref{francesi} and Proposition \ref{prp:Poincare_infinite}, that is by inequalities for the weighted Poincar\'e-Wirtinger constant. To conclude it remains to handle with
\begin{displaymath}
\sup_{\gamma \in V_0^{(n)}} \int_{\ps}\int_{\ss} \Big\|\mathfrak D_{\theta} \frac{\nabla_x f(x\,|\,\theta)}{f(x\,|\,\theta)}\Big\|^2 \gamma(\ud x) \pi_n^{\ast}(\ud \theta\,|\,\gamma), 
\end{displaymath}
which is expected to be bounded with respect to $n$, in regular situations. To deal with this term, a possible strategy consists in obtaining inequality of the form
\begin{displaymath}
\int_{\ss} \Big\|\mathfrak D_{\theta} \frac{\nabla_x f(x\,|\,\theta)}{f(x\,|\,\theta)}\Big\|^2 \gamma(\ud x) \leq C_{\gamma} W(\theta)
\end{displaymath}
for a suitable constant $C_{\gamma}$ and a suitable function $W$. This particular point will be made more precise in Section \ref{sect:illustrations} with respect to some specific statistical models.


\section{Applications}\label{sect:illustrations}

\subsection{Regular parametric models} \label{sect:regularmodels}

Consider the case of  dominated Bayesian statistical models with a finite-dimensional parameter $\theta \in \ps \subset \rd$. Accordingly, we start by considering the set of Assumptions \ref{ass:thm1}, 
with $d \in \naturals$, along with the hypotheses of Theorem \ref{main_thm1}. In this setting, the Kullback-Leibler divergence $\kksf(\theta\,|\,\theta_0)$ is a $C^2$ function, whose Hessian at $\theta_0$ just coincides with the Fisher information 
matrix at $\theta_0$. Whence,
\begin{equation} \label{TaylorKullback}  
\kksf(\theta\,|\, \theta_0) = \frac 12\  ^t(\theta - \theta_0) \mathrm{I}[\theta_0] (\theta - \theta_0) + o(|\theta -\theta_0|^2)
\end{equation}
as $\theta \rightarrow \theta_0$. Finally, we assume \eqref{Kullback_deltapositive}. Therefore, we can apply Proposition \ref{prp:Laplace_finite} to get
\begin{equation} 
\Wp(\pi_n^{\ast}(\ud\theta | \mu_0), \delta_{\theta_0}) = \left(\frac{\int_{\ps}  |\theta - \theta_0|^p e^{-n \kksf(\theta\, |\, \theta_0)} \pi(\ud \theta)}{\int_{\ps} e^{-n \kksf(\theta\, |\, \theta_0)} \pi(\ud \theta)} 
\right)^{1/p} = O\big(n^{-1/2}\big)
\end{equation}
as $n \rightarrow +\infty$. Now, we discuss the behavior of the constant $L_0^{(n)}$, as $n$ goes to infinity. First, we would like to stress that there are plenty of conditions that entail
\begin{displaymath}
[\mathfrak C_2(\pi_n^{\ast}(\cdot\,|\,\gamma))]^2 \leq \frac{C(\gamma)}{n}
\end{displaymath}
for every $n \in \naturals$ and some positive constant $C(\gamma)$. We consider the double integral 
$$
\int_{\ps}\int_{\ss} \Big\|\mathrm D_{\theta} \frac{\nabla_x f(x\,|\,\theta)}{f(x\,|\,\theta)}\Big\|^2\gamma(\ud x) \pi_n^{\ast}(\ud \theta\,|\,\gamma)\ .
$$
First, if $f(x\,|\,\theta) = \exp\{\langle\theta, T(x)\rangle - M(\theta)\}$, that is the model is an element of the exponential family in the canonical form, then we notice that
$D_{\theta} \frac{\nabla_x f(x|\theta)}{f(x|\theta)}$ reduces to a $d\times m$ matrix whose entries are given by $\partial_{x_j} T_i(x)$, for $j = 1, \dots, m$ and $i = 1, \dots, d$. Therefore, the study of the above double integral boils down to that of the much simpler expressions $\int_{\ss} |\partial_{x_j} T_i(x)|^2 \gamma(\ud x)$, which are independent of $n$. More generally, we can  reduce the problem by resorting to the Laplace method for approximating probability integrals, from which we have that
\begin{displaymath}
\int_{\ps}\Big\|\mathrm D_{\theta} \frac{\nabla_x f(x\,|\,\theta)}{f(x\,|\,\theta)}\Big\|^2 \pi_n^{\ast}(\ud \theta|\gamma) \sim \Big\|\mathrm D_{\theta} \frac{\nabla_x f(x\,|\,\theta)}{f(x\,|\,\theta)}\Big\|^2 \ _{\Big|\ \theta = \theta^*(\gamma)}
\end{displaymath}
as $n \to +\infty$, where $\theta^*(\gamma)$ denotes a maximum point of the mapping $\theta \mapsto \int_{\ss} \log f(y|\theta) \gamma(\ud y)$. Therefore, if the above right-hand side proves to be positive, 
a reasonable plan to prove global boundedness of $L_0^{(n)}$ with respect to $n$ can be based on the following to steps. First, we check the validity of an inequality like
\begin{displaymath}
\sup_{n\in \naturals} \int_{\ps}\Big\|\mathrm D_{\theta} \frac{\nabla_x f(x\,|\,\theta)}{f(x|\theta)}\Big\|^2 \pi_n^{\ast}(\ud \theta\,|\,\gamma) \leq C \Big\|\mathrm D_{\theta} 
\frac{\nabla_x f(x\,|\,\theta)}{f(x\,|\,\theta)}\Big\|^2 \ _{\Big|\ \theta = \theta^*(\gamma)}
\end{displaymath}
for every $\gamma$ belonging to a $\Wdue^{(\pms)}$-neighborhood of $\mu_0$, where $C$ is a positive constant possibly depending on the fixed neighborhood. Second, we prove global boundedness (for $\gamma$
varying in the neighborhood) of the following integral
\begin{displaymath}
\int_{\ss} \Big\|\mathrm D_{\theta} \frac{\nabla_x f(x\,|\,\theta)}{f(x\,|\,\theta)}\Big\|^2 \ _{\Big|\ \theta = \theta^*(\gamma)}\gamma(\ud x) < +\infty\ .
\end{displaymath}
To fix ideas in a more concrete way, we consider the Gaussian case, where $\theta = (\mu, \Sigma)$ and
\begin{displaymath}
f(x\,|\,\theta) = (2\pi)^{-m/2} \frac{1}{\sqrt{\textrm{det}(\Sigma)}} \exp\left\{-\frac 12 (x-\mu)^t \Sigma^{-1}(x-\mu)\right\} \qquad x \in \reals^m\ .
\end{displaymath}
Note that the mapping $\theta \mapsto \int_{\ss} \log f(y\,|\,\theta) \gamma(\ud y)$ depends on $\gamma$ only through its moments of order 1 and 2. Thus, the above strategy reduces to an ordinary 
finite-dimensional maximization problem, very similar to the question of finding the maximum likelihood estimator. Finally, the last term in on the right-hand side of \eqref{PCR_main1} can be treated as in Corollary \ref{coro1}, by studying the asymptotic behavior of some posterior $r$-moment as in \eqref{bound_moment_posterior}. We state two propositions that summarize the above considerations. The former result holds when Theorem \ref{new_therem_PCR} can be applied and gives the optimal rate, while the latter result ensues from Theorem \ref{main_thm1}. 
\begin{prp}
Assume that there exist a separable Banach space $\B$ with dual $\B^{\ast}$ and two measurable maps $\beta : \ss \to \B$ and $g : \ps \to \B^{\ast}$ for which \eqref{mu_n_exp} is in force.    
If the assumptions of \emph{Theorem} \ref{new_therem_PCR} and \emph{Propositions} \ref{prp:Laplace_finite}, \ref{prop:Bayes_Club} and \ref{francesi} are met, then as $n\rightarrow+\infty$
\begin{displaymath}
\epsilon_n = O\big(n^{-1/2}\big),
\end{displaymath}
which is the optimal rate. 
\end{prp}
\begin{prp}
Assume that the model $\{f(\cdot|\cdot)\}_{\theta \in \ps}$ and the prior $\pi$ satisfy 
\emph{Assumption} \ref{ass:thm1} along with \emph{Propositions} \ref{prp:Laplace_finite}, \ref{prop:Bayes_Club} and \ref{francesi}. Then, as $n \to +\infty$
$$
\epsilon_n = O\big( \varepsilon_{n,2}(\ss,\mu_0) \big),
$$ 
which is the optimal rate, at least when $m=1$. 
\end{prp}

\subsection{Multinomial models}
Consider the case in which the observations, i.e. both the sequence $\{X_i\}_{i \geq 1}$ and the sequence $\{\xi_i\}_{i \geq 1}$, take values in the finite set, say $\{a_1, \dots, a_N\}$. It is easy to check that $\ps$ can be assumed to coincide with the interior of the $(N-1)$-dimensional simplex
$$
\Delta_{N-1} := \left\{\theta = (\theta_1, \dots, \theta_{N-1}) \in [0,1]^{N-1}\ \Big|\ \sum_{i=1}^{N-1} \theta_i \leq 1\right\} 
$$
and
$$
\pi_n(\ud\theta\,|\, x_1, \dots, x_n) = \frac{\left[ \prod_{i=1}^{N} \theta_i^{\nu_{n,i}(x)} \right] \pi(\ud \theta)}{\int_{\Delta_{N-1}} \left[ \prod_{i=1}^{N} t_i^{\nu_{n,i}(x)} \right] \pi(\ud t) }
$$
where $\theta = (\theta_1, \dots, \theta_{N-1})$, $t = (t_1, \dots, t_{N-1})$, $\theta_N := 1 - \sum_{i=1}^{N-1} \theta_i$, $t_N := 1- \sum_{i=1}^{N-1} t_i$ and 
$$
\nu_{n,i}(x) := \sum_{j=1}^n \mathds{1}\{x_j = a_i\} \qquad\qquad i = 1, \dots, N.
$$
Of course, if we put $\ss = \{a_1, \dots, a_N\}$, we can not directly apply Theorem \ref{main_thm1}. Nonetheless, we can resort 
to a reinterpretation of the data, in terms of the frequencies $\nu_{n,i}$, that we now explain, that allows the use of our theorem. We consider
$$
\pi_n^{\ast}(\ud\theta\,|\, p) := \frac{\left[ \prod_{i=1}^{N} \theta_i^{np_i} \right] \pi(\ud \theta)}{\int_{\Delta_{N-1}} \left[ \prod_{i=1}^{N} t_i^{np_i} \right] \pi(\ud t) }
$$
defined for $p = (p_1, \dots, p_{N-1}) \in \Delta_{N-1}$ with the usual proviso that $p_N := 1 - \sum_{i=1}^{N-1} p_i$. Whence,
$$
\pi_n(\ud\theta\,|\, x_1, \dots, x_n) = \pi_n^{\ast}\left(\ud\theta\,\Big|\, \left(\frac{\nu_{n,1}(x)}{n}, \dots, \frac{\nu_{n,N-1}(x)}{n}\right)\right)  \ .
$$ 
The problem of consistency, and the allied question of finding a PCR, can be now reformulated as follows. After fixing $\theta_0 \in \Delta_{N-1}$, we consider the sequence $\{\xi_i\}_{i \geq 1}$ of i.i.d. random variables,
each taking values in $\{a_1, \dots, a_N\}$, with $\P[\xi_1 = a_i] = \theta_{0,i}$, for $i = 1, \dots, N$. An analogous version of Lemma \ref{lm:Wpcr} states that 
$$
\epsilon_n = \E[\Wdue^{(\mathcal P(\ps))}(\pi_n(\ud\theta\,|\, \xi_1, \dots, \xi_n); \delta_{\theta_0}) ]
$$
provides a PCR at $\theta_0$. Now, we reformulate Theorem \ref{main_thm1} as follows. First of all, we have that
$$
\E\left[\left| \left( \frac{\nu_{n,1}(\xi)}{n}, \dots, \frac{\nu_{n,N-1}(\xi)}{n} \right) - \theta_0 \right| \right] \leq \sqrt{\sum_{i=1}^{N-1} \E\left[\left|\frac{\nu_{n,1}(\xi)}{n} - \theta_{0,i} \right|^2\right]} \leq \frac{1}{\sqrt{n}}
$$
replaces the speed of the mean Glivenko-Cantelli convergence. Then, we have that
$$
\kksf(\theta\,|\, \theta_0) = \sum_{i=1}^N \theta_{0,i} \log\left(\frac{\theta_{0,i}}{\theta_{i}}\right)\ .
$$
The relation analogous to that in \eqref{PCR_main1}, which gives a PCR at $\theta_0$, reads as follows
\begin{align} \label{PCR_multinomial}
\epsilon_n &=  3\left( \frac{\int_{\Delta_{N-1}} |\theta-\theta_0|^2 e^{-n \kksf(\theta|\theta_0)} \pi(\ud\theta)}{\int_{\ps} e^{-n \kksf(\theta\,|\,\theta_0)} \pi(\ud\theta)} \right)^{1/2} \\
&\quad+ L_0^{(n)}(\delta_n) \E\left[\left| \left( \frac{\nu_{n,1}(\xi)}{n}, \dots, \frac{\nu_{n,N-1}(\xi)}{n} \right) - \theta_0 \right| \right] \nonumber \\
&\quad+ 2|\theta_0| \P\left[\left| \left( \frac{\nu_{n,1}(\xi)}{n}, \dots, \frac{\nu_{n,N-1}(\xi)}{n} \right)- \theta_0 \right| > \delta_n\right] \nonumber \\
&\quad+ \E\Bigg[\left( 2 \frac{\int_{\Delta_{N-1}} |\theta|^2 \left[ \prod_{i=1}^{N} \theta_i^{\nu_{n,i}(\xi)} \right] \pi(\ud\theta)}{\int_{\Delta_{N-1}} \left[ \prod_{i=1}^{N} \theta_i^{\nu_{n,i}(\xi)} \right]\pi(\ud\theta)} 
\right)^{1/2} \times \nonumber \\
&\quad\quad\times \mathds{1}\left\{\left| \left( \frac{\nu_{n,1}(\xi)}{n}, \dots, \frac{\nu_{n,N-1}(\xi)}{n} \right)- \theta_0 \right| > \delta_n\right\} \Bigg]\nonumber  
\end{align}
where $\{\delta_n\}_{n \geq 1}$ provides a sequence of positive numbers and $L_0^{(n)} $ is defined as follows
$$
L_0^{(n)} := n \sup_{|p - \theta_0| \leq \delta_n} [\mathfrak C_2(\pi_n^{\ast}(\cdot\,|\, p))]^2 \left(
\int_{\Delta_{N-1}} \Big| \nabla_{\theta} \nabla_p \kksf(\theta\,|\, p) \Big|^2 \pi_n^{\ast}(\ud \theta\,|\, p) \right)^{1/2}\ .
$$

We show that the PCR in \eqref{PCR_multinomial} reduces to a simpler expression. Indeed, the first term on the right-hand side of \eqref{PCR_multinomial} is similar to the one already studied in the previous section. By resorting to the same theorems from \citet{Bre(94)}, recalling that the mapping $\theta \mapsto \kksf(\theta\,|\,\theta_0)$ is minimum when 
$\theta=\theta_0$, we get, as $n \to +\infty$,
$$
\left(\frac{\int_{\Delta_{N-1}}  |\theta - \theta_0|^2 e^{-n \kksf(\theta \,|\, \theta_0)} \pi(\ud \theta)}{\int_{\Delta_{N-1}} e^{-n \kksf(\theta\, |\, \theta_0)} \pi(\ud \theta)} \right)^{1/2} = O\big(n^{-1/2}\big)
$$
provided that $\pi$ has full support. As for the second terms on the right-hand side of \eqref{PCR_multinomial}, we have already shown that the expectation is controlled by $1/\sqrt{n}$. Apropos of the constant
$L_0^{(n)}(\delta_n)$, we can easily show that it is bounded, at least whenever $\theta_0$ is fixed in the interior of $\Delta_{N-1}$. In fact, $\delta_n$ can be chosen equal to any positive constant $\delta$ less than the 
distance between $\theta_0$ and the boundary of $\Delta_{N-1}$. In particular, by exploiting the convexity of the mapping $\theta \mapsto \kksf(\theta\,|\, p)$, we can resort to Proposition \ref{francesi}, upon assuming more regularity on the prior distribution $\pi$, in order to obtain $[\mathfrak C_2(\pi_n^{\ast}(\cdot\,|\, p))]^2 \leq C(\delta)/n$, with a positive constant $C(\delta)$ which is independent 
of $p$. Then, by means of a direct computation, under the above conditions on $\theta_0$ and $\delta$ we can show that the integral 
$$
\int_{\Delta_{N-1}} \Big| \nabla_{\theta} \nabla_p \kksf(\theta\,|\, p) \Big|^2 \pi_n^{\ast}(\ud \theta\,|\, p) 
$$
can be bounded uniformly in $n$. To conclude the analysis of the terms on the right-hand side of \eqref{PCR_multinomial}, we only need to exploit the boundedness of $|\theta|$, as $\theta$ varies in $\Delta_{N-1}$, 
to show that the third and the fourth terms are both bounded by a multiple of
$$
\P\left[\left| \left( \frac{\nu_{n,1}(\xi)}{n}, \dots, \frac{\nu_{n,N-1}(\xi)}{n} \right)- \theta_0 \right| > \delta_n\right]\ . 
$$
Thus, if $\theta_0$ is in the interior of $\Delta_{N-1}$ and $\delta_n = \delta$, for the same $\delta$ as above, it is well-known from the theory of large deviations that this probability goes to zero exponentially fast. See \citet[Chapter 2]{DZ(98)} for a detailed account. To conclude, we state a proposition that summarizes the above considerations. 

\begin{prp}
Let $N \geq 2$ be an integer. Let $\pi$ be a prior on $\Delta_{N-1}$. If $\pi$ has a density $q$ (with respect to the Lebesgue measure) such that $q \in \mathrm{C}^1(\overline{\Delta_{N-1}})$ and $q(\theta) = 0$ for any $\theta \in\partial \Delta_{N-1}$, then as $n\rightarrow+\infty$
$$
\epsilon_n = O\big(n^{-1/2}\big),
$$ 
which is the optimal rate.  
\end{prp}

\subsection{Finite-dimensional logistic-Gaussian model}

Consider a class of dominated statistical models specified by density functions of the form

\begin{equation} \label{logistic_fin}
f(x\,|\,\theta) = \frac{e^{\theta \cdot \Gamma_N(x)}}{\int_0^1 e^{\theta \cdot \Gamma_N(x)} \lambda(\ud y)} \qquad x \in \ss, \theta \in \ps
\end{equation}
where, for simplicity, we have fixed $N \in \N$, $\ps = \R^N$, $\ss = [0,1]$, $\ssa = \mathscr{B}([0,1])$, $\lambda = \mathcal{L}^1_{[0,1]}$, that is the one-dimensional Lebesgue measure restricted to [0,1], and
$$
\Gamma_N(x) := \left(\sin \pi x, \sin 2\pi x, \dots, \sin N\pi x \right)\ .
$$
Of course, the expression $\theta \cdot \Gamma_N(x)$ represents a Fourier polynomial and, for sufficiently large $N$, can approximate very well any smooth function, in various norm. This model has been studied 
in connection with the problem of density estimation (\citet{Crain(76a),Crain(76b),Lenk(88),Lenk(91)}), essentially as a toy model. In the following section, we will analyze its infinite dimensional generalization, which is a more flexible statistical model, even if more complex from a mathematical point of view. 

To apply Theorem \ref{new_therem_PCR}, we start by fixing $\theta_0 \in \ps$, so that  $\mu_0(\ud x) = f(x | \theta_0) \ud x$, where $x \mapsto f(x | \theta_0)$ is a continuous and bounded density function on $[0,1]$. Then, we let $\{\xi_i\}_{i \geq 1}$ be a sequence  of independent random variables identically distributed as the probability law $\mu_0$. The model \eqref{logistic_fin} satisfies Definition \ref{def:Reg_Exp_Fam} with $\B = \ps$, $\B^{\ast} = \ps$ (by Riesz's representation theorem) and $\Gamma = \B^{\ast}$, with $ _{\B^{\ast}}\langle \cdot, \cdot \rangle_{\B}$ being identified with the standard (Euclidean) scalar product of $\R^N$. The function $\beta$ coincides with $\Gamma_N(x)$, while $g$ is the identity function. Finally, we have that
$$
M(\theta) = \log\left( \int_0^1 e^{\theta \cdot \Gamma_N(x)} \ud x \right)
$$
which proves to be a convex function, steep, of class $\mathrm{C}^{\infty}(\ps)$ and analytic. Therefore, we have a regular exponential family, in canonical form. As for the prior distribution, besides the multivariate (non-degenerate) Gaussian distribution $\mathcal N(m, Q)$, with $m \in \ps$ and $Q$ being a symmetric and positive-definite $N\times N$ matrix, any other distribution of log-concave form like $\pi(\ud\theta) \propto \exp\{-U(\theta)\} \ud\theta$ fits our assumptions, provided that $U$ is of class $\mathrm{C}^2(\ps)$ and strongly convex. 

Coming back to the application of  Theorem \ref{new_therem_PCR}, we check the validity of the assumptions. First, $|\Gamma_N(x)| \leq \sqrt{N}$ for all $x \in [0,1]$, so that \eqref{first_mom_exp} is in force. Then, \eqref{int_exp_post} and $\int_{\ps} \|\theta\|^{ap} \pi(\ud\theta) < +\infty$ for some $a > 1$ hold, because of the assumptions on the prior distribution. Thus, the bound \eqref{Main_PCR1} provides the desired PCR, so that we proceed further by analyzing the various terms as in Section \ref{sect:statementsNew}. Since
$$
\mathrm{I}(\theta_0) = \mathrm{Hess}[M](\theta_0) = \mathsf{Cov}_{\theta_0}(\Gamma_N(\xi_1))
$$
is strictly positive definite, we can apply Proposition \ref{prp:Laplace_finite}. We conclude that the first term on the right-hand side of \eqref{Main_PCR1} goes to zero as $\frac{1}{\sqrt n}$. Then, the boundedness condition $|\Gamma_N(x)| \leq \sqrt{N}$ for all $x \in [0,1]$ entails that the second and the third terms on the right-hand side of \eqref{Main_PCR1} go to zero exponentially fast by means of classical concentration inequalities, like Bernstein's inequality for instance (\citet{BLM(13),DZ(98)}). Finally, we consider the last term on the right-hand side of \eqref{Main_PCR1}. In particular, by Jensen's inequality
$$
\E[\|\hat{S}_n- S_0\|_{\B}] \leq \left(\E[\|\hat{S}_n- S_0\|_{\B}^2]\right)^{1/2} = O\big(n^{-1/2}\big)\ ,
$$ 
as $n \to +\infty$, where $\hat{S}_n = n^{-1} \sum_{i=1}^n \Gamma_N(\xi_i)$ and 
$$
S_0 = \E_{\theta_0}[\Gamma_N(\xi_1)] = \int_0^1 \Gamma_N(x) \mu_0(\ud x) \ .
$$ 
Then, we apply Proposition \ref{prop:Bayes_Club}. Since $\mathfrak{D}_{\theta}[g]$ coincides with the identity operator, we conclude that 
$$
\int_{\ps} \|\mathfrak{D}_{\theta}[g]\|_{\ast}^2\ \pi_n^{\ast}(\ud\theta\,|\, S)= 1
$$ 
for all $S \in \B$. Whence, 
$$
L_0^{(n)} = n \sup_{S \in \mathcal{U}_{\delta_n}(S_0)} \left\{ \Cpw[\pi_n^{\ast}(\cdot\,|\, S)] \right\}^2 \ .
$$
We conclude our analysis by estimating the weighted Poincar\'e-Wirtinger constant by means of Proposition \ref{francesi}. Indeed, a common feature of these logistic models is that the behavior of the Kullback-Leibler $\kksf(\theta|\theta_0)$ is twofold: it is quadratic as $\theta$ varies around $\theta_0$, while it is linear as $|\theta| \to +\infty$. Thus, the strong Bakry-Emery condition does not apply here, and we resort to the boundedness condition
\begin{equation} \label{Lower_Kullback_francesi}
\theta \cdot \nabla_{\theta} \kksf(\theta\,|\,\theta_0) \geq c|\theta|
\end{equation}
for all $|\theta| \ge R$. To check the validity of this lower bound, we fix for simplicity $\theta_0 = 0$, to get 
$$
\theta \cdot \nabla_{\theta} \kksf(\theta\,|\,\theta_0) = \frac{\int_0^1 \phi(x) e^{\phi(x)} \ud x}{\int_0^1 e^{\phi(x)} \ud x} - \int_0^1 \phi(x) \ud x
$$
with $\phi(x) := \theta \cdot \Gamma_N(x)$. By fixing a unitary vector $\sigma \in S^{N-1}$ and considering $\theta = t\sigma$, the Laplace approximation yields that the above right-hand side is asymptotic to 
$$
t \left[\max_{x \in [0,1]} \sigma \cdot \Gamma_N(x) - \sigma \cdot \int_0^1\Gamma_N(x)\ud x \right]
$$
as $t \to +\infty$. At this stage, the function 
\begin{equation} \label{Fourier_trick}
\sigma \mapsto \max_{x \in [0,1]} \sigma \cdot \Gamma_N(x) - \sigma \cdot \int_0^1\Gamma_N(x)\ud x
\end{equation}
proves to be continuous and non-negative on $S^{N-1}$. The minimum of such a function must be positive, otherwise there would exist $\hat{\sigma} \in S^{N-1}$ for which the map $x \mapsto \hat{\sigma} \cdot \Gamma_N(x)$
turns out to be constant. But this contradicts the linear independence of the Fourier basis $\Gamma_N$, yielding that the function in \eqref{Fourier_trick} must be strictly positive. This fact validates \eqref{Lower_Kullback_francesi}. Thus, by point (1) of Proposition \ref{francesi}, and the square of the weighted Poincar\'e-Wirtinger constant is asymptotic to $1/n$.  All the above considerations can be summarized in the following proposition.
\begin{prp}  
Let $\pi(\ud\theta) \propto \exp\{-U(\theta)\} \ud\theta$ be any prior on $\ps = \R^N$, with $U$ is of class $\mathrm{C}^2(\ps)$ and stongly convex. Then, 
for the finite-dimensional logistic-Gaussian model as in \eqref{logistic_fin} with $N \in \N$, $\ss = [0,1]$, $\ssa = \mathscr{B}([0,1])$ and $\lambda = \mathcal{L}^1_{[0,1]}$, the PCR $\epsilon_n$ satisfies
$$
\epsilon_n = O\big(n^{-1/2}\big)
$$
as $n \to +\infty$, which is the optimal one. 
\end{prp}

\subsection{Infinite-dimensional logistic-Gaussian model}\label{Inf_Exp_Fam}

Consider a class of dominated statistical models specified by density functions of the form
\begin{equation} \label{logistic_inf}
f(x\,|\,\theta) = \frac{e^{\theta(x)}}{\int_0^1 e^{\theta(y)} \lambda(\ud y)} \qquad x \in \ss, \theta \in \ps
\end{equation}
where we have fixed $\ss = [0,1]$, $\ssa = \mathscr{B}([0,1])$, $\lambda = \mathcal{L}^1_{[0,1]}$, that is the one-dimensional Lebesgue measure restricted to [0,1]. As for the
parameter space $\ps$, we set
\begin{equation} \label{Theta_Sobolev}
\ps = \Hunoz := \{\phi \in \Huno\ |\ \phi(0)=0\}
\end{equation}
thought of as an infinite-dimensional Hilbert space endowed with scalar product 
\begin{equation} \label{scalarP_Hunoz}
\langle \phi, \psi \rangle = \int_0^1 \phi'(z) \psi'(z) \ud z
\end{equation}
and norm
\begin{displaymath}
\|\phi\|_{\Hunoz} := \left(\int_0^1 [\phi'(z)]^2 \ud z\right)^{1/2}\ .
\end{displaymath}
Here, the well-known Sobolev embedding theorem (\citet{Mazya(11)}) states that $\Hunoz$ is continuously embedded in $\mathrm{C}^0[0,1]$, and therefore the above notations $\theta(x)$ and $\phi(0)$ are referred to the continuous representatives of $\theta$ and $\phi$, respectively. 

The infinite-dimensional logistic-Gaussian model is typically considered in connection with the fundamental problem of  density estimation (\citet{Crain(76a),Crain(76b),Lenk(88),Lenk(91)}). Under the assumption that the prior is a Gaussian measure, Bayesian consistency is investigated in \citet{TL(07)}, whereas PCRs are provided in \citet{GiNi(11),RivRou(12),Scri(06)} and \citet{vdVvZ(08)}. These results consider the set $\ps$ to be the space of all density functions on $[0,1]$, typically endowed with the total variation distance, the Hellinger distance, some $\mathrm{L}^p$ norm or the Kullback-Leibler divergence. Our approach to PCRs relies on the choice  \eqref{Theta_Sobolev}, so that our PCRs refers to Definition \ref{def:consistency} with $\ud_{\ps}$ equal to the $\Hunoz$-norm. This metric is generally stronger, since a Sobolev norm is, for suitable exponents, grater that the $\mathrm{L}^r$ norm considered in \citet{GiNi(11)} and, in turn, greater than the (squared) Hellinger distance, as proved in \citet[Lemma A.1]{Scri(06)}. In connection with the statistical model \eqref{logistic_inf}, the work of \citet{Fuk(09)} provides an implicit Riemannian structure on the space of densities which is modeled on the metric of the underlying space $\ps$, that is the Riemannian distance between two densities $f(\cdot | \theta_1)$ and $f(\cdot | \theta_2)$ turns out to be locally equivalent to $\|\theta_1 - \theta_2\|_{\Hunoz}$. Another (geometrical) view
of the set $\{f(\cdot\,|\,\theta)\}_{\theta \in \ps}$, which is simply thought as a differential manifold, is provided in \citet{PisRo(99)}.

We provide PCRs for the model \eqref{logistic_inf} on the basis of Theorem \ref{new_therem_PCR}. We start by fixing $\theta_0 \in \ps$, with $\ps$ being the same as in \eqref{Theta_Sobolev}. Whence, $\mu_0(\ud x) = f(x | \theta_0) \ud x$, where $x \mapsto f(x | \theta_0)$ is a continuous and bounded density function on $[0,1]$. Then, we let $\{\xi_i\}_{i \geq 1}$ be a sequence 
of independent random variables identically distributed with probability law $\mu_0$. At this stage, we notice that the model \eqref{logistic_inf} satisfies Definition \ref{def:Reg_Exp_Fam} with $\B = \ps$, $\B^{\ast} = \ps$ (by Riesz's representation theorem) and $\Gamma = \B^{\ast}$. For completeness, we specify that also the pairing $ _{\B^{\ast}}\langle \cdot, \cdot \rangle_{\B}$ is identified with the scalar product $\langle\cdot, \cdot\rangle$ as in \eqref{scalarP_Hunoz}, again by Riesz's representation theorem.
In this setting, we deduce that the function $\beta$ in Definition \ref{def:Reg_Exp_Fam} coincides with the Riesz representative of the $\delta_x$ functional, for any $x \in [0,1]$, that is $\beta_x(z) := z\ind_{[0,x]}(z) + x\ind_{(x, 1]}(z)$ for $z \in [0,1]$, since, for any $\theta \in \ps$, $\theta(x) = \langle \theta, \beta_x\rangle$. Lastly, we fix $g$ as the identity map on $\ps$, so that \eqref{mu_n_exp} is satisfied and 
\begin{equation} 
M(\theta) = \log \int_0^1 e^{\theta(y)} \ud y .
\end{equation}
As for the prior $\pi$, we assume that it is a Gaussian measure on $\ps$, with mean $m \in \ps$ and covariance operator $Q : \ps \to \ps$. We recall that $Q$ is a trace operator with eigenvalues $\{\lambda_k\}_{k \geq 0}$ that satisfy $\sum_{k=0}^{\infty} \lambda_k < +\infty$. See \citet{dP(14),DaPrato(06)}, and references therein, for a review on Gaussian measures on Hilbert spaces. 

Now, we check the validity of the assumptions of Theorem \ref{new_therem_PCR}. First, we have that 
\begin{equation} \label{norm_beta}
\|\beta_x\|_{\B} = \left(\int_0^1 \ind_{[0,x]}(z) \ud z\right)^{1/2} = \sqrt{x} \in [0,1]
\end{equation}
yielding that \eqref{first_mom_exp} is trivially satisfied. In particular, the element $S_0$ is given by
$$
z \mapsto S_0(z) = \int_0^1 \beta_x(z) \mu_0(\ud x) = z\mu_0([z,1]) + \int_0^z x \mu_0(\ud x) \in \ps
$$
and $\hat{S}_n = n^{-1}\sum_{i=1}^n \beta_{\xi_i}$. We recall that $\hat{S}_n \rightarrow S_0$ as $n \rightarrow +\infty$, in both $\P$-a.s. and $\mathrm{L}^2$ sense, by the Laws of Large Numbers
in Hilbert spaces (\citet[Corollary 7.10]{LT(91)}). Now, we observe that \eqref{int_exp_post} boils down to write that
$$
\int_{\ps} \exp\{n\ _{\B^{\ast}}\!\langle g(\tau), b\rangle_{\B} \} \pi(\ud\tau) = \exp\left\{n\langle m, b\rangle + \frac{n^2}{2}\langle Q[b], b\rangle \right\} < +\infty
$$ 
for all $n \in \N$ and $b \in \B = \ps$. See \citet[Proposition 1.15]{DaPrato(06)}. Then, condition iii) of Theorem \ref{new_therem_PCR} is trivially satisfied. Finally, with regards to iv), we mention that any sequence $\delta_n \sim n^{-q}$ with $q \in [0,1/2)$, as $n \rightarrow +\infty$, is valid as far as we verify the validity of \eqref{Lipschitz_kernel},
as we will do just below. After these preliminaries, we start analyzing the four terms on the right-hand side of \eqref{Main_PCR1}.

As for the first term on the right-hand side of \eqref{Main_PCR1}, we study \eqref{Gibbs}. We observe that
\begin{align*}  
\mathrm{K}(\theta\, |\, \theta_0) &= \int_0^1 [\theta_0(y) - M(\theta_0) - \theta(y) + M(\theta)] f(y\,|\,\theta_0) \ud y \\
&= \int_0^1 \langle \theta_0 - \theta, \beta_y\rangle f(y\,|\,\theta_0) \ud y + M(\theta) - M(\theta_0)  \\
&= \langle \mathrm D_{\theta_0} M, \theta_0 - \theta\rangle + M(\theta) - M(\theta_0) \\
&= \frac 12 \langle \mathrm{Hess}_{\theta_0}[M][\theta_0 - \theta], \theta_0 - \theta\rangle + o(\|\theta_0 - \theta\|^2) \qquad (\text{as}\  \theta \to \theta_0)
\end{align*}
where $\mathrm D_{\theta_0} M$ represents the (Riesz representative of) the Fr\'echet differential of $M : \ps \rightarrow \R$ at $\theta_0$, while $\mathrm{Hess}_{\theta_0}[M] = \mathrm{I}(\theta_0)$ stands for the Hessian operator of $M$ 
at $\theta_0$, which coincides with the Fisher information operator $\mathrm{I}(\theta_0)$. In particular, in the last identity we have used the Taylor expansion of $M$ around $\theta_0$. In view of a more concrete characterization of
$\mathrm D_{\theta_0} M$ and $\mathrm{I}(\theta_0)$, we write that
\begin{align*}
&M(\theta_0 + h) - M(\theta_0)\\
&\quad = \int_{0}^1 h(y)\mu_0(\ud y) + \frac 12\left[\int_{0}^1 h^2(y)\mu_0(\ud y) - \left(\int_{0}^1 h(y)\mu_0(\ud y)\right)^2\right] + \mathcal R(h; \mu_0)
\end{align*}
where $\|\mathcal R(h; \mu_0)\| \leq C(\mu_0) \|h\|_{\ps}^3$ for $\|h\|_{\ps} \leq 1$, with some suitable constant $C(\mu_0)$ depending solely on $\mu_0$. In particular, a straightforward integration by parts shows that 
$$
\int_{0}^1 h(y)\mu_0(\ud y) = \int_{0}^1 h'(y)\Phi'_0(y) \ud y = \langle h, \Phi_0 \rangle
$$
with $\Phi_0(y) := \int_{0}^y [1 - F_0(z)] \ud z$ and $F_0(z) := \mu_0([0,z])$. Whence, $\Phi_0 = \mathrm D_{\theta_0} M$, by means of Riesz's representation. Moreover, with the same technique, we obtain
\begin{equation} \label{Fisher_logistic}
\mathrm{Hess}_{\theta_0}[M][h](y) = 2\int_0^y h(z) [1 - F_0(z)] \ud z - \langle h, \Phi_0 \rangle \Phi_0(y)
\end{equation}
for any $y \in [0,1]$ and $h \in \ps$. Tthe above left-hand side should be read as follows: first, the operator $\mathrm{Hess}_{\theta_0}[M]$, applied to $h \in \ps$, gives a new element of $\ps$, called $\mathrm{Hess}_{\theta_0}[M][h]$; second, this new object, as a continuous function evaluated at $y$, coincides with the right-hand side. 
Finally, integration by parts entails that 
$$
\int_{0}^1 h^2(y)\mu_0(\ud y) - \left(\int_{0}^1 h(y)\mu_0(\ud y)\right)^2 = \langle h, \mathrm{Hess}_{\theta_0}[M][h]\rangle
$$
for any $h \in \ps$. The way is now paved for the application of Proposition \ref{prp:Laplace_infinite} and Remark \ref{rmk:Carlen}. As first step, we check that the operator in \eqref{Fisher_logistic}, from $\ps$ to itself, is compact. As for the term $\langle h, \Phi_0 \rangle \Phi_0$, it defines a finite-rank operator, which is of course compact. As for the term $2\int_0^y h(z) [1 - F_0(z)] \ud z$, it is enough to pick a bounded sequence, say $\{h_n\}_{n \geq 1}$, in $\ps$, and study the sequence $\{\Psi_n\}_{n \geq 1}$ given by $\Psi_n(y) := 2\int_0^y h_n(z) [1 - F_0(z)] \ud z$.  Now, from the well-known properties of weak topologies of separable Hilbert spaces, we can extract a subsequence $\{h_{n_j}\}_{j \geq 1}$, which converges weakly to some $h_{\ast} \in \ps$. Whence, $h_{n_j}$ converges uniformly
(i.e. in the strong topology of $\mathrm{C}^0[0,1]$) to $h_{\ast}$, by the Rellich-Kondrachov embedding theorem. Consequently, it is trivial to get that the sequence $\{\Psi_{n_j}\}_{j \geq 1}$ converges strongly in $\ps$ to $\Psi_{\ast}(y) := 2\int_0^y h_{\ast}(z) [1 - F_0(z)] \ud z$ since
$$
\|\Psi_{n_j} - \Psi_{\ast}\|_{\ps}^2 = 4 \int_0^1 |h_{n_j}(z) - h_{\ast}(z)|^2 [1 - F_0(z)]^2 \ud z \leq 4 \| h_{n_j} - h_{\ast} \|_{\infty}^2 \rightarrow 0
$$
as $j \rightarrow +\infty$. This proves that the operator in \eqref{Fisher_logistic}, from $\ps$ to itself, is a compact operator, even if it is not self-adjoint. Then, we resort to Remark \ref{rmk:Carlen}, noticing that 
\begin{align*}
&\int_0^1 h^2(y)\mu_0(\ud y) - \left(\int_{0}^1 h(y)\mu_0(\ud y)\right)^2 \\
&= \int_{0}^1 \left[h(x) - \int_{0}^1 h(y)\mu_0(\ud y)\right]^2 f(x\,|\,\theta_0) \ud x\\
&\geq \exp\{- \text{osc}(\theta_0)\} \int_{0}^1 \left[h(x) - \int_{0}^1 h(y)\mu_0(\ud y)\right]^2 \ud x\\
&\geq \exp\{- \text{osc}(\theta_0)\} \int_{0}^1 \left[h(x) - \int_{0}^1 h(y)\ud y\right]^2 \ud x\
\end{align*}
where $\text{osc}(\theta_0) := \max_{x \in [0,1]} \theta_0(x) - \min_{x \in [0,1]} \theta_0(x)$ is the oscillation. Therefore, we can set
$$
\mathrm{I}^{\dagger}[h] := \exp\{- \text{osc}(\theta_0)\} \mathrm{Hess}_{0}[M][h]
$$
where $\mathrm{Hess}_{0}[M][h]$ is defined by \eqref{Fisher_logistic} with $\theta_0 \equiv 0$, to re-write the above relations  as 
$$
\langle h, \mathrm{Hess}_{\theta_0}[M][h]\rangle \geq \langle h, \mathrm{I}^{\dagger}[h] \rangle
$$
or simply as $\mathrm{Hess}_{\theta_0}[M] \geq \mathrm{I}^{\dagger}$. By means of the above argument, $\mathrm{I}^{\dagger}$, as a linear operator from $\ps$ to itself, is again compact, but not
self-adjoint. By a straightforward integration by part, we find that a self-adjointized version of $\mathrm{I}^{\dagger}$ is given by
$$
\mathrm{I}^{\ast}[h](x) := \exp\{- \text{osc}(\theta_0)\} \left\{\int_0^1 \beta_x(y) h(y) \ud y - \left(x - \frac{x^2}{2}\right) \int_0^1 h(y) \ud y\right\}
$$
with $x \in [0,1]$ and $h \in \ps$. The relation $\mathrm{Hess}_{\theta_0}[M] \geq \mathrm{I}^{\ast}$ is, of course, still in force. We can now invoke the spectral theorem for compact, self-adjoint operators on separable Hilbert spaces to deduce the existence of a Fourier basis (complete orthonormal system) $\{\mathbf e_k\}_{k \geq 1}$ for $\ps$ which diagonalizes $\mathrm{I}^{\ast}$. 
With reference to \eqref{eigenvalues_logistic2}, we call $\{\gamma_k^{\ast}\}_{k \geq 1}$ the sequence of the relative eigenvalues, for which we have that $\gamma_k^{\ast} \rightarrow 0$ as $k \rightarrow +\infty$, again by the spectral theorem. An explicit derivation of $\mathbf e_k$ could be drawn from the following integral-differential Cauchy problem
$$
\begin{cases}
-\mathbf e_k(x) + \displaystyle{\int_{0}^1} \mathbf e_k(y)\ud y = \gamma_k^{\ast} \mathbf e^{''}_k(x)\ ,  \qquad x \in [0,1] \\
\mathbf e_k(0) = \mathbf e^{'}_k(0) = 0
\end{cases}
$$ 
which is obtained by differentiating twice the relation $\mathrm{I}^{\ast}[\mathbf e_k] = \gamma_k^{\ast} \mathbf e_k$. Explicit solutions are
\begin{equation} \label{Fourier_Basis}
\mathbf e_k(x) = \frac{\sqrt{2}}{k\pi}(1 - \cos(k\pi x)) \qquad\qquad \gamma_k^{\ast} = \frac{\exp\{- \text{osc}(\theta_0)\}}{(k\pi)^2} 
\end{equation}
with $x \in [0,1]$ and $k\in \N$. After having fixed the Fourier basis $\{\mathbf e_k\}_{k \geq 1}$, we can further specify the prior distributions in terms of the probability laws of the random elements $\Xi$, with values in $\ps$, of the form (Karhunen-Lo\`eve representation)
$$
\Xi = \sum_{k=1}^{\infty} Z_k \mathbf e_k\ . 
$$
Here, $\{Z_k\}_{k \geq 1}$ is a sequence of independent real-valued random variables with $Z_k \sim \mathcal N(m_k, \lambda_k)$, for suitable sequences $m := \{m_k\}_{k \geq 1} \subset \R$ and 
$\{\lambda_k\}_{k \geq 1} \subset (0,+\infty)$ with $\{m_k\}_{k \geq 1} \in \ell^2$ and $\{\lambda_k\}_{k \geq 1} \in \ell^1$. Thus, if  $\pi(B) := \P\left[\ \Xi \in B \right]$ for any $B \in \mathscr{B}(\ps)$, it is straightforward to check that $\pi$ is a Gaussian measure on $(\ps,\mathscr{B}(\ps))$ with mean $m$ and covariance operator $Q$ satisfying 
$Q[\mathbf e_k] = \lambda_k \mathbf e_k$. Whence, \eqref{eigenvalues_logistic2} is verified. To justify the validity of \eqref{Laplace2_ast}, we check the 
remaining assumptions of Proposition \ref{prp:Laplace_infinite}. First, it is trivial to check that $\theta \mapsto \kksf(\theta|\theta_0)$ belongs to $\mathrm{C}^{\infty}(\ps)$, so that we can put $q=1$. 
Then, we consider points i)--iv). For simplicity, we again fix $\theta_0 \equiv 0$, with no real loss of generality. We start with the definition of the space $\Kbb$, expressed as the closure of $\ps$ with 
respect to the norm
\begin{equation} \label{Norm_HilbertK}
\|\theta\|_{\Kbb} := \sup_{\substack{\psi \in \ps \\ \|\psi\|_{\ps} \leq 1}} \int_0^1 \left[\theta(x) - \int_0^1 \theta(y)\ud y \right] \psi(x) \ud x
\end{equation}
which represents, plainly speaking, a dual Sobolev norm of the function $x \mapsto  \theta(x) - \int_0^1 \theta(y)\ud y$. The embedding $\ps \subset \Kbb$, with dense and continuous inclusion, follows from the Poincar\'e-Wirtinger inequality. Then, we notice that the function
\begin{align*}
\theta \mapsto \kksf(\theta| \theta_0) &= \log\left(\int_0^1 e^{\theta(y)} \ud y\right) - \int_0^1 \theta(y) \ud y \\
&= \log\left(\int_0^1 \exp\left\{\theta(x) - \int_0^1 \theta(y)\ud y\right\} \ud x\right)
\end{align*}
has two different behaviors according on whether the norm of $\theta$ is small or large. To be more precise, we fix $\sigma \in \ps$ with $\|\sigma\|_{\ps} =1$ and then we set $\theta = t\sigma$ for any $t\in (0,+\infty)$.
In particular, as $t \to 0$, a straightforward argument based on Taylor expansions of the exponential and the logarithmic functions shows that  
$$
\kksf(t\sigma| \theta_0) = \frac{t^2}{2} \int_0^1 \left[ \sigma(x) - \int_0^1 \sigma(y)\ud y\right]^2 \ud x + o(t^2)\ .
$$
On the other hand, as $t \to +\infty$, by means of a direct application of the Laplace method of approximation (\citet[Theorem 1.II]{Won(01)}), we obtain the following expansion
$$
\kksf(t\sigma| \theta_0) \sim t \max_{x \in [0,1]} \left(\sigma(x) - \int_0^1 \sigma(y)\ud y\right)_+
$$
with $(a)_+ := \max\{a, 0\}$. Upon denoting by $\mathrm{H}^{-1}_{\ast}(0,1)$ the dual space of $\ps$, we can exploit that $\mathrm{L}^1(0,1) \subset \mathrm{H}^{-1}_{\ast}(0,1)$, with continuous 
dense embedding, to obtain that
$$
\max_{x \in [0,1]} \left(\sigma(x) - \int_0^1 \sigma(y)\ud y\right)_+ \geq \frac 12 \int_0^1 \left| \sigma(x) - \int_0^1 \sigma(y)\ud y \right| \ud x \gtrsim \|\sigma\|_{\Kbb} \ .
$$
Therefore, \eqref{lower_kullback}-\eqref{lower_fisher} are fulfilled with the above choice of the space $\Kbb$, and some $\phi : [0,+\infty) \to [0,+\infty)$ which behaves quadratically for small arguments
and linearly for large arguments, like $\phi(x) = x^2 \ind_{[0,1]}(x) + x\ind_{(1,+\infty)}(x)$. Then, the choice of $q=1$ entails that $r \in (1, \frac 32)$. 
Further insights on inequalities like  
\eqref{lower_kullback} can be found in \citet{BMRS(21)}, while 
 properties of homogeneous spaces like $\Kbb$ have been recently investigated in \citet{BGV(21)}. As for the validity of the interpolation inequality \eqref{interpolation}, we can fix, for example, 
 $r= 4/3$, $s=4$ and start from the following specific version of the Gagliardo-Nirenberg interpolation inequality
$$
\|f\|_{\mathrm{H}^1(0,1)} \lesssim \|f\|_{\mathrm{L}^2(0,1)}^{7/8} \|f\|_{\mathrm{H}^8(0,1)}^{1/8}
$$
where $\mathrm{H}^m(0,1)$ denotes the standard (Hilbertian) Sobolev space of order $m$ \citep[Corollary 5.1]{BM(18)}. Applying this inequality to $f(x) = \theta(x) - \int_0^1 \theta(y)\ud y$, 
we get 
\begin{equation} \label{GaNi}
\|\theta\|_{\ps} \lesssim \left\|\theta - \int_0^1 \theta(y)\ud y\right\|_{\mathrm{L}^2(0,1)}^{7/8} \cdot \left\|\theta - \int_0^1 \theta(y)\ud y\right\|_{\mathrm{H}^8(0,1)}^{1/8}
\end{equation}
for all $\theta \in \ps$ such that $\frac{\ud^8}{\ud x^8} \theta(x) \in \mathrm{L}^2(0,1)$. Now, we define the Hilbert space $\Vbb$ as the subspace of $\ps$ formed by those $\theta \in \ps$ such that 
$\frac{\ud^8}{\ud x^8} \theta(x) \in \mathrm{L}^2(0,1)$, with the norm
\begin{equation} \label{Norm_HilbertV}
\|\theta\|_{\Vbb} := \left\|\theta - \int_0^1 \theta(y)\ud y\right\|_{\mathrm{H}^8(0,1)}\ .
\end{equation}
The inclusion $\Vbb \subset \ps$ with continuous and dense embedding follows by means of the usual Sobolev embedding theorem \citep{Mazya(11)}. At this stage, we make use of the other  
specific version of the Gagliardo-Nirenberg interpolation inequality given by
$$
\|f\|_{\mathrm{H}^2(0,1)} \lesssim \|f\|_{\mathrm{H}^1(0,1)}^{6/7} \|f\|_{\mathrm{H}^8(0,1)}^{1/7}
$$
to deduce that 
$$
\|u\|_{\mathrm{L}^2(0,1)} \lesssim \|u\|_{\Kbb}^{6/7} \|u\|_{\Vbb}^{1/7}
$$
holds for any $u \in \mathrm{C}^{\infty}_c(0,1)$ with $\int_0^1 u(x)\ud x = 0$. By combining this inequality with \eqref{GaNi}, we finally deduce \eqref{interpolation} with $r= 4/3$ and $s=4$. To guarantee that $\pi(\Vbb) = 1$, we can resort to the standard Kolmogorov three-series criterion to obtain that $\P\left[\ \Xi \in \Vbb \right] =1$ provided that
$m_k = O(k^{-8-\delta})$ and $\lambda_k = O(k^{-16-\delta})$ as $k \to +\infty$, for some $\delta > 0$. Finally, since we have $\|\mathbf e_k\|_{\mathbb V} = O(k^7)$ as $k \to +\infty$, then
\begin{align*}
\int_{\Vbb} e^{t \|\theta\|_{\Vbb}} \pi(\ud\theta) = \E\left[e^{t \|\Xi\|_{\Vbb}} \right] &\leq \E\left[\exp\left\{t \sum_{k=1}^{\infty} |Z_k| \cdot \|\mathbf e_k\|_{\Vbb} \right\}\right] \\
&\leq \exp\left\{t \sum_{k=1}^{\infty} |m_k| \cdot \|\mathbf e_k\|_{\Vbb} + \frac{t^2}{2} \sum_{k=1}^{\infty} \lambda_k \|\mathbf e_k\|^2_{\Vbb}  \right\} < +\infty
\end{align*} 
holds for any $t > 0$.

By proceeding with the analysis of the other terms on the right-hand side of \eqref{Main_PCR1}, we observe that the boundedness condition \eqref{norm_beta} entails a direct application of results in \citet{Pin(86)} and \citet{Yur(70)}, yielding that
$$
\E[\|\hat{S}_n- S_0\|_{\B}] \leq \left(\E[\|\hat{S}_n- S_0\|_{\B}^2]\right)^{1/2} = O\big(n^{-1/2}\big)
$$ 
and in addition that, for any sequence $\{\delta_n\}_{n \geq 1}$ such that $\delta_n \sim n^{-q}$ with $q \in [0,1/2)$,
$$
\P\left[\hat{S}_n \not\in \mathcal{U}_{\delta_n}(S_0) \right] \leq 2\exp\{-C n^{1-2q}\}
$$
for a positive constant $C$ that depends only on $\mu_0$. It remains to deal with the asymptotic behavior of $L_0^{(n)}$ by combining Propositions \ref{prop:Bayes_Club} and Proposition \ref{prp:Poincare_infinite}. 
Here, we exploit once again the fact that $g$ coincides with identity function, so that 
$$
\int_{\ps} \|\mathfrak{D}_{\theta}[g]\|_{\ast}^2\ \pi_n^{\ast}(\ud\theta\,|\, S)= 1
$$ 
for all $S \in \ps$. Whence, 
$$
L_0^{(n)} = n \sup_{S \in \mathcal{U}_{\delta_n}(S_0)} \left\{ \Cpw^{(M)}[\pi_n^{\ast}(\cdot\,|\, S)] \right\}^2, 
$$
where we have indicated our preference for a weighted Poincar\'e-Wirtinger constant, with respect to the Malliavin derivative. Indeed, we can argue as in the finite-dimensional setting, exploiting the key observation that the Kullback-Leibler divergence $\kksf(\theta|\theta_0)$ behaves quadratically if $\theta$ varies around $\theta_0$, while it is linear as $\|\theta\| \to +\infty$. To be more precise, we can use the same arguments developed above to show that 
the choice $\mathrm{G}_0 = \mathrm{I}^*$ fits the requirements of Proposition \ref{prp:Poincare_infinite}. Thus, the eigenfunctions $\{\mathbf e_k\}_{k \geq 1}$ are the same as above, and
$\eta_k = \gamma_k^{\ast}$. In order to exploit point (1) of Proposition \ref{prp:Poincare_infinite}, we can mimic the same arguments already used in the previous section to prove \eqref{Lower_Kullback_francesi}. Actually, it works in the same way, with the sole difference that the function in \eqref{Fourier_trick} is now replaced by 
\begin{equation} \label{compact_Malliavin}
\sigma \mapsto \max_{x \in [0,1]} Q^{1/2}[\sigma](x) - \int_0^1 Q^{1/2}[\sigma](x)\ud x
\end{equation}
with $\|\sigma\|_{\ps} = 1$, because of the fact that the gradient is replaced by the Malliavin derivative. See \citet[Section 2.3]{dP(14)}. Since $Q^{1/2}$ is a compact operator, the image of the bounded set 
$\{\sigma \in \ps\ |\ \|\sigma\|_{\ps} = 1\}$ through $Q^{1/2}$ is sequentially compact. Thus, the infimum of the function in \eqref{compact_Malliavin} cannot be equal to zero. Finally, with the application of 
\eqref{tesi_francesi_infinite}, which provides the rate of the weighted Poincar\'e-Wirtinger constant, the discussion is completed. To conclude our analysis, we state a proposition that summarizes all the above considerations. 
\begin{prp} \label{prop:PCR_logistic}
In connection with the model \eqref{logistic_inf}, let $\ss = [0,1]$ and $\ps = \mathrm{H}^1_{\ast}(0,1)$. Let $\theta_0 \in \ps$ be fixed. Assume that $\pi = \mathcal N(m,Q)$ with $m \in \ps$ and $Q$ a non-degenerate trace-class operator satisfying \eqref{eigenvalues_logistic}. Fix the eigenfunctions $\{\mathbf e_k\}_{k \geq 1}$ and the spaces $\Kbb$ and $\Vbb$ as in \eqref{Fourier_Basis}, \eqref{Norm_HilbertK} and
\eqref{Norm_HilbertV}, respectively. Finally, set $\gamma^{\ast}_k$ as in \eqref{Fourier_Basis}, 
$\eta_k = \gamma^{\ast}_k$ and $\omega_k$ according to the Fourier representation $\theta_0 - m = \sum_{k=1}^{\infty} \omega_k \mathbf{e}_k$. 
If $m_k = O(k^{-8-\delta})$ and $\lambda_k = O(k^{-16-\delta})$ as $k \to +\infty$, for some $\delta > 0$, then  
points i)-iv) of Proposition \ref{prp:Laplace_infinite} are valid, along with the assumptions of point (1) of Proposition \ref{prp:Poincare_infinite}. In conclusion,
it holds
$$
\epsilon_n = O\left(\sqrt{\sum_{k=1}^{\infty} \frac{\lambda_k}{n \lambda_k \gamma_k^{\ast} + 1}} + \sqrt{\sum_{k=1}^{\infty} \frac{\omega_k^2}{(n \lambda_k \gamma_k^{\ast} + 1)^2}} +  
\sqrt{n} \max_{k \in \N} \left\{\frac{\lambda_k}{n\lambda_k \eta_k + 1} \right\}\right)\ . 
$$
\end{prp}

To provide some hints on the optimality of our PCRs, it is useful to recall the discussion at the end of Section 3.1. At least in the simpler case when $m = \theta_0$, the above rate has the form 
$O\big(n^{-\frac{a-1}{2(a + 3)}}\big)$ when $\lambda_k = O(k^{-(1+a)})$. The parameter $a$ can 
be interpreted as a smoothness parameter, in the sense that it measures the analytical regularity of the trajectories of the prior $\pi$. By way of example, supposing $m_k = 0$ for all $k$ for simplicity, 
a precise statement is as follows: if $a > 1$ and $\varepsilon \in \left(0, \frac{a-1}{2}\right)$, then the trajectories of the random process $\Xi$ belong to $\mathrm H^{1+\varepsilon}(0,1)$ almost surely. We notice that our rate is just slightly slower than the standard rate 
$n^{-\frac{\alpha}{2\alpha + 1}}$ which is proved in \citet{GiNi(11),RivRou(12),Scri(06)}, where $\alpha$ is characterized by the fact that the random process $\Xi$ belongs to $\mathrm H^{\alpha}(0,1)$ almost surely.
This slight discrepancy makes sense since our reference norm (i.e., the Sobolev norm of $\mathrm{H}^1_{\ast}$) is larger than any $\mathrm{L}^p$ norm, for any $p \in [1,+\infty]$. To the best of our knowledge, 
our rate does not admit a fair comparison with any other known rate of consistency, neither Bayesian nor classical, because of the different choice of the loss function. The only fair comparison could be made with the rates 
obtained in \citet{SFGHK(17)}, which are nonetheless relative to distinguished classical estimators (see, in particular, Theorem 7, point $ii)$ therein). Since these classical rates are slower than $n^{-1/3}$, 
we notice, in support of the optimality of our approach, that our rate is: i) arbitrarily close to the optimal (parametric) rate $n^{-1/2}$ if $a \to +\infty$; ii) faster than $n^{-1/3}$ as soon as $a>9$, a condition which is surely met in the framework presented in Proposition \ref{prop:PCR_logistic}, where $a = 15 + \delta$. Hence, a Bayesian estimator, that shares our PCR as rate of consistency, performs better that the minimum-distance estimator proposed in \citet{SFGHK(17)}.


\subsection{Infinite-dimensional linear regression} \label{sect:linear_regression}

Consider a statistical model that arises from the popular linear regression. The observed data are the collection of pairs $(u_1,v_1), \dots, (u_n, v_n)$, such that: i) the $u_{i}$ vary in an interval $[a,b] \subset \reals$, and are modeled as i.i.d. random variables, say $U_1, \dots, U_n$, with a known distribution, say $\varpi(\ud u) = h(u)\ud u$, on $([a,b], \Bcr([a,b]))$; the $v_{i}$'s vary in $\reals$, and are modeled as i.i.d. random variables $V_1, \dots, V_n$. The $V_i$'s are stochastically dependent of the $U_i$'s according 
to the relation
\begin{equation} \label{lin_regression}
V_i = \theta(U_i) + E_i\ , \qquad\qquad i=1, \dots, n\ ,
\end{equation}
where $E_1, \dots, E_n$ are i.i.d. random variables with Normal $\mathcal N(0,\sigma^2)$ distribution, while $\theta : [a,b] \to \reals$ is an unknown continuous function. Assuming for simplicity that $\sigma^2 > 0$ is known,
the statistical model is characterized by probability densities $f(\cdot|\theta)$ on $[a,b]\times \reals$, with respect to the Lebesgue measure, given by
\begin{equation} \label{density_regression}
f(x|\theta) = f((u,v)|\theta) = \frac{1}{\sqrt{2\pi\sigma^2}} \exp\left\{- \frac{[\theta(u) - v]^2}{2\sigma^2}\right\} h(u) \ .
\end{equation}
The space $\ps$ is chosen, as in the previous section, as a Sobolev space $H^s(a,b)$ with $s > 1/2$, which is continuously embedded in $C^0[a,b]$. Whence, upon fixing $\theta_0 \in \ps$, 
$$
\mu_0(\ud u \ud v) = \frac{1}{\sqrt{2\pi\sigma^2}} \exp\left\{- \frac{[\theta_0(u) - v]^2}{2\sigma^2}\right\} h(u) \ud u \ud v\ .
$$
On the other hand, from the Bayesian point of view, upon fixing a prior distribution $\pi$ on $(\ps, \psa)$ and resorting to the Bayes formula, the posterior takes on the form
\begin{align*}
\pi_n(\ud \theta|x_1, \dots, x_n) &= \pi_n(\ud \theta|(u_1, v_1) \dots, (u_n, v_n)) \\
&= \frac{\exp\left\{- \frac{1}{2\sigma^2} \sum_{i=1}^n [\theta(u_i) - v_i]^2\right\} \pi(\ud\theta)}{\int_{\ps} \exp\left\{- \frac{1}{2\sigma^2} \sum_{i=1}^n [\tau(u_i) - v_i]^2\right\} \pi(\ud\tau)}\ .
\end{align*}
Whence, for any probability measure $\gamma \in \mathcal P_2([a,b]\times\reals)$, we can write the following
$$
\pi_n^{\ast}(\ud \theta|\gamma) = \frac{\exp\left\{- \frac{n}{2\sigma^2} \int_{[a,b]\times\reals} [\theta(u) - v]^2\gamma(\ud u \ud v) \right\} \pi(\ud\theta)}{\int_{\ps} \exp\left\{- \frac{n}{2\sigma^2} 
\int_{[a,b]\times\reals} [\tau(u) - v]^2\gamma(\ud u \ud v)\right\} \pi(\ud\tau)}\ .
$$
Lastly, as for the Kullback-Leibler divergence, a straightforward computation yields
$$
\kksf(\theta|\theta_0) = \frac{1}{2\sigma^2} \int_a^b [\theta(u) - \theta_0(u)]^2 h(u) \ud u\ .
$$
This statistical model is particularly versatile with respect to our theory, because it can be studied as either an infinite-dimensional exponential family or by means of  
Theorem \ref{main_thm1} and Corollary \ref{coro1}. For example, to see that we can use the theory of infinite-dimensional exponential families, it suffices to consider the identities
\begin{align*}
(\theta(u) - v)^2 &= \int_a^b\int_{\R} (\theta(x)-y)^2 \delta_{(u,v)}(\ud x\ud y) \\
&= -\int_a^b\int_{\R} [\Delta((\theta(x)-y)^2)] \mathcal G(x,y;u,v) \ud x\ud y
\end{align*}
where $\mathcal G(x,y;u,v)$ stands for the Green function of the set $[a,b]\times \R$. If $\theta$ varies in a sufficiently regular space, that is if $s$ is sufficiently large, then $[\Delta((\theta(x)-y)^2)]$ is still a function, which can be
set equal to $g(\theta)$. On the other hand, $\mathcal G(x,y;u,v)$ represents the function $\beta$ in the theory of exponential families. 

As for the assumptions of Theorem \ref{main_thm1}, we can prove their validity if, for instance, $h$ belongs to $C^0[a,b] \cap C^2(a,b)$ and it is bounded away from zero. The assumption $\int_{\ss} |x|^q \mu_0(\ud x) < +\infty$ is valid for any $q>0$ and \eqref{bound_moment_posterior} holds if we assume, for instance, a Gaussian prior $\pi$. As for Corollary \ref{coro1}, we can check the validity of \eqref{bound_Kullback} as a consequence of the Gagliardo-Nirenberg interpolation inequality (\citet[Section 12.3]{Mazya(11)}). Being $\kksf(\theta|\theta_0)$ equivalent to the squared $L^2$-norm,
$$
\|\theta - \theta_0\|_{H^s(a,b)} \leq \|\theta - \theta_0\|_{L^2(a,b)}^{1-\alpha}  \|\theta - \theta_0\|_{H^{s'}(a,b)}^{\alpha}
$$
for any $s' > s$, where $\alpha := s/s'$. Therefore, choosing a prior distribution that is supported on $H^{s'}(a,b)$, such as for instance a Gaussian type prior, and recalling that $H^{s'}(a,b)$ is dense in $H^s(a,b)$, it is enough to consider the neighborhood $\|\theta - \theta_0\|_{H^{s'}(a,b)} \leq 1$ of $\theta_0$ and check that the interpolation inequality immediately yields \eqref{bound_Kullback}. Whence, $\beta = 2/(1-\alpha)$. These considerations show that Proposition \ref{prp:Laplace_infinite} is applicable, provided that the prior is Gaussian with a covariance matrix that satisfies \eqref{eigenvalues_logistic}. In any case, both the methods end up by highlighting the main terms that figures on the right-hand sides of \eqref{Main_PCR1} and \eqref{PCR_main1}.

Now, for the sake of brevity, we confine ourselves on the application of Theorem \ref{new_therem_PCR}. Apropos of the first term on the right-hand side of \eqref{Main_PCR1}, we notice that $\mathrm{I}(\theta_0)$ is independent of $\theta_0$, and is equivalent to the identity operator. In view of a 
straightforward coercivity, we can apply the results in Section 3.3 of \citet{AS(99)} to obtain that the first term on the right-hand side of \eqref{Main_PCR1} is asymptotic to $\frac{1}{\sqrt{n}}$, as $n \to
+\infty$. Then, the second and the third terms are exponentially small, and hence asymptotically negligible.  
To complete the treatment, we are left to discuss the asymptotic behavior of the constant $L_0^{(n)}$. Apropos of the Poincar\'e constant $[\mathfrak C_2(\pi_n^{\ast}(\cdot|\gamma))]^2$, here it is trivial to notice that the mapping $\theta \mapsto \int_{[a,b]\times\reals} [\theta(u) - v]^2\gamma(\ud u \ud v)$ is twice Frech\'et-differentiable with respect to $\theta$. Therefore, the Bakry-Emery criterion applies and, if $\pi$ is Gaussian, results in \citet[Chapters 10--11]{DaPrato(06)} show that $[\mathfrak C_2(\pi_n^{\ast}(\cdot|\gamma))]^2 = O(1/n)$, as $n \to \infty$. As for the term $\mathfrak D_{\theta} \frac{\nabla_x f(x|\theta)}{f(x|\theta)}$, we first notice that
$$
\frac{\nabla_x f(x|\theta)}{f(x|\theta)} = -\frac{1}{\sigma^2} \left( \theta'(u)(\theta(u) - v), -(\theta(u) - v) \right) \ . 	
$$
This is the sum of the terms $\sigma^{-2}\left( v\theta'(u), \theta(u) - v\right)$ and $\frac{1}{\sigma^2} \left( -\theta'(u)\theta(u), 0 \right)$, where the former vector is a linear functional of $\theta$. Thus, the Fr\'echet derivative of the first term with respect to $\theta$ is given by the vector $\sigma^{-2}(vS_u, T_u)$, where $T_u$ ($S_u$, respectively) stands for the Riesz representative of the functional $\delta_u$ ($-\delta'_u$, respectively). It is useful to observe that such a derivative, being independent of $\theta$, does not contribute asymptotically in the expression of the double integral, as we have already discussed in the previous section. Finally, the Fr\'echet derivative of the second term is $-\sigma^{-2}\left( T_u\theta'(u) + S_u\theta(u), 0 \right)$. 
At this stage, we can see that the study of the double integral can be reduced, through the use of Sobolev inequalities, to the study of the corresponding posterior moments. To conclude, we state a proposition that summarizes the above considerations.

\begin{prp}
In connection with the model \eqref{lin_regression}, let $\ss = [a,b] \times \R$ and 
$\ps = \mathrm{H}^s(0,1)$ for some $s \geq 3$. Suppose that $h$ is a smooth density on $[a,b]$, satisfying $1/c \leq h(x) \leq c$ for any $x \in [a,b]$ and some $c>0$. 
Let $\theta_0 \in \ps$ be fixed. Assume that $\pi = \mathcal N(m,Q)$ with any $m \in \ps$ and $Q$ a non-degenerate trace-class
operator. Then, it holds that
$$
\epsilon_n = O\big(n^{-1/2}\big)
$$
as $n \rightarrow +\infty$, which represents the optimal rate. 
\end{prp}


\section{Discussion}\label{sec:discuss}
We conclude our work by discussing some directions for future research. The flexibility of the Wasserstein distance is promising when considering non-regular Bayesian statistical models, even in a finite-dimensional setting. One may consider the problem of dealing with dominated statistical models that have moving supports, i.e. supports that depend on $\theta$. The prototypical example is the family of Pareto distributions, which is characterized by a density function
\begin{displaymath}
f(x\,|\,\alpha, x_0) = \frac{\alpha x_0^{\alpha}}{x^{1+\alpha}} \mathds{1}\{x \geq x_0\},
\end{displaymath}
where $\theta = (\alpha, x_0) \in (0,+\infty)^2$. Under this model, by rewriting the posterior distribution to obtain the representation  \eqref{representation_kernel}, we observe that the empirical distribution can be replaced by the minimum of the observations, which is the  maximum likelihood estimator. In doing this, we expect to parallel the proof of Theorem \ref{main_thm1}, with the minimum playing  the role of the sufficient statistic, instead of the empirical measure. In particular, we expect that the term $\varepsilon_{n,p}(\ss, \mu_0)$ should be replaced by other rates typically involved in limit theorems of order statistics. The theoretical framework for such an extension of our results is developed in the work of \citet{DM(20b)}, where it is shown how the continuity equation yields a specific boundary-value problem of Neumann type.

As for the infinite-dimensional setting covered by Theorem \ref{new_therem_PCR} and Theorem \ref{main_thm1}, an interesting development of our approach to PCRs is represented by the possibility of finding, for general statistical models, explicit sufficient statistics belonging to Banach spaces of functions. To be more precise, we hint at a constructive version of the well-known Fisher-Neyman factorization lemma. This result would pave the way for a suitable rewriting of the statistical model, that allows for the use of our approach. By way of example, one may consider the identity $\log f(x | \theta) = \int_{\X} \log f(y | \theta) \delta_x(\ud y)$, and exploit an integration-by-part formula to obtain an identity like \eqref{mu_n_exp}, with respect to a suitable measure $\lambda$ on $(\X, \mathscr X)$. Such a procedure is at the basis for the development of our approach to PCRs in the context of popular nonparametric models, not considered in this paper, such as the Dirichlet process mixture model (\citet[Chapter 5]{GV(00)}), the random histograms (\citet[Example 5.11]{GV(00)}) and P\'olya trees (\citet[Section 3.7]{GV(00)}).

Another promising line of research consists in extending Theorem \ref{main_thm1} to metric measure spaces. The theoretical ground for this development may be found in the seminal works of \citet{Gig(09)}, \citet{Gig(12)}, \citet{AGS(14)}, \citet{OV(00)} and \citet{vRS(09)}. In such a context, it is of interest the treatment of the relative entropy-functional in the Wasserstein space. It is well-known that the Hessian of the relative entropy-functional, i.e. the Kullback-Leibler divergence, generalizes by using techniques from infinite-dimensional Riemannian geometry (\citet{OV(00)}). From the statistical side, the possibility of choosing a parameter space that coincides with a space of measures allows to re-consider, from a different point of view, popular Bayesian statistical models such as Dirichlet process mixture models, which are defined as
\begin{displaymath}
 f(x\,|\,\pfrak) = \int \tau(x\,|\,y) \pfrak(\ud y)
\end{displaymath}
where $\tau$ is a kernel parameterized by $y$, and $\pfrak$ is a random probability measure with a Dirichlet process prior (\citet{Fer(73)}). The goal should be that of considering PCRs relative to Wasserstein neighborhoods of a given true distribution, say $\pfrak_0$. This approach is again different from the nonparametric framework considered in \citet{Bert(21)}, and seems still unexplored. 


\appendix

\section{Proofs}

\subsection{Proof of Lemma \ref{lm:Wpcr}} \label{proof:lm1}
By a standard measure theoretic argument, any two solutions $\pi_n(\cdot|\cdot)$ and $\pi'_n(\cdot|\cdot)$ of \eqref{disintegration} satisfy $\pi_n(\cdot|x^{(n)}) = \pi'_n(\cdot|x^{(n)})$ as elements of
$\pms$, for all $x^{(n)} \in \ss^n\setminus N_n$, where $N_n$ is a $\alpha_n$-null set. The assumption $\mu_0^{\otimes_n} \ll \alpha_n$ entails that $\xi^{(n)} := (\xi_1, \dots, \xi_n)$ takes values in $N_n$ with $\P$-probability zero, yielding 
the desired well-definiteness. 

Then, if $\pi \in \mathcal P_p(\ps)$, any solution $\pi_n(\cdot|\cdot)$ of \eqref{disintegration} satisfies $\pi_n(\mathcal P_p(\ps)|x^{(n)}) = 1$ for $\alpha_n$-almost every $x^{(n)}\in\ss^n$. Since $\mu_0^{\otimes_n} \ll \alpha_n$,
it follows that $\pi_n(\mathcal P_p(\ps)|\xi^{(n)}) = 1$ with $\P$-probability 1. Whence, 
$$
\Wp(\pi_n(\cdot|\xi^{(n)}); \delta_{\theta_0}) = \left( \int_{\ps} [\ud_{\ps}(\theta,\theta_0)]^p \pi_n(\ud\theta|\xi^{(n)}) \right)^{1/p}
$$
is a random variable, which proves to be finite $\P$-a.s.. Combining Markov's and Lyapunov's inequalities, it follows that
$$
\pi_n\left(\left\{\theta \in \ps\ :\ \ud_{\ps}(\theta,\theta_0) \geq M_n\epsilon_n\right\} \big| \xi^{(n)} \right) \leq \frac{\left( \int_{\ps} [\ud_{\ps}(\theta,\theta_0)]^p \pi_n(\ud\theta|\xi^{(n)}) \right)^{1/p}}{M_n \epsilon_n}
$$
holds $\P$-a.s.. Now, taking expectation of both sides and taking account of \eqref{Wpcr} yields
$$
\E\left[ \pi_n\left(\left\{\theta \in \ps\ :\ \ud_{\ps}(\theta,\theta_0) \geq M_n\epsilon_n\right\} \big| \xi^{(n)} \right) \right] \leq \frac{1}{M_n} \rightarrow 0\ . 
$$
Thus, the convergence indicated in \eqref{post_consistency} holds in $\mathrm L^1\probabilityspace$ and, hence, in $\P$-probability. The proof is complete.


\subsection{Proof of identity \eqref{representation_kernel_kullback}} \label{proof:lm2}

In view of \eqref{kernel_exp}, it is enough to prove that
\begin{equation} \label{kullback_Appendix}
_{\B^{\ast}}\!\langle g(\theta), b\rangle_{\B} - M(\theta) = -\kksf(\theta\ |\ \theta_b) + \mathsf{H}(b)
\end{equation}
holds for all $\theta \in \ps$ and $b$ in the range of $\mathcal S \circ g$, with some suitable function $\mathsf{H} : \text{Range}(\mathcal S \circ g) \to \R$. Then, \eqref{mu_n_exp} yields
\begin{align*}
-\kksf(\theta\ |\ \theta_b) &= -\int_{\ss} \log\left(\frac{\varphi(x\ |\ g(\theta_b))}{\varphi(x\ |\ g(\theta))}\right) \mu(\ud x\ |\ \theta_b) \\
&= M(\theta_b) - M(\theta) + \int_{\ss}\ _{\B^{\ast}}\!\langle g(\theta), \beta(x) \rangle_{\B}\ \mu(\ud x\ |\ \theta_b) \\
&-  \int_{\ss}\  _{\B^{\ast}}\!\langle g(\theta_b), \beta(x) \rangle_{\B}\ \mu(\ud x\ |\ \theta_b)\ .
\end{align*}
Combining the above identity with \eqref{definition_S} and observing that $g(\theta_b) = \mathcal S^{-1}(b)$, it follows that
\begin{align*}
-\kksf(\theta\ |\ \theta_b) &= M(\theta_b) - M(\theta) + _{\B^{\ast}}\!\langle g(\theta), \mathcal S(g(\theta_b)) \rangle_{\B} -\ _{\B^{\ast}}\!\langle g(\theta_b), \mathcal S(g(\theta_b)) \rangle_{\B} \\
&=\  _{\B^{\ast}}\!\langle g(\theta), b\rangle_{\B} - M(\theta) + [M(\theta_b) -\  _{\B^{\ast}}\!\langle g(\theta_b), b\rangle_{\B} ]
\end{align*}
is valid for all $\theta \in \ps$ and $b$ in the range of $\mathcal S \circ g$. Then, the validity of \eqref{kullback_Appendix} follows by putting $\mathsf{H}(b) :=\  _{\B^{\ast}}\!\langle g(\theta_b), b\rangle_{\B} - M(\theta_b)$, 
completing the proof.


\subsection{Proof of Theorem \ref{new_therem_PCR}} \label{proof:new_thm}

Under the assumptions of the Theorem, Lemma \ref{lm:Wpcr} is valid, and a PCR at $\theta_0$ is given by \eqref{Wpcr}. Moreover, \eqref{representation_kernel} is valid with $\mathfrak{S}_n(\xi_1, \dots, \xi_n) = \hat{S}_n$,
where $\hat{S}_n$ is given by \eqref{S_n_exp}, $\Sd = \B$ endowed with the distance ensuing from the norm $\|\cdot\|_{\B}$, and the kernel $\pi_n^{\ast}(\cdot|\cdot)$ is. given by \eqref{kernel_exp}. The triangle inequality 
for $\Wp^{(\mathcal P(\ps))}$ gives
$$
\Wp^{(\mathcal P(\ps))}(\pi_n^{\ast}(\cdot| \hat{S}_n); \delta_{\theta_0}) \leq \Wp^{(\mathcal P(\ps))}(\pi_n^{\ast}(\cdot| S_0); \delta_{\theta_0}) + \Wp^{(\mathcal P(\ps))}(\pi_n^{\ast}(\cdot| \hat{S}_n); \pi_n^{\ast}(\cdot| S_0))
$$
with the same $S_0$ as in \eqref{S_0_exp}. See \cite[Chapter 7]{AGS(08)} for information about the aforesaid triangle inequality. Then, take the expectation of both sides above to obtain
\begin{align} \label{triangular1}
\epsilon_n &
\leq \Wp^{(\mathcal P(\ps))}(\pi_n^{\ast}(\cdot| S_0); \delta_{\theta_0}) + 
\E\left[\Wp^{(\mathcal P(\ps))}(\pi_n^{\ast}(\cdot| \hat{S}_n); \pi_n^{\ast}(\cdot| S_0))\right] \nonumber \\
&= \Wp^{(\mathcal P(\ps))}(\pi_n^{\ast}(\cdot| S_0); \delta_{\theta_0}) \nonumber \\
&+ \E\left[ \Wp^{(\mathcal P(\ps))}(\pi_n^{\ast}(\cdot| S_0); \pi_n^{\ast}(\cdot| \hat{S}_n)) \mathds{1}\{\hat{S}_n \in \mathcal{U}_{\delta_n}(S_0)\}\right] \nonumber\\
&+ \E\left[ \Wp^{(\mathcal P(\ps))}(\pi_n^{\ast}(\cdot| S_0); \pi_n^{\ast}(\cdot| \hat{S}_n)) \mathds{1}\{\hat{S}_n \not\in \mathcal{U}_{\delta_n}(S_0)\}\right] \ .
\end{align}
At this stage, the first summand on the last member of \eqref{triangular1} is exactly equal to the first summand on the right-hand side of \eqref{Main_PCR1}, thanks to identity \eqref{Wasserstein_nonrandom}. 
For the second summand on the last member of \eqref{triangular1}, invoke \eqref{Lipschitz_kernel} to conclude that such term is majorized by the last summand on the right-hand side of \eqref{Main_PCR1}. It remains to handle 
the third summand on the last member of \eqref{triangular1}. Exploit the fact that, for any two elements $\mu, \nu \in \mathcal P_p(\ps)$ there holds
$$
\Wp^{(\mathcal P(\ps))}(\mu; \nu) \lesssim \left[\int_{\ps} \|\theta\|_{\ps}^p \mu(\ud\theta)\right]^{\frac 1p} + \left[\int_{\ps} \|\theta\|_{\ps}^p \nu(\ud\theta)\right]^{\frac 1p}
$$
to obtain that 
\begin{align} \label{triangular2}
& \E\left[ \Wp^{(\mathcal P(\ps))}(\pi_n^{\ast}(\cdot| S_0); \pi_n^{\ast}(\cdot| \hat{S}_n)) \mathds{1}\{\hat{S}_n \not\in \mathcal{U}_{\delta_n}(S_0)\}\right] \nonumber \\
&\lesssim \E\left[\mathds{1}\{\hat{S}_n \not\in \mathcal{U}_{\delta_n}(S_0)\} \left( \int_{\ps} \|\theta\|_{\ps}^p \pi_n^{\ast}(\ud\theta| \hat{S}_n) \right)^{\frac 1p}\right] \nonumber \\
&+ \left( \int_{\ps} \|\theta\|_{\ps}^p \pi_n^{\ast}(\ud\theta| S_0) \right)^{\frac 1p} \P\left[ \hat{S}_n \not\in \mathcal{U}_{\delta_n}(S_0) \right]\ . 
\end{align}
Now, the first summand on the right-hand side of \eqref{triangular2} can be bounded by means of a combination of H\"older's and Lyapunov's inequalities, yielding
\begin{align*}
&\E\left[\mathds{1}\{\hat{S}_n \not\in \mathcal{U}_{\delta_n}(S_0)\} \left( \int_{\ps} \|\theta\|_{\ps}^p \pi_n^{\ast}(\ud\theta| \hat{S}_n) \right)^{\frac 1p}\right] \\
&\leq \left( \E\left[ \int_{\ps} \|\theta\|_{\ps}^{ap} \pi_n^{\ast}(\ud\theta\,|\, \hat{S}_n) \right] \right)^{\frac{1}{ap}} \left(\P\left[\hat{S}_n \not\in \mathcal{U}_{\delta_n}(S_0)\right]\right)^{1 - \frac{1}{ap}}\ . 
\end{align*}
For the second summand on the right-hand side of \eqref{triangular2} just exploit the triangular inequality to obtain  
$$
\left( \int_{\ps} \|\theta\|_{\ps}^p \pi_n^{\ast}(\ud\theta| S_0) \right)^{\frac 1p} \lesssim \|\theta_0\|_{\ps} + \left( \int_{\ps} \|\theta - \theta_0\|_{\ps}^p \pi_n^{\ast}(\ud\theta| S_0) \right)^{\frac 1p}\ . 
$$ 
Re-organizing the terms just obtained yields the right-hand side of \eqref{Main_PCR1}, concluding the proof. 


\subsection{Proof of Proposition \ref{prp:Laplace_infinite}} \label{proof:Laplace_infinite}

Start by fixing $\varepsilon$ in the interval $\big(2(r-1), q\big)$, which is possible since $0 < 2(r-1) < q$.
Then, let $\{\eta_n\}_{n \geq 1}$ be a sequence of positive numbers such that $\eta_n = O(n^{-1/(2+\varepsilon)})$ as $n \to +\infty$. 
Let $B_{\eta_n}(\theta_0)$ denote the open ball in $\ps$ with radius $\eta_n$, centered at $\theta_0$. Without loss of generality, assume that $\theta_0 \in \Vbb$. 
Otherwise, by density of $\Vbb$, pick a sequence $\{\theta_{0,n}\}_{n \geq 1} \subset \Vbb$ such that 
$\|\theta_{0,n} - \theta_0\|_{\ps} \to 0$ sufficiently fast, and replace $\theta_0$ by $\theta_{0,n}$.

The proof is divided into four steps, according to typical operations in the theory of Laplace approximation. First, let us prove that
\begin{equation} \label{Laplace_first_step}
\int_{\ps} \|\theta - \theta_0\|_{\ps}^2 \pi_n^{\ast}(\ud\theta| S_0) \sim  
\frac{ \int_{B_{\eta_n}(\theta_0)} \|\theta - \theta_0\|_{\ps}^2 \exp\{-n \kksf(\theta | \theta_0)\}\pi(\ud\theta)}{\int_{B_{\eta_n}(\theta_0)} \exp\{-n \kksf(\theta | \theta_0)\}\pi(\ud\theta)}
\end{equation}
as $n \to +\infty$ where, for any pair of sequences $\{a_n\}_{n\geq 1}$ and $\{b_n\}_{n\geq 1}$ of positive numbers, the notation $a_n \sim b_n$ means that $\lim_{n \to +\infty} a_n/b_n = 1$. 
To this aim, it is enough to show that the integrals on the exterior of $B_{\eta_n}(\theta_0)$ are exponentially small, and hence irrelevant in the global asymptotic expansion. 
Exploiting \eqref{interpolation}--\eqref{lower_kullback}, one gets
\begin{align*}
&\int_{B_{\eta_n}(\theta_0)^c} \exp\{-n \kksf(\theta | \theta_0)\}\pi(\ud\theta)\\
&\quad\leq \int_{B_{\eta_n}(\theta_0)^c} \exp\{-n \phi(\|\theta-\theta_0\|_{\Kbb})  \}\pi(\ud\theta) \\
&\quad\leq \int_{B_{\eta_n}(\theta_0)^c} \exp\left\{-n \phi\left( \left[\frac{\|\theta-\theta_0\|_{\ps}}{\|\theta-\theta_0\|_{\Vbb}^{1/s}}\right]^r \right)\right\}\pi(\ud\theta) \\
&\quad\leq \int_{\Vbb} \exp\left\{-n \phi\left( \left[\frac{\eta_n}{\|\theta-\theta_0\|_{\Vbb}^{1/s}}\right]^r \right)\right\}\pi(\ud\theta) \ .
\end{align*}
Now, let $B_{\rho_n}^{(\Vbb)}(\theta_0)$ denote the open ball in $\Vbb$, with radius $\rho_n$ and centered at $\theta_0$. Thus, the last integral can be bounded from above by
\begin{equation} \label{Bernstein_Exterior}
\pi\left(B_{\rho_n}^{(\Vbb)}(\theta_0)^c \right) + \exp\left\{-n \phi\left(\left[\frac{\eta_n}{\rho_n^{1/s}}\right]^r\right)\right\}
\end{equation}
which can be made an exponentially small quantity after choosing properly the sequence $\{\rho_n\}_{n \geq 1}$. Actually, it is enough to fix that $\rho_n = O(n^h)$ as $n \to +\infty$, for some 
$h$ satisfying
\begin{equation} \label{ineq_rteps}
0 < h < \left(\frac{1}{2r} - \frac{1}{2+\varepsilon} \right) \frac{r}{r-1} \ .
\end{equation}
Of course, this is possible in view of the bound $2+\varepsilon > 2r$. Now, $h > 0$ entails that $\eta_n \rho_n^{-1/s} \to 0$ as $n \to +\infty$. Then, exploiting that 
$\phi(x) = O(x^2)$ as $x \to 0^+$, one gets
$$
n \phi\left(\left[\eta_n \rho_n^{-1/s} \right]^r\right) \sim n^{1 - 2r\left(\frac{1}{2+\varepsilon} + \frac{h}{s}\right)}\ .
$$
Lastly, combination of the identity $s = \frac{r}{r-1}$ with the inequality \eqref{ineq_rteps} entails that
$$
c := 1 - 2r\left(\frac{1}{2+\varepsilon} + \frac{h}{s}\right) > 0\ .
$$
This argument shows that the second summand in \eqref{Bernstein_Exterior} goes to zero like $e^{-n^c}$, making it a negligible quantity. Finally, the first summand in \eqref{Bernstein_Exterior}
is also bounded by a term that goes to zero like $e^{-n^h}$, thanks to a straightforward combination of Markov's inequality with 
the assumption that $\int_{\Vbb} e^{t \|\theta\|_{\Vbb}} \pi(\ud\theta) < +\infty$ for some $t > 0$.

As for the term
$$
\int_{B_{\eta_n}(\theta_0)^c} \|\theta - \theta_0\|_{\ps}^2 \exp\{-n \kksf(\theta | \theta_0)\}\pi(\ud\theta)
$$
the argument to prove that it is also exponentially small is similar. Indeed, it is enough to get rid of the term $\|\theta - \theta_0\|_{\ps}^2$ by a straightforward application of H\"older inequality.
This proves \eqref{Laplace_first_step}.

After reducing both the integrals on $B_{\eta_n}(\theta_0)$, exploit the regularity of the map $\theta \mapsto \kksf(\theta|\theta_0)$ by showing that it can be replaced by its second order
Taylor polynomial, which reads
$$
\frac 12 \langle \theta-\theta_0, \mathrm{I}(\theta_0)[\theta-\theta_0] \rangle
$$
because $\kksf(\theta_0|\theta_0) = 0$ and $\mathrm{D}_{\theta}\kksf(\theta|\theta_0) _{|{\theta=\theta_0}} = 0$. By the assumptions of the proposition,
$$
\left| \kksf(\theta_0|\theta_0) - \frac 12 \langle \theta-\theta_0, \mathrm{I}(\theta_0)[\theta-\theta_0] \rangle \right| \sim n \left(\frac{1}{\sqrt n}\right)^{2+q} \to 0
$$
which entails that the two integrals
\begin{align*}
\int_{B_{\eta_n}(\theta_0)} \|\theta - \theta_0\|_{\ps}^2 \left| \exp\{-n \kksf(\theta | \theta_0)\} - \exp\left\{ - \frac{n}{2} \langle \theta-\theta_0, \mathrm{I}(\theta_0)[\theta-\theta_0] \rangle \right\} \right| \pi(\ud\theta) \\
\int_{B_{\eta_n}(\theta_0)} \left| \exp\{-n \kksf(\theta | \theta_0)\} - \exp\left\{ - \frac{n}{2} \langle \theta-\theta_0, \mathrm{I}(\theta_0)[\theta-\theta_0] \rangle \right\} 
\right| \pi(\ud\theta)
\end{align*}
go to zero faster than their respective counterparts
\begin{align*}
\int_{B_{\eta_n}(\theta_0)} \|\theta - \theta_0\|_{\ps}^2 \exp\{-n \kksf(\theta | \theta_0)\}\pi(\ud\theta) \\
\int_{B_{\eta_n}(\theta_0)}  \exp\{-n \kksf(\theta | \theta_0)\}\pi(\ud\theta) \ .
\end{align*}
Whence,
\begin{align} \label{Laplace_second_step}
&\frac{ \int_{B_{\eta_n}(\theta_0)} \|\theta - \theta_0\|_{\ps}^2 \exp\{-n \kksf(\theta | \theta_0)\}\pi(\ud\theta)}{\int_{B_{\eta_n}(\theta_0)} \exp\{-n \kksf(\theta | \theta_0)\}\pi(\ud\theta)}\\
&\notag\quad \sim\frac{ \int_{B_{\eta_n}(\theta_0)} \|\theta - \theta_0\|_{\ps}^2 \exp\{- \frac n2 \langle \theta-\theta_0, \mathrm{I}(\theta_0)[\theta-\theta_0] \rangle\}\pi(\ud\theta)}{\int_{B_{\eta_n}(\theta_0)} 
\exp\{- \frac n2 \langle \theta-\theta_0, \mathrm{I}(\theta_0)[\theta-\theta_0] \rangle \} \pi(\ud\theta)}  \ . 
\end{align}
This concludes the second step. The third step is similar to the first one, the goal being to show that
\begin{align} \label{Laplace_third_step}
&\frac{ \int_{B_{\eta_n}(\theta_0)} \|\theta - \theta_0\|_{\ps}^2 \exp\{- \frac n2 \langle \theta-\theta_0, \mathrm{I}(\theta_0)[\theta-\theta_0] \rangle\}\pi(\ud\theta)}{\int_{B_{\eta_n}(\theta_0)} 
\exp\{- \frac n2 \langle \theta-\theta_0, \mathrm{I}(\theta_0)[\theta-\theta_0] \rangle \} \pi(\ud\theta)}\\
&\notag\quad  \sim 
\frac{ \int_{\ps} \|\theta - \theta_0\|_{\ps}^2 \exp\{- \frac n2 \langle \theta-\theta_0, \mathrm{I}(\theta_0)[\theta-\theta_0] \rangle\}\pi(\ud\theta)}{\int_{\ps} 
\exp\{- \frac n2 \langle \theta-\theta_0, \mathrm{I}(\theta_0)[\theta-\theta_0] \rangle \} \pi(\ud\theta)}\ . 
\end{align}
The argument is similar, just utilize \eqref{lower_fisher} instead of \eqref{lower_kullback}. This concludes the third step. Lastly, observe that the right-hand side of \eqref{Laplace_third_step}
coincides with the ratio of two Gaussian integrals, which are factorized in view of the assumption \eqref{eigenvalues_logistic}. An explicit computation now gives
\begin{align} \label{Laplace_fourth_step}
&\frac{ \int_{\ps} \|\theta - \theta_0\|_{\ps}^2 \exp\{- \frac n2 \langle \theta-\theta_0, \mathrm{I}(\theta_0)[\theta-\theta_0] \rangle\}\pi(\ud\theta)}{\int_{\ps} 
\exp\{- \frac n2 \langle \theta-\theta_0, \mathrm{I}(\theta_0)[\theta-\theta_0] \rangle \} \pi(\ud\theta)}\\
&\notag\quad =\sum_{k=1}^{\infty} \frac{\lambda_k}{n \lambda_k \gamma_k + 1} + \sum_{k=1}^{\infty} \frac{\omega_k^2}{(n \lambda_k \gamma_k + 1)^2}\ .
\end{align}
This concludes the fourth step and the proof. 


\subsection{Proof of Proposition \ref{prop:Bayes_Club}} \label{proof:Bayes_Club}

The main issue is to prove the validity of \eqref{Lipschitz_kernel}. Thus, fix $S_0 \in \B$ and $S' \in \mathcal{U}_{\delta_n}(S_0)$. 
For $t$ varying in $[0,1]$, let $S_t = S_0 + t(S' - S_0)$ denote the line-segment joining $S_0$ with $S'$.
Use the kernel $\pi_n^{\ast}(\cdot|\cdot)$ defined in \eqref{kernel_exp} to lift the line-segment $[S_t]_{t \in [0,1]}$ to $\mathcal P_2(\ps)$, by means of the new curve
$$
\mu_t^{\ast}(\cdot) := \pi_n^{\ast}(\cdot| S_t)
$$
which joins $\pi_n^{\ast}(\cdot| S_0)$ with $\pi_n^{\ast}(\cdot| S')$. Here, we apply the Benamou-Brenier representation introduced in Section \ref{sect:main_problem} with $\mathbb M = \overline{\Theta}$, to get
$$
\Wdue^2(\pi_n^{\ast}(\cdot| S_0); \pi_n^{\ast}(\cdot| S')) \leq  \int_0^1 \!\!\!\!\int_\ps \!\! \| \D_{\theta} u^{\ast}(\theta,t)\|^2\,\mu_t^{\ast}(\ud\theta)\,\ud t = \int_0^1 \!\!\! \|u^{\ast}(\cdot, t)\|^2_{1,\mu_t^{\ast}} \ud t
$$
where $u^{\ast}(\cdot,t)$ is the (unique) solution in $\mathrm{H}^1_m(\ps; \mu^{\ast}_t)$ of \eqref{cont_eq} with $\gamma_t = \mu^{\ast}_t$. Here, $\mathrm{H}^1_m(\ps; \mu^{\ast}_t)$ is defined 
as the completion of the space
$$
\left\{\psi \in C^1_b(\overline{\ps})\ \Big|\ \int_{\ps} \psi(\theta) \mu^{\ast}_t(\ud\theta) = 0 \right\}
$$
with respect to the norm 
$$
\|\psi\|_{1,t} := \left( \int_{\ps} \|\mathrm D_{\theta} \psi(\theta)\|^2 \mu^{\ast}_t(\ud\theta) \right)^{1/2}
$$ 
associated with the scalar product
$$
\langle \varphi, \psi \rangle_{1,t} := \int_{\ps} \langle\mathrm D_{\theta}\ \varphi(\theta), \mathrm D_{\theta} \ \psi(\theta)\rangle \mu^{\ast}_t(\ud\theta)\ .
$$
Then, rewrite \eqref{cont_eq} as
$$
\T_t[\psi] = \langle \psi, u^{\ast}(\cdot, t) \rangle_{1, \mu^{\ast}_t}
$$
with
$$
\T_t[\psi] := \frac{\ud}{\ud t}\int_\ps \psi(\theta)\,\mu^{\ast}_t(\ud\theta) = \frac{\ud}{\ud s}\int_\ps \psi(\theta)\,\mu_s^{\ast}(\ud\theta) \ _{\big{|}_{s=t}}\ . 
$$
By Riesz representation, we get 
$$
\|u^{\ast}(\cdot, t)\|_{1,\mu_t^{\ast}} = \sup_{\|\psi\|_{1,\mu_t^{\ast}} \leq 1} | \T_t[\psi] |\ .
$$
Now, take the derivative inside the integral in the expression of $\T_t$, consider the expression of $\mu_t^{\ast}$ and apply the Leibnitz rule, as follows.
\begin{align*}
    |\T_t[\psi]| &= \left| \int_\ps \psi(\theta)\, \left[ \frac{\partial}{\partial t} \frac{\exp\{n[\langle g(\theta), S_t\rangle] - M(\theta)}{\int_{\ps} \exp\{n[\langle g(\tau), S_t\rangle 
    - M(\tau)]\} \pi(\ud\tau)} \right]\pi(\ud\theta)\right| \\
    & \stackrel{\text{(Leibnitz)}}{=} n \ \left|\textsf{Cov}_{\mu_t^{\ast}} \Big( \psi(\cdot), \langle g(\cdot), S' - S_0\rangle \Big)\right| \\
    & \stackrel{\text{(Cauchy-Scwartz)}}{\leq} n \sqrt{\textsf{Var}_{\mu_t^{\ast}}( \psi(\cdot))} \sqrt{\textsf{Var}_{\mu_t^{\ast}} \Big( \langle g(\cdot), S' - S_0\rangle \Big)} \\
    & \stackrel{\text{(Poincar\'e-Wirtinger)}}{\leq} n \{\Cpw[\mu_t^{\ast}]\}^2 \underbrace{
    \| \psi\|_{1, \mu_t^{\ast}}}_{\leq 1} \|\langle g(\cdot), S' - S_0\rangle\|_{1, \mu_t^{\ast}} \\
    &\stackrel{\text{(duality)}}{\leq} 
    n \{\Cpw[\mu_t^{\ast}]\}^2 \|S' - S_0\|_{\B}
    \left( \int_{\ps} \|\mathfrak{D}_{\theta}[g]\|_{\ast}^2\ \mu_t^{\ast}(\ud\theta) \right)^{1/2} \\
    &\leq \|S' - S_0\|_{\B}\ \underbrace{n \sup_{S \in \mathcal{U}_{\delta_n}(S_0)} \{\Cpw[\pi_n^{\ast}(\cdot| S)] \}^2 \left( \int_{\ps} \|\mathfrak{D}_{\theta}[g]\|_{\ast}^2\ \pi_n^{\ast}(\ud\theta| S) 
    \right)^{1/2}}_{= L_0^{(n)}}
    \end{align*}
where in the inequality with the super-script ``duality'' we have used the fact that, for any $b \in \B$, it holds
$$
\|\D_{\theta}\left[\langle g(\cdot), b \rangle\right]\|_{\ps} = |\langle \mathfrak{D}_{\theta}[g], b\rangle|\ \leq  \|\mathfrak{D}_{\theta}[g]\|_{\ast} \|b\|_{\B}\  .
$$
This proves inequality \eqref{L0_bound1}. Finally, \eqref{L0_bound2} follows trivially from \eqref{L0_bound1}, in view of the boundedness condition \eqref{Laplace_rafforzato}. 


\subsection{Proof of Proposition \ref{prp:Poincare_infinite}} \label{proof:Poincare_infinite}

The first step of the proof is to provide a result analogous to \citet[Theorem 1.4]{BBCG(08)}. To this aim, we need the concept of Lyapunov function, as done in that paper. Therefore, let $V : \ps \to \R$ a $\mathrm{C}^2$ function bounded from below. Define the probability measure $\mu_{V,\pi}$ in Gibbsean form as
\begin{equation} \label{Gibbs_Vpi}
\mu_{V,\pi}(\ud\theta) = \frac{e^{-V(\theta)}\pi(\ud\theta)}{\int_{\ps} e^{-V(\tau)}\pi(\ud\tau)}\ . 
\end{equation}
Then, define the differential operator $\mathfrak{L}_{V,\pi} := \mathfrak{L}_{\pi} - \Dcal_{\theta}[V] \cdot \Dcal_{\theta}$, where $\Dcal$ and $\mathfrak{L}_{\pi}$ denote the Malliavin derivative and the 
Malliavin-Laplace operator associated to $\pi$, respectively. See \citet[Chapter 2]{dP(14)} for definition and properties of these differential operators. In particular, here it is enough to recall the following integration-by-parts formula 
that links these operator together:
$$
\int_{\ps} \left\{- \mathfrak{L}_{\pi}[\phi](\theta) \right\} \psi(\theta) \pi(\ud\theta) = \int_{\ps} \langle \Dcal_{\theta}[\phi] \cdot \Dcal_{\theta}[\psi] \rangle \pi(\ud\theta)
$$
for arbitrary $\mathrm{C}^2$ functions $\phi, \psi : \ps \to \R$. Then, we shall say that $W : \ps \to \R$ is a Lyapunov function if $W$ belongs to $\mathrm{C}^2(\ps)$, $W(\theta) \geq 1$ and 
\begin{equation} \label{Lyapunov}
\mathfrak{L}_{V,\pi}[W](\theta) \leq - aW(\theta) + b\ind_{B_R}
\end{equation}
hold for all $\theta \in \ps$, for some suitable constants $a > 0$, $b \geq 0$ and $R > 0$, with $B_R := \{\|\theta\| < R\}$. We notice that, for any function $f \in \mathrm{C}^1_b(\ps)$, we have
\begin{align*}
& \int_{\ps} \frac{-\mathfrak{L}_{V,\pi}[W](\theta)}{W(\theta)} [f(\theta)]^2 e^{-V(\theta)}\pi(\ud\theta) \\
&= \int_{\ps} \frac{-\mathfrak{L}_{\pi}[W](\theta)}{W(\theta)} [f(\theta)]^2 e^{-V(\theta)}\pi(\ud\theta) +
\int_{\ps} \Dcal_{\theta}[V] \cdot \Dcal_{\theta}[W] \frac{[f(\theta)]^2}{W(\theta)} e^{-V(\theta)}\pi(\ud\theta) \\
&= \int_{\ps} \Dcal_{\theta}[W] \cdot \Dcal_{\theta}\left[\frac{[f(\theta)]^2}{W(\theta)} e^{-V(\theta)}\right]\pi(\ud\theta) +
\int_{\ps} \Dcal_{\theta}[V] \cdot \Dcal_{\theta}[W] \frac{[f(\theta)]^2}{W(\theta)} e^{-V(\theta)}\pi(\ud\theta) \\
&= \int_{\ps} \Dcal_{\theta}[W] \cdot \Dcal_{\theta} \left[ \frac{[f(\theta)]^2}{W(\theta)}\right] e^{-V(\theta)} \pi(\ud\theta) \\
&= 2\int_{\ps} \frac{f(\theta)}{W(\theta)} \Dcal_{\theta}[W] \cdot \Dcal_{\theta}[f] e^{-V(\theta)} \pi(\ud\theta) - \int_{\ps} \left[\frac{f(\theta)}{W(\theta)}\right]^2 \|\Dcal_{\theta}[W]\|^2 e^{-V(\theta)} \pi(\ud\theta) \\
&= \int_{\ps} \|\Dcal_{\theta}[f]\|^2 e^{-V(\theta)} \pi(\ud\theta) - \int_{\ps}  \left\| \Dcal_{\theta}[f] - \frac{f(\theta)}{W(\theta)}\Dcal_{\theta}[W]\right\|^2 e^{-V(\theta)} \pi(\ud\theta) \\
&\leq \int_{\ps} \|\Dcal_{\theta}[f]\|^2 e^{-V(\theta)} \pi(\ud\theta)\ .
\end{align*}
At this stage, we can follow the same exact steps in \citet{BBCG(08)} to conclude that 
\begin{equation} \label{Teo14francesi}
[ \mathfrak C_2^{(M)}(\mu_{V,\pi})]^2 \leq \frac 1a(1 + b \kappa_R)  
\end{equation}
where the constants $a,b$ are the same as in \eqref{Lyapunov}, while $\kappa_R$ denotes the weighted Poincar\'e-Wirtinger constant (relative to the Malliavin derivative) of the measure $\mu_{V,\pi}$ restricted on the ball $\|\theta\| < R$.

After these preliminaries, let us consider point {\it (1)}. We put $V:=nG$ in \eqref{Gibbs_Vpi}. Let $W$ be a $C^2(\ps)$ function such that $W\ge1$ on $\ps$ and such that $W(\theta)=e^{\|\theta\|}$ if $|\theta|\ge R$. 
Let $C_R:=\sup_{B_r} [|W| + \|\Dcal W\| + |\mathfrak{L}_{\pi}[W]|]$. The above operator $\mathfrak{L}_{V,\pi}$ now becomes $\mathfrak{L}_{\pi} - n\Dcal_{\theta}[G] \cdot \Dcal_{\theta}$.
A computation shows that if $\|\theta\|\ge R$ there holds
\begin{equation} \label{mainini1}
\mathfrak{L}_{V,\pi}[W](\theta) = \left[\sum_{k=1}^{\infty} \lambda_k \left(
\frac{\theta_k^2 + \|\theta\| - \theta_k^2/\|\theta\|}{\|\theta\|^2} -n \frac{\theta_k \Gamma_k}{\|\theta\|} \right)\right]\,W(\theta)
\end{equation}
where $\theta_k$ and $\Gamma_k$ denote the $k$-th coordinate of $\theta$ and $\Dcal_{\theta}[G]$, respectively, with respect to the basis $\{\mathbf e_k\}$. 
Let $\tau_n:= cn - \mathrm{Tr}[Q](1 + 1/R)$, so that $\tau_n>0$ as soon as $n>\mathrm{Tr}[Q](1 + 1/R)/c$. If $|\theta|\ge R$, a combination of \eqref{mainini1} with the assumption 
$\theta\cdot\Dcal_{\theta} G\ge c\|\theta\|_{\ps}$ yields $\mathfrak{L}_{V,\pi}[W](\theta) \le -\tau_n W(\theta)$. If $|\theta|\le R$, we  estimate as
\[
\begin{aligned}
\mathfrak{L}_{V,\pi}[W](\theta) &=-\tau_n W(\theta)+\tau_n W(\theta)+\mathfrak{L}_{V,\pi}[W](\theta)\\
&\le -\tau_nW(\theta)+\tau_n W(\theta) + |\mathfrak{L}_{\pi}[W]| + n\|\Dcal_{\theta} W\| \cdot \|\Dcal_{\theta} G\|\\
&\le -\tau_n W(\theta) + C_R(1+\tau_n+ nG_R).
\end{aligned}
\]
Thus, we have $W(\theta)\ge 1$ and \eqref{Lyapunov} holds for every $\theta\in\ps$, with $a=\tau_n$ and $b = C_R(1+\tau_n+nG_R)$. At this stage, the conclusion follows from \eqref{Teo14francesi}
after noticing that 
\begin{equation} \label{BakryEmery_Factor}
\kappa_R \lesssim \max_{k \in \N} \left\{\frac{\lambda_k}{n\lambda_k \eta_k + 1}\right\}
\end{equation}
which is valid in view of a combination of the assumption \eqref{factor_Poincare_Infinite} with the classical Bakry-Emery criterion and the well-known tensorization property of the Poincar\'e-Wirtinger constant 
(see \cite[Proposition 4.3.1]{BGL(14)}. This complete the proof of point {\it (1)}.


Let us consider point {\it (2)}. Let $W(\theta)=\exp\{G(\theta)-\inf_{\ps} G\}$, $\theta\in\ps$. By direct computation, $\mathfrak{L}_{\pi}[W] = W\{\|\Dcal G \|^2+ \mathfrak{L}_{\pi}[G]\}$. Putting $V:=nG$ in \eqref{Gibbs_Vpi} as above, we have again that the operator $\mathfrak{L}_{V,\pi}$ becomes $\mathfrak{L}_{\pi} - n\Dcal_{\theta}[G] \cdot \Dcal_{\theta}$. Whence, 
\begin{equation}\label{computa}
\mathfrak{L}_{V,\pi}[W] =\left\{ (1-n) \| \Dcal_{\theta}[G] \|^2+ \mathfrak{L}_{\pi}[G]\right\}W \ .
\end{equation}
Thanks to assumption \eqref{growth_francesi}, we have 
$$
(n-1) \| \Dcal_{\theta}[G] \|^2 -  \mathfrak{L}_{\pi}[G] \ge 2(n-1)c_1 + [(n-1)c_2 - 1] \left(\mathfrak{L}_{\pi}[G]\right)_+ \ge nc_1
$$ 
whenever $|\theta|\ge R$ and $n>1+1/c_2$. Thus, if $n>1+1/c_2$, \eqref{computa} entails 
$$
\mathfrak{L}_{V,\pi}[W] \le - c_1n\,W
$$ 
whenever $|\theta|\ge R$. On the other hand, if $|\theta|\le R$, we easily deduce from \eqref{computa} that
$$
\mathfrak{L}_{V,\pi}[W] \le -c_1n\,W(\theta)+ e^{\omega_R}(c_1n +G_R^*).
$$
Thus, we have $W(\theta)\ge 1$ and $\mathfrak{L}_{V,\pi}[W](\theta) \le -c_1n\,W(\theta)+\tilde b_n\chi_{B_R}(\theta)$ for every $\theta\in\ps$, where $\tilde b_n:= e^{\omega_R}(c_1n+G_R^*)$. 
to conclude, we resort to \eqref{Teo14francesi}, which holds with $a=c_1n$ and $b = \tilde b_n$, in combination with \eqref{BakryEmery_Factor}.


\subsection{Proof of Theorem \ref{main_thm1}} \label{proof:thm1}
To establish \eqref{Bayes_measure}, we start from the Bayes formula
$$
\pi_n(\ud\theta|x^{(n)}) = \frac{\left[\prod_{i=1}^n f(x_i|\theta)\right]}{\int_{\ps} \left[\prod_{i=1}^n f(x_i|\tau)\right] \pi(\ud\tau)} \pi(\ud\theta)
$$
and we observe that the regularity of the mapping $x \mapsto f(x|\theta)$ allows us to write $\prod_{i=1}^n f(x_i|\theta)$ as 
$$
\exp\left\{\sum_{i=1}^n \log f(x_i|\theta)\right\} = \exp\left\{ n\int_{\ss} \log f(y|\theta) \left(\frac 1n \sum_{i=1}^n \delta_{x_i}(\ud y) \right)\right\}.
$$
Then, the bound \eqref{buond_statmodel} entails that the integral $\int_{\ss} \log f(y|\theta)\gamma(\ud y)$ is well-defined and finite for any $\gamma \in \mathcal P_2(\ss)$ and $\theta \in \ps$. Now, recalling that $\mu_0 \in \mathcal P_2(\ss)$, 
let $V_0^{(n)}$ be the $\Wdue^{(\pms)}$-neighborhood of $\mu_0$ for which \eqref{Ln} is in force. Let $\zeta$ be a fixed element of such a neighborhood and let 
$\{\zeta_t\}_{t \in [0,1]}$ be a $\mathcal W_2$-constant speed geodesic connecting $\mu_0$ with $\zeta$. In particular,  $[0,1]\ni t\mapsto \zeta_t$ is an absolutely continuous curve in  $\mathcal P_2(\ss)$.
The map $\pi_n^{\ast}$ allows the construction of a lifting of this path, in the sense that $\{\pi_n^{\ast}(\cdot| \zeta_t)\}_{t \in [0,1]}$ is a path in $\mathcal{P}_2(\ps)$ connecting 
$\pi_n^{\ast}(\cdot|\mu_0)$ with $\pi_n^{\ast}(\cdot|\zeta)$, with $\pi_n^{\ast}(\cdot| \zeta_t)$ having full support in $\ps$ for any $t \in [0,1]$. The Benamou-Brenier formula discussed in Section \ref{sec_mainres} 
shows that
\begin{equation} \label{BenBre}
\left[ \Wdue^{(\mathcal{P}_2(\ps))}(\pi_n^{\ast}(\cdot|\mu_0), \pi_n^{\ast}(\cdot|\zeta))\right]^2 \leq \int_0^1 \int_{\ps} \|\mathrm D_{\theta} u(\theta, t)\|^2 \pi_n^{\ast}(\ud\theta|\zeta_t) \ud t
\end{equation}
where, for almost every $t \in (0,1)$, $u(\cdot, t)$ denotes the solution of the elliptic problem. The weak formulation of the elliptic problem reads as the following problem
\begin{equation} \label{Elliptic}
 \int_{\ps} \langle\mathrm D_{\theta}\ u(\theta, t), \mathrm D_{\theta}\ \psi(\theta)\rangle \pi_n^{\ast}(\ud\theta|\zeta_t)=\frac{\ud}{\ud s} \int_{\ps} \psi(\theta) \pi_n^{\ast}(\ud\theta\ |\ \zeta_s)\ _{\big| s= t}\ , 
\quad\forall\ \psi \in C^1_b(\ps). 
\end{equation}
The right space for the solution of this problem is, for fixed $t \in (0,1)$, the weighted Sobolev space $\mathrm{H}^1_m(\ps; \pi_n^{\ast}(\cdot|\zeta_t))$, defined as the completion of the space
$$
\left\{\psi \in C^1_b(\overline{\ps})\ \Big|\ \int_{\ps} \psi(\theta) \pi_n^{\ast}(\ud\theta|\zeta_t) = 0 \right\}
$$
with respect to the norm 
$$
\|\psi\|_{1,t} := \left( \int_{\ps} \|\mathrm D_{\theta} \psi(\theta)\|^2 \pi_n^{\ast}(\ud\theta|\zeta_t) \right)^{1/2}
$$ 
associated with the scalar product
$$
\langle \varphi, \psi \rangle_{1,t} := \int_{\ps} \langle\mathrm D_{\theta}\ \varphi(\theta), \mathrm D_{\theta} \ \psi(\theta)\rangle \pi_n^{\ast}(\ud\theta|\zeta_t)\ .
$$
Now, since the equation displayed in \eqref{Elliptic} can be re-written in the abstract form as
$$
T_t[\psi] = \langle u(\cdot, t), \psi \rangle_{1,t}
$$
where 
\begin{equation} \label{operatorTt}
T_t[\psi] := \frac{\ud}{\ud s} \int_{\ps} \psi(\theta) \pi_n^{\ast}(\ud\theta|\zeta_s)\ _{\big| s= t}\ ,
\end{equation}
existence, uniqueness and regularity for the solution of \eqref{Elliptic} would follow from the Riesz representation theorem, provided that the functional $T_t$ belongs to the dual of 
$\mathrm{H}^1_m(\ps; \pi_n^{\ast}(\cdot|\zeta_t))$. Whence, again by Riesz theorem, we have
$$ 
\|u(\cdot, t)\|_{1,t} = \left( \int_{\ps} \|\mathrm D_{\theta} u(\theta, t)\|^2 \pi_n^{\ast}(\ud\theta|\zeta_t) \right)^{1/2} = \sup_{\substack{\psi \in \mathrm{H}^1_m(\ps; \pi_n^{\ast}(\cdot|\zeta_t))
\\ \|\psi\|_{1,t} \leq 1}} T_t[\psi] \ . 
$$
Accordingly, by combining the Riesz representation with \eqref{BenBre}, we obtain that
\begin{equation} \label{stimaRiesz}
\left[ \Wdue^{(\mathcal{P}_2(\ps))}(\pi_n^{\ast}(\cdot|\mu_0), \pi_n^{\ast}(\cdot|\zeta))\right]^2 \leq 
\int_0^1 \left[ \sup_{\substack{\psi \in \mathrm{H}^1_m(\ps; \pi_n^{\ast}(\cdot|\zeta_t))
\\ \|\psi\|_{1,t} \leq 1}} T_t[\psi] \right]^2 \ud t\ .
\end{equation}
Now, in order to obtain further estimates, we introduce the following function
$$
G(\theta, t) := \int_{\ss} \log f(y|\theta) \zeta_t(\ud y)
$$
and we indicate by $G'(\theta, t)$ the partial derivative of $G(\theta, t)$ with respect to $t$. 
Coming back to the expression of the operator $T_t$, after justifying the exchange of derivatives with integrals by the regularity assumptions on the mapping $(x,\theta) \mapsto f(x|\theta)$,
the Leibniz rule for the derivative of a quotient gives
\begin{align} 
T_t[\psi] &= n \frac{\int_{\ps} \psi(\theta) G' e^{n G} \ud\pi \int_{\ps} e^{n G} \ud\pi - \int_{\ps} \psi(\theta) e^{n G} \ud\pi \int_{\ps} G'
e^{n G} \ud\pi}{\left( \int_{\ps} e^{n G} \ud\pi \right)^2} \nonumber \\
&= n \left[ \int_{\ps} \psi(\theta) G'(\theta, t) \pi_n^{\ast}(\ud\theta|\zeta_t) - \int_{\ps} \psi(\theta) \pi_n^{\ast}(\ud\theta|\zeta_t) 
\int_{\ps} G'(\theta, t) \pi_n^{\ast}(\ud\theta|\zeta_t) \right]\nonumber \\
&= n \int_{\ps} \left[\psi(\theta) - \int_{\ps} \psi(\tau) \pi_n^{\ast}(\ud\tau|\zeta_t)\ \right] 
\left[ G'(\theta, t) - \int_{\ps} G'(\tau, t) \pi_n^{\ast}(\ud\tau|\zeta_t)\ \right] \pi_n^{\ast}(\ud\theta|\zeta_t)\ . \label{operatorTt2}
\end{align}
The last term in the above chain of inequalities can be interpreted as a covariance operator, so that the Cauchy-Schwartz inequality entails the following
\begin{align} 
\{T_t[\psi]\}^2\leq n^2 \int_{\ps}& \left[\psi(\theta) - \int_{\ps} \psi(\tau) \pi_n^{\ast}(\ud\tau|\zeta_t)\ \right]^2 \pi_n^{\ast}(\ud\theta|\zeta_t) \times \nonumber\\
& \times \int_{\ps} \left[ G'(\theta, t) - \int_{\ps} G'(\tau, t) \pi_n^{\ast}(\ud\tau|\zeta_t)\ \right]^2 \pi_n^{\ast}(\ud\theta|\zeta_t) \ . \label{CauchySchwartz}
\end{align}
In order to obtain further bounds, we now recall the definition of the Poincar\'e-Wirtinger constant $\mathfrak C_2[\cdot]$, which is given in Section \ref{sect:main_problem}. Thus, \eqref{CauchySchwartz} directly gives 
\begin{equation} \label{Poincare2}
\{T_t[\psi]\}^2\leq n^2 \{\mathfrak C_2[\pi_n^{\ast}(\cdot|\zeta_t)]\}^4 \int_{\ps} \|\mathrm D_{\theta} \psi(\theta)\|^2 \pi_n^{\ast}(\ud\theta|\zeta_t)
\int_{\ps} \|\mathrm D_{\theta} G'(\theta, t)\|^2 \pi_n^{\ast}(\ud\theta|\zeta_t) \ . 
\end{equation}
Now, we provide another expression for $G'(\theta, t)$, exploiting the fact the $\zeta_t$ is a Wasserstein constant speed geodesic: indeed, applying again the Benamou-Brenier representation, in this case there exist 
$\mathbf w\in L^1((0,1);L^2_{\zeta_t}(\mathbb X))$ such that
\begin{equation} \label{BenBreX}
\left[ \Wdue^{(\pms)}(\mu_0, \zeta) \right]^2 = \int_0^1 \int_{\ss} |\mathbf  w(x, t)|^2 \zeta_t(\ud x) \ud t
\end{equation}
where, for almost every $t \in (0,1)$, $\mathbf w(\cdot, t)$ satisfies 
\begin{equation} \label{EllipticX}
\frac{\ud}{\ud s} \int_{\ss} \phi(x) \zeta_s(\ud x)\ _{\big| s= t} = \int_{\ss} \mathbf w(x, t) \cdot \nabla_x\ \phi(x) \zeta_t(\ud x)
\quad\quad \forall\ \phi \in C^1_b(\ss)\ .
\end{equation}
See \citet[Proposition 3.30]{AG}. At this stage, in view of a density argument and in view of \eqref{buond_statmodel}, by replacing $\phi$ by $\log f(\cdot|\theta)$ in \eqref{EllipticX} yields
\begin{equation} \label{Gprimo}
G'(\theta, t) = \int_{\ss} \mathbf  w(x, t) \cdot \frac{\nabla_x\ f(x|\theta)}{f(x|\theta)}\zeta_t(\ud x).
\end{equation}
Thus, in view of \eqref{Poincare2}, the squared supremum in \eqref{stimaRiesz} can be bounded as follows
\begin{equation} \label{ennio}
\sup_{\substack{\psi \in \mathrm{H}^1_m(\ps; \pi_n^{\ast}(\cdot|\zeta_t))\\ \|\psi\|_{1,t} \leq 1}} \{T_t[\psi]\}^2 = n^2 \{\mathfrak C_2[\pi_n^{\ast}(\cdot|\zeta_t)]\}^4
\int_{\ps} \|\mathrm D_{\theta}\ G'(\theta, t)\|^2 \pi_n^{\ast}(\ud\theta|\zeta_t) \ .
\end{equation}
By \eqref{Gprimo}, after justifying the exchange of the gradient with the integral, we can write
$$
\nabla_{\theta}\ G'(\theta, t) = \int_{\ss} \mathbf w(x, t) \cdot \mathrm{D}_\theta \frac{\nabla_x\ f(x|\theta)}{f(x|\theta)}\zeta_t(\ud x)
$$
so that, again by Cauchy-Schwartz, 
$$
\|\mathrm{D}_{\theta}\ G'(\theta, t)\|^2 \leq \int_{\ss} |\mathbf w(x, t)|^2 \zeta_t(\ud x) \int_{\ss}  \Big\|\mathrm{D}_\theta \frac{\nabla_x\ f(x|\theta)}{f(x\ |\ \theta)}\Big\|^2 \zeta_t(\ud x)\ .
$$
Then, by a direct combination of \eqref{stimaRiesz} and \eqref{ennio} with this last inequality we obtain
\begin{align} 
\left[ \Wdue^{(\mathcal{P}_2(\ps))}(\pi_n^{\ast}(\cdot|\mu_0), \pi_n^{\ast}(\cdot|\zeta))\right]^2 &\leq n^2 \int_0^1 
\{\mathfrak C_2[\pi_n^{\ast}(\cdot|\zeta_t)]\}^4 \left( \int_{\ss} |\mathbf  w(x, t)|^2 \zeta_t(\ud x) \right) \times \nonumber \\
&\times 
\left(\int_{\ps} \int_{\ss}  \Big\|\mathrm{D}_\theta \frac{\nabla_x\ f(x|\theta)}{f(x|\theta)}\Big\|^2 \zeta_t(\ud x) \pi_n^{\ast}(\ud\theta|\zeta_t) \right) \ud t\ .
\label{final1}
\end{align}
We invoke assumption \eqref{Ln} to conclude that
$$
\left[ \Wdue^{(\mathcal{P}_2(\ps))}(\pi_n^{\ast}(\cdot|\mu_0), \pi_n^{\ast}(\cdot|\zeta))\right]^2 \leq \{L_0^{(n)}\}^2 \int_0^1\int_{\ss} |\mathbf w(x, t)|^2 \zeta_t(\ud x)
$$
which, in view of \eqref{BenBreX}, coincides with \eqref{Lipschitz_kernel} to be proved. To get \eqref{PCR_main1}, we start from considering the right-hand side of the last inequality in \eqref{splitPCR}. For the first summand,
\begin{align*}
&\Wdue^{(\mathcal P(\ps))}(\pi_n^{\ast}(\cdot| \mu_0); \delta_{\theta_0}) \\
&= \left( \frac{\int_{\ps} \|\theta -\theta_0\|^2  \exp\left\{ n\int_{\ss} [\log f(y|\theta)] f(y|\theta_0) \ud y\right\}
\pi(\ud\theta)}{\int_{\ps} \exp\left\{ n\int_{\ss} [\log f(y|\theta)] f(y|\theta_0) \ud y\right\}\pi(\ud\theta)} \right)^{1/2}
\end{align*}
so that it is enough to observe that
\begin{align*}
&\int_{\ss} [\log f(y|\theta)] f(y|\theta_0) \ud y \\
&= \int_{\ss} [\log f(y|\theta_0)] f(y|\theta_0) \ud y + \int_{\ss} \left[ \log \left(\frac{f(y|\theta)}{f(y|\theta_0)}\right) \right] f(y|\theta_0) \ud y \\
&= \mathrm H(\theta_0) - \mathrm K(\theta|\theta_0)\ . 
\end{align*}
Whence,
$$
\Wdue^{(\mathcal P(\ps))}(\pi_n^{\ast}(\cdot| \mu_0); \delta_{\theta_0}) = \left( \frac{\int_{\ps} \|\theta -\theta_0\|^2  e^{n[\mathrm H(\theta_0) - \mathrm K(\theta|\theta_0)]}
\pi(\ud\theta)}{\int_{\ps} e^{n[\mathrm H(\theta_0) - \mathrm K(\theta|\theta_0)]}\pi(\ud\theta)} \right)^{1/2}\ .
$$
Then, the second summand on the right-hand side of \eqref{PCR_main1} is already provided by the second summand on the the right-hand side of the last inequality in \eqref{splitPCR}.
Finally, the last two terms on the the right-hand side of \eqref{PCR_main1} comes from the last summand on the the right-hand side of \eqref{splitPCR}, after noticing that we have
$$
\Wdue^{(\mathcal P(\ps))}(\pi_n^{\ast}(\cdot| \mu_0); \pi_n^{\ast}(\cdot| \empiric^{(\xi)})) \leq \sqrt{2\int_{\ps} \|\theta\|^2 \pi_n^{\ast}(\ud\theta|\mu_0)} + 
\sqrt{2\int_{\ps} \|\theta\|^2 \pi_n^{\ast}(\ud\theta|\empiric^{(\xi)})}\ .
$$
Indeed, the last term on the above right-hand side yields immediately the last term on the right-hand side of \eqref{PCR_main1}. Lastly, we just observe that
$$
\sqrt{2\int_{\ps} \|\theta\|^2 \pi_n^{\ast}(\ud\theta|\mu_0)} \leq 2\left[ \sqrt{\int_{\ps} \|\theta - \theta_0\|^2 \pi_n^{\ast}(\ud\theta|\mu_0)}  +\|\theta_0\| \right]
$$
so that the result follows.

\subsection{Proof of Corollary \ref{coro1}} \label{proof:coro1}
With respect to the first term on the right-hand side of Equation \eqref{PCR_main1}, we exploit the inequality \eqref{bound_Kullback} to obtain that
\begin{align*}
\frac{\int_{\ps} \|\theta-\theta_0\|^2 e^{-n \kksf(\theta|\theta_0)} \pi(\ud\theta)}{\int_{\ps} e^{-n \kksf(\theta|\theta_0)} \pi(\ud\theta)}
&\leq \frac{\int_{\ps} [\kksf(\theta|\theta_0)]^{2/\beta} e^{-n \kksf(\theta|\theta_0)} \pi(\ud\theta)}{\int_{\ps} e^{-n \kksf(\theta|\theta_0)} \pi(\ud\theta)} \\
&= \frac{\int_0^{\infty} z^{2/\beta} e^{-nz} \mu(\ud z)}{\int_{\ps} e^{-n z} \mu(\ud z)} \sim \left(\frac{1}{n}\right)^{2/\beta}
\end{align*}
where the last asymptotic relation comes from a straightforward application of the Laplace method for approximating exponential integrals. See, e.g., \citet[Theorems 41 and Theorem 43]{Bre(94)}. 
The third term is obtained by just inserting the bound borrowed from \citet[Theorem 2.7]{BGV(07)}. For the last term, we start from applying, in combination, H\"older and Lyapunov's inequalities to get
\begin{align*}
&\ee\left[\left(\frac{\int_{\ps} \|\theta\|^2 \left[\prod_{i=1}^n f(\xi_i|\theta) \right] \pi(\ud\theta)}{\int_{\ps} \left[\prod_{i=1}^n f(\xi_i|\theta) \right]\pi(\ud\theta)} \right)^{1/2} 
\mathds{1}\{\empiric^{(\xi)}\not\in V_0^{(n)}\} \right] \\
& \leq \left\{ \ee\left[\int_{\ps} \|\theta\|^r \pi_n(\ud\theta|\xi_1, \dots, \xi_n)\right] \right\}^{1/r} \cdot \left\{ \P[\empiric \not\in V_0^{(n)}] \right\}^{(r-1)/r}\ .
\end{align*}
The conclusion of the proof then follows by using \eqref{bound_moment_posterior}, again by a direct combination with respect to the bound borrowed from \citet[Theorem 2.7]{BGV(07)}. 


\section*{Acknowledgement}

The authors are grateful to an anonymous Referee and the Associate Editor for their constructive
remarks during the review process. E.D. and S.F. are grateful to Ismael Castillo and Matteo Giordano for helpful discussions. E.D. and S.F. received funding from the European Research Council (ERC) under the European Union's Horizon 2020 research and innovation programme under grant agreement No 817257. E.D. and S.F. gratefully acknowledge the financial support from the Italian Ministry of Education, University and Research (MIUR), ``Dipartimenti di Eccellenza" grant 2018-2022. E.D. and E.M. acknowledge support from the MIUR-PRIN Research Project No 2017TEXA3H. E.M. acknowledges support from the INdAM-GNAMPA 2019 Research Project ``Trasporto ottimo per dinamiche con interazione"


\end{document}